\documentclass[12pt]{article}

\usepackage{amsopn,amssymb,amsthm,amsmath,bm}
\usepackage{paralist}
\usepackage{algorithm,algorithmic}
\usepackage{color}
\usepackage{graphicx}
\usepackage{comment}
\usepackage{hyperref}
\hypersetup{
     colorlinks   = true,
     linkcolor    = blue,
     citecolor    = green
}
\usepackage{multirow}
\usepackage{multicol}
\usepackage{float}
\usepackage{subfigure}
\usepackage{enumitem}
\usepackage{booktabs}
\usepackage[style=numeric,natbib=true,backend=bibtex, maxbibnames=9, maxcitenames=2,sortcites=true]{biblatex}

\usepackage{thmtools}
\usepackage{thm-restate}
\usepackage{cleveref}

\bibliography{references}
\usepackage{framed}

\newlength{\defbaselineskip}
\usepackage[a4paper , total={8.5in, 11in}, margin=1in]{geometry}

\usepackage{enumitem}
\setlist[enumerate,1]{leftmargin=*,wide=0em, label = {\bfseries \roman*.}}
\setlist[itemize,1]{leftmargin=*,wide=0em, itemsep=0pt,topsep=0pt}

\newcommand\bbR{\ensuremath{\mathbb{R}}} 

\newcommand{\Abs}[1]{\left |#1\right|}
\newcommand{\norm}[1]{{\left\|#1\right\|}}

\newcommand*\lin[1]{\langle #1\rangle}



\newif\ifwithcomments
\withcommentstrue
\ifwithcomments
\newcommand{\xxx}[1]{{\color{red} (#1)}}
\else
\newcommand{\xxx}[1]{}
\fi

\DeclareMathOperator{\argmin}{argmin}

\def\lmin{\lambda_{\min}}
\newcommand{\vv}[1] {\mathbf{#1}}

\def\s{\vv{s}}

\def\x{\vv{x}}

\def\w{\vv{w}}

\def\a{\vv{a}}

\def\H{\vv{H}}

\begin{document}
\title{\Large Second-Order Optimization for Non-Convex Machine Learning: An Empirical Study}
\author{
	Peng Xu
	\thanks{
		Institute for Computational and Mathematical Engineering,
		Stanford University,
		Email: pengxu@stanford.edu
	}
	\and
	Farbod Roosta-Khorasani
	\thanks{
		School of Mathematics and Physics, University of Queensland, Brisbane, Australia, and 
		International Computer Science Institute, Berkeley, USA,     
		Email: fred.roosta@uq.edu.au
	}
	\and
	Michael W. Mahoney
	\thanks{
	International Computer Science Institute and Department of Statistics, 
	University of California at Berkeley,    
	Email: mmahoney@stat.berkeley.edu
}
}

\date{\today}
\maketitle
	

\begin{abstract}
While first-order optimization methods such as stochastic gradient descent (SGD) are popular in machine learning (ML), they come with well-known deficiencies, including relatively-slow convergence, sensitivity to the settings of hyper-parameters such as learning rate, stagnation at high training errors, and difficulty in escaping flat regions and saddle points.
These issues are particularly acute in highly non-convex settings such as those arising in neural networks.
Motivated by this, there has been recent interest in second-order methods that aim to alleviate these shortcomings by capturing curvature information. 
In this paper, we report detailed empirical evaluations of a class of Newton-type methods, namely sub-sampled variants of trust region (TR) and adaptive regularization with cubics (ARC) algorithms, for non-convex ML problems.
In doing so, we demonstrate that these methods not only can be computationally competitive with hand-tuned SGD with momentum, obtaining comparable or better generalization performance, but also they are highly robust to hyper-parameter settings. Further, in contrast to SGD with momentum, we show that the manner in which these Newton-type methods employ curvature information allows them to seamlessly escape flat regions and saddle points. 
\end{abstract}

\section{Introduction}
\label{sec:intro}

The large-scale nature of many modern ML problems poses computational challenges which have rendered many classical optimization methods (developed for scientific computing and other related areas) inefficient or inapplicable. 
In this light, first-order optimization methods, such as SGD and its variants, have been the workhorse for ML applications due to their simplicity and versatility. 
On the other hand, since they consider only first-order gradient information, these methods come with well known deficiencies. 
These include relatively-slow convergence and sensitivity to hyper-parameter settings  such as learning rate. Furthermore, without dedicated ``baby-sitting'' of these methods, they can stagnate at high training loss and have difficulty in escaping saddle points and flat regions.
These deficiencies are particularly problematic in highly non-convex ML problems such as those that arise in neural network applications.
By incorporating second-order information, i.e., curvature, second-order optimization methods hold the promise to solve these well-known deficiencies.
When implemented na\"{\i}vely, however, second-order methods are clearly not computationally competitive with first order alternatives. 
This, in turn, has unfortunately lead to the conventional wisdom that second-order methods are not appropriate for large-scale ML applications.

Within the scientific computing community, more so than the ML community, it is well-known that not only stochastic Newton-type methods in general, and Gauss-Newton in particular, can be made scalable \cite{rodoas1,rodoas2,doas12}, but more importantly, and unlike first-order methods, they are also very \emph{resilient} to a variety of adversarial effects~\cite{roosta2016sub1,roosta2016sub2,xu2016sub}.
Perhaps the most well known example is their resilience to \emph{ill-conditioning}. 
A subtle, yet potentially more severe, example is that the success of most first-order methods is tightly intertwined with \emph{fine-tunning} (often many) \emph{hyper-parameters}, most importantly the step size~\cite{berahas2017investigation}. 
It is rare that these methods exhibit acceptable performance on first try, and it often takes many trials and errors before one can see reasonable results. In fact, the ``true training time'', which almost always includes the time it takes to appropriately tune these parameters, can be frustratingly long. This sort of brute force hyper-parameter tuning is naturally computationally as well as financially expensive.
In contrast, second-order optimization algorithms involve much less parameter tuning, and are less sensitive to specific choices of their hyper-parameters~\cite{berahas2017investigation}. 
A third example that has received attention due to the popularity of non-convex deep learning problems~\cite{lecun2015deep,goodfellow2016deep} has to do with avoiding (possibly degenerate) \emph{saddle points} and finding a \emph{local minimum}. 
While some first-order algorithms can guarantee convergence to an approximate \emph{second-order critical point}\footnote{where gradient is very small and Hessian is almost positive semi-definite.}, e.g.,~\cite{jin2017escape,ge2015escaping, levy2016power}, the vast majority of them lack such performance guarantees. 
Instead, their convergence can, at best, only be ensured to \emph{first-order critical points}, which include saddle points.  
It has been argued that converging to saddle points can be undesirable for obtaining low generalization errors~\cite{dauphin2014identifying,choromanska2015loss,saxe2013exact,lecun2012efficient}. In addition, important cases have been demonstrated where SGD stagnates at such high training error points~\cite{he2016deep}. 
In contrast, employing the curvature information in the form of Hessian, in addition to the advantages mentioned above, can help with achieving convergence to second-order criticality. 

\vspace{-1mm}
\subsubsection*{Our Objective} 


Here, we aim to provide an empirical evaluation of variants of Newton-type methods for non-convex ML problems, and study whether they can address, in a computationally competitive manner, the main aforementioned challenges associated with first order methods, i.e., relatively-slow convergence, sensitivity to hyper-parameter settings, stagnation, and entrapment near saddle points.
To do so, we exploit recent theoretical developments in~\citet{xuNonconvexTheoretical2017}, and focus on variants of trust-region~\cite{conn2000trust,xuNonconvexTheoretical2017} and (adaptive) cubic regularization~\cite{cubic1981Griewank,cartis2011adaptiveI,xuNonconvexTheoretical2017} algorithms, in which the Hessian is suitably approximated. 
More specifically, in the context of several non-convex machine learning applications, we study the empirical performance of \emph{sub-sampled} versions of these algorithms and set out to paint a more complete picture of their \emph{practical impact}. 
In the process, we highlight the above-mentioned shortcomings of first-order methods and, in their light, assess the \emph{effectiveness of employing curvature information as remedy}.%

\vspace{-1mm}
\subsubsection*{Main Questions}
To accomplish our objective, we set out to answer the following questions, all of which are designed to shed light on various \emph{practical} aspects of this paper's underlying thesis.
\begin{enumerate}[itemsep=-1pt,topsep=-1pt,label = {\bfseries Q.\arabic*}]
    \item \label{question:efficient} ({\it Computational Efficiency})
    Can these sub-sampled Newton-type methods be \emph{computationally efficient} enough to be competitive with hand-tuned SGD with momentum?
	\item \label{question:parameter} ({\it Robustness to Hyper-parameters})
	Does the performance of such Newton-type methods exhibit \emph{robustness} to \emph{hyper-parameter tuning}?
    \item \label{question:saddle} ({\it Escaping Saddle Point})
	Does employing Hessian information help with \emph{avoiding saddle points} and converging to \emph{lower training errors} in highly non-convex problems?
	\item \label{question:generalization} ({\it Generalization Performance})
    Can second-order methods be beneficial for obtaining \emph{low generalization error} in machine learning problems?
    \item \label{question:sampling} ({\it Benefits of Sub-sampling})
       How does \emph{sub-sampling} in general, and various sub-sampling \emph{schemes} in particular, affect the performance of TR and ARC methods? 
    \item \label{question:comparison} ({\it Comparison Among Second-Order Methods})
       How do sub-sampled TR and ARC \emph{compare with other second-order methods} like (sub-sampled) Gauss-Newton (GN) and limited-memory BFGS (L-BFGS) methods?

\end{enumerate}

\vspace{-1mm}
\subsubsection*{Methodology}
To minimize the effect of many confounding factors involved in almost all empirical evaluations, rather than chasing to beat the state-of-the-art performances on benchmark tasks, we focus our attention on two simple, yet illustrative, classes of non-convex ML problems, i.e., \emph{multi-layer perceptron (MLP) networks} and \emph{non-linear least squares (NLS)}. More specifically, to address \ref{question:efficient}, \ref{question:parameter}, \ref{question:saddle}, and \ref{question:generalization}, which are concerned with comparisons among first and second-order methods, we consider training (deep) MLP networks with no additional bells and whistles; and to address \ref{question:sampling} and \ref{question:comparison}, which involve comparisons among various second-order methods, we consider a simpler NLS problem.
These questions, references to their theoretical studies, and sections of this paper involving their empirical treatment with relevant figures, are all gathered in Table~\ref{tab:questions}.

\begin{table*}[htb]
\caption{Fundamental questions related to the underlying thesis of this paper.}
	\label{tab:questions}
	\centering
	\begin{tabular}{|c|c|c|c|}
		\hline
		Question & Theoretical Results & Empirical Results & Figure  \\ \hline
		\ref{question:efficient} & \citep[Theorems 6-8]{xuNonconvexTheoretical2017} & Sections~\ref{sec:example_1layer} and~\ref{sec:example_encoder} & \ref{fig:cifar},\ref{fig:curves},\ref{fig:mnist}\\
		\ref{question:parameter} & \citep[Theorems 1-3 and 6-8]{xuNonconvexTheoretical2017} & Sections~\ref{sec:example_1layer} and~\ref{sec:example_encoder} & \ref{fig:cifar},\ref{fig:curves},\ref{fig:mnist} \\
		\ref{question:saddle} & -- & Sections~\ref{sec:example_1layer} and~\ref{sec:example_encoder} & \ref{fig:cifar},\ref{fig:curves},\ref{fig:mnist} \\
		\ref{question:generalization} & -- & Sections~\ref{sec:example_1layer}, and~\ref{sec:example_encoder} & \ref{fig:cifar},\ref{fig:curves},\ref{fig:mnist}\\
		\ref{question:sampling} & \citep[Lemmas 4 and 5]{xuNonconvexTheoretical2017} & Section~\ref{sec:example_nls} & \ref{fig:blc}\\
		\ref{question:comparison} & - & Section~\ref{sec:example_nls} & \ref{fig:blc}\\
		\hline
	\end{tabular}
\end{table*} 

In Section~\ref{sec:example_deep_learning}, we address~\ref{question:efficient}--\ref{question:generalization}.
%
In particular, in the context of image classification (Section~\ref{sec:example_1layer}) and deep auto-encoder (Section~\ref{sec:example_encoder}), we 
study the efficiency of sub-sampled TR method, which incorporates inexactness in both Hessian and the sub-problem solver, as compared with hand-tuned SGD with momentum. 
%
We answer \cref{question:efficient} by measuring the convergence speed over total computational cost. 
We then treat~\ref{question:parameter} by demonstrating the resiliency/sensitivity of various algorithms with respect to their main hyper-parameter. We do this through multiple simulations of all these examples with several choices of the main hyper-parameter for each algorithm. The main hyper-parameter is understood as the one for which, in practice, there is no ``typical'' value. For SGD with momentum, the \emph{learning rate} is considered as the main parameter, since the momentum parameter is typically set to $ \approx 0.9 $. For trust region, the main hyper-parameter is the \emph{initial trust region}, as there are typical values for other parameters of the algorithm.
%
For~\ref{question:saddle}, we consider various initialization schemes, including those that are close to high-level saddle points. 
We study the behavior of SGD-based methods near such regions and evaluate whether an appropriate use of curvature can indeed help with escaping such undesirable saddle points and making continued progress towards areas with lower training error. 
To address~\ref{question:generalization}, we present test performances for all our experiments and assess the generalization errors obtained by all the methods considered. 

In Section~\ref{sec:example_nls}, we turn our attention to~\ref{question:sampling}--\ref{question:comparison}, where  
we consider the NLS problem arising from binary classification task with least-squares loss. On several real datasets, we then  demonstrate the effect of various sub-sampling strategies on the performance of sub-sampled TR and ARC methods, as compared with other second-order algorithms.


\vspace{-2mm}
\section{Background}
\label{sec:background}
In this section, we give a brief review of the general formulation of the optimization problems considered in our study as well as the Newton-type algorithms in question.

\vspace{-1mm}
\subsection{Non-Convex Finite-Sum Minimization}
\label{sec:finite_sum}
Following many machine learning applications, we consider the ``finite-sum'' optimization problem
\vspace{-1mm}
\begin{equation}\vspace{-1mm}
\label{eq:obj_sum}
\min_{\x \in \bbR^d} F(\x) \triangleq \frac{1}{n}\sum_{i=1}^n f_i(\x) \tag{{\bf P1}},
\end{equation}
where each $ f_{i}: \mathbb{R}^{d} \rightarrow \mathbb{R} $ is a smooth but possibly non-convex function. 
Many machine learning and scientific computing applications involve optimization problems of the form~\eqref{eq:obj_sum} where each $ f_{i}$  is a loss (or misfit) function corresponding to $i^{th}$ observation (or measurement), e.g.,~\cite{roszas, doas12, tibshirani1996regression, friedman2001elements, kulis2012metric, bottou2016optimization, sra2012optimization}. In particular, in machine learning applications, $F$ in~\eqref{eq:obj_sum} corresponds to the \emph{empirical risk}~\cite{shalev2014understanding} and the goal of solving~\eqref{eq:obj_sum} is to obtain a solution with small generalization error, i.e., high predictive accuracy on ``unseen'' data. 

\vspace{-1mm}
\subsection{Main Algorithms: Sub-Sampled TR and ARC}
\label{sec:subsampled_TR_CR}
Arguably, line-search is the most straightforward approach for \emph{globalization} of many Newton-type algorithms. However, near saddle points where the gradient magnitude can be small, traditional line search methods can be very ineffective and in fact produce iterates that can get stuck at a saddle point~\cite{nocedal2006numerical}. Trust region ~\cite{sorensen1982newton,conn2000trust} and cubic regularization methods \cite{cubic1981Griewank,nesterov2006cubic,cartis2011adaptiveI,cartis2011adaptiveII}  are two elegant globalization alternatives that, specially recently, have attracted much attention. 

In large-scale non-convex settings, where the application of exact Hessian can be computationally infeasible, recently,  \citet{xuNonconvexTheoretical2017} theoretically studied the variants of TR and ARC algorithms in which the Hessian is suitably approximated. The details of the resulting TR and ARC algorithms are give in Algorithms \ref{alg:STR_fg} and \ref{alg:SCR_fg}, respectively. Iterations of these algorithms involve the following sub-problems:
\vspace{-1mm}
\begin{subequations}
\begin{align}
&\text{\textit{\underline{TR Sub-Problem:}}~}  \label{eq:STR_subp}\\[-5pt]
&\s_t \approx \underset{\|\s\|\le \Delta_t}{\argmin}~ m_t(\s) \triangleq \lin{\nabla F(\x_{t}), \s} + \frac{1}{2}\lin{\s, \H_t \s}, \nonumber\\[5pt]
&\text{\textit{\underline{ARC Sub-Problem:}}~} \label{eq:SARC_subp}\\[-5pt]
&\s_t \approx \underset{\s \in \mathbb{R}^{d}}{\argmin}~ m_t(\s) \triangleq \lin{\nabla F(\x_{t}), \s} + \frac{1}{2} \lin{\s, \H_t \s} + \frac{\sigma_{t}}{3} \|\s\|^{3}, \nonumber 
\end{align}
where $ \H $ is an approximation of the the exact Hessian.
\end{subequations}


\begin{figure*}[!htb]
\vspace{-5mm}
\begin{minipage}[t]{.45\textwidth}
\begin{algorithm}[H]
	\caption{Sub-Sampled TR}
	\label{alg:STR_fg}
	\begin{algorithmic}[1]
		\STATE {\bf Input:} 
		\begin{itemize}[wide=0em, itemsep=-5pt,topsep=-5pt, label = \textbf{-}]
			\item Starting point: $\x_0$
			\item Initial trust-region radius: $0 < \Delta_{0}  < \infty$
			\item Additional parameters: $\epsilon_{g}, \epsilon_{H}$, $0<\eta_1\le \eta_2\le 1, \gamma_2 \ge \gamma_1 > 1$. 
		\end{itemize}
		\FOR{$ t = 0,1,\ldots $}
		\STATE Set the approximate Hessian $\H_t$, as in~\eqref{eq:subsampled_H} \label{step:STR_step}
		\IF{$\norm{\nabla F(\x_{t})} \le \epsilon_g, ~\lmin(\H_t) \ge -\epsilon_H\;$}  
		\STATE  Return $\x_t$.
		\ENDIF
		\STATE (Approximately) solve {\bf TR Sub-Problem} \eqref{eq:STR_subp}.
		\STATE Set $\rho_t \triangleq \dfrac{F(\x_t) - F(\x_t + \s_t)}{-m_t(\s_t)}$,
		\IF {$\rho_t \ge \eta_2$}
		\STATE $\x_{t+1} = \x_t + \s_t$ \;and\; $\Delta_{t+1} = \gamma_2 \Delta_t$
		\ELSIF {$\rho_t \ge \eta_1$}
		\STATE $\x_{t+1} = \x_t + \s_t$ \;and\; $\Delta_{t+1} = \gamma_1 \Delta_t$
		\ELSE
		\STATE $\x_{t+1} = \x_t$ \;and\; $\Delta_{t+1} = \Delta_t/\gamma_2$
		\ENDIF
		\ENDFOR
		\STATE {\bf Output:} $\x_{t}$
	\end{algorithmic}
\end{algorithm}
\end{minipage}%
\hfill
\begin{minipage}[t]{.45\textwidth}
\begin{algorithm}[H]
    \caption{Sub-Sampled ARC}
	\label{alg:SCR_fg}
	\begin{algorithmic}[1]
		\STATE {\bf Input:} 
		\begin{itemize}[wide=0em, itemsep=-5pt,topsep=-5pt, label = \textbf{-}]
			\item Starting point: $\x_0$
			\item Initial regularization parameter: $0 < \sigma_{0}  < \infty$ 
			\item Additional parameters: $\epsilon_{g}, \epsilon_{H}$, $0<\eta_1\le \eta_2\le 1, \gamma_2 \ge \gamma_1 > 1$.
		\end{itemize}
		\FOR{$ t = 0,1,\ldots $}
		\STATE Set the approximate Hessian $\H_t$, as in~\eqref{eq:subsampled_H} \label{step:STR_step}
		\IF{$\norm{\nabla F(\x_{t})} \le \epsilon_g, ~\lmin(\H_t) \ge -\epsilon_H\;$}  
		\STATE  Return $\x_t$.
		\ENDIF
		\STATE (Approximately) solve {\bf ARC Sub-Problem} \eqref{eq:SARC_subp}
		\STATE Set $\rho_t \triangleq \dfrac{F(\x_t) - F(\x_t + \s_t)}{-m_t(\s_t)}$ 
		\IF {$\rho_t \ge \eta_2$}
		\STATE $\x_{t+1} = \x_t + \s_t$ \;and\; $\sigma_{t+1} = \sigma_t/\gamma_2$
		\ELSIF {$\rho_t \ge \eta_1$}
		\STATE $\x_{t+1} = \x_t + \s_t$ \;and\; $\sigma_{t+1} = \sigma_t/\gamma_1$
		\ELSE
		\STATE $\x_{t+1} = \x_t$ \;and\; $\sigma_{t+1} = \gamma_2\sigma_t$
		\ENDIF
		\ENDFOR
		\STATE {\bf Output:} $\x_{t}$
	\end{algorithmic}
\end{algorithm}
\end{minipage}%
\vspace{-2mm}
\end{figure*}

\noindent 
\textbf{Hessian Sub-Sampling} We consider~\eqref{eq:obj_sum} in large-scale regime where $n , d \gg 1$. In such settings, the mere evaluations of the Hessian and the gradient increase linearly in $ n $. As studied in \citet{xuNonconvexTheoretical2017}, given a sampling distribution $\{p_i\}_{i=1}^n$ over the set of indices $\{1,2,\cdots,n\}$, the sub-sampled Hession has the form 
\begin{equation}\vspace{-1mm}
\H(\x) \triangleq \frac{1}{n |\mathcal{S}|} \sum_{j \in \mathcal{S}} \frac{1}{p_{j}}\nabla^{2} f_{j}(\x),
\label{eq:subsampled_H}
\end{equation}
where $\mathcal S\subset\{1,2,\cdots, n\}$ is a random sample collection. When $\Abs{\mathcal S} \ll n$, sub-sampling can offer significant computational savings; see \citet{roosta2016sub1,roosta2016sub2,xu2016sub,bollapragada2016exact} for examples of studies in convex settings. Recently~\citet{xuNonconvexTheoretical2017} theoretically showed that randomized sub-sampling can also be seamlessly extended to non-convex settings. 

\citet{xuNonconvexTheoretical2017} further showed that, in certain settings, one can construct more ``informative'' distributions over the indices in $\{1,2,\ldots,n\}$, as opposed to oblivious uniform sampling. Indeed, it is typically advantageous to bias the probability distribution towards picking indices corresponding to those $ f_{i} $'s which are more \emph{relevant}, in certain sense, in forming the Hessian. 
One such setting where this is possible is the finite-sum optimization of the form, 
\begin{equation}\vspace{-1mm}
\label{eq:obj_sum_ERM}
\min_{\x\in\bbR^d} F(\x)  \triangleq \frac{1}{n} \sum_{i=1}^n f_i(\a_i^T\x), \tag{{\bf P2}}
\end{equation}
for some given data vectors $ \{\a_{i}\}_{i=1}^{n} \subset \mathbb{R}^{d} $. Problems of the form~\eqref{eq:obj_sum_ERM}, which is a special case of~\eqref{eq:obj_sum}, arise often in many machine learning problems~\cite{shalev2014understanding}, e.g., logistic regression with least squares loss as in Example~\ref{sec:example_nls}. For problems of the form~\eqref{eq:obj_sum_ERM}, one can construct a more informative sampling scheme by considering the sampling distribution as $ p_i = {|f_i''(\a_i^T\x)|\|\a_i\|^2}/( \sum_{j=1}^n |f_j''(\a_j^T\x)|\|\a_j\|^2 ) $.
In~\citet{xuNonconvexTheoretical2017}, it was shown that, in order to obtain similar approximation guarantee, such nonuniform sub-sampling scheme  yields a sample size which can be significantly smaller than that required by oblivious uniform sampling.

\noindent 
\textbf{Inexact Sub-problem Solver} In Algorithms~\ref{alg:STR_fg} and \ref{alg:SCR_fg}, it is imperative that the sub-problems~\eqref{eq:STR_subp} and~\eqref{eq:SARC_subp}, respectively, are solved only approximately. Indeed, in large-scale problems, where the exact solution of sub-problems is a computational bottleneck, this relaxation is crucial; see~\citet{xuNonconvexTheoretical2017} for precise definitions as well as  ways to obtain the approximate solution of the sub-problems~\eqref{eq:STR_subp} and~\eqref{eq:SARC_subp}.

\vspace{-3mm}
\section{Numerical Experiments}
\label{sec:experiments}
We are now ready to empirically evaluate the performance of the Newton-type methods considered in this paper (and studied theoretically in~\citet{xuNonconvexTheoretical2017}) in several settings. In particular, we study the answers to \ref{question:efficient}--\ref{question:comparison} posed at the outset in Section~\ref{sec:intro}  in the context of two simple, yet illustrative, classes of non-convex optimization problems, i.e. \emph{(deep) MLP networks} (Section ~\ref{sec:example_deep_learning}) and \emph{NLS} (Section~\ref{sec:example_nls}). See the supplementary materials for more experiments.


\vspace{-2mm}
\subsection{General Experimental Settings}
\paragraph{Complexity Measure} In all of our experiments, we plot various quantities vs.\ total \emph{number of propagations}~\cite{goodfellow2016deep}, which is equivalent to measuring the number of oracle calls of function, gradient and Hessian-vector product. This is so since comparing algorithms in terms of ``wall-clock'' time can be highly affected by their particular implementation details as well as system specifications. In contrast, counting the number of propagations, as an implementation and system independent unit of complexity, is most appropriate and fair. Specifically, in neural nets, for a given data at the input layer, evaluation of network's output layer involves one forward propagation. Performing one additional backward propagation gives the corresponding gradient. Each of the Hessian-vector products, required to solve the respective sub-problems of the second-order methods, is equivalent to two gradient evaluations, i.e., compared to the gradient, it involves one additional forward and backward propagations~\cite{pearlmutter1994fast}.
Combining the batch size in each iteration, we summarize the number of propagations per iteration for each algorithm in Table \ref{tab:props}.

\begin{table}[!htb]
\caption{Total number of propagations per iteration for various algorithms. ``$ r $'' denotes the number of Hessian-vector products for solving the respective subproblems, and hence can be different for each algorithm. In the ``Full'' cases (e.g., Full TR), $\Abs{\mathcal S} = n$.}
\label{tab:props}
\centering
\begin{tabular}{ccc}
\toprule
\sc Sub-Sampled TR,ARC,GN  &  L-BFGS & SGD \\
\hline
$ 2 \left( n + |\mathcal{S}| r \right)  $  & $ 2 n  $  & $ 2 |\mathcal{S}|  $\\
\bottomrule
\end{tabular}
\end{table} 

In all of the following experiments, to approximately solve the sub-problems \eqref{eq:STR_subp} and \eqref{eq:SARC_subp}, respectively, we use CG-Steihaug method \cite{steihaug1983conjugate,nocedal2006numerical} and the generalized Lanczos method~\cite{cartis2011adaptiveI} with a maximum of 250 Lanczos iterations.

\vspace{-2mm}
\subsection{Hessian-Free Optimization for MLPs}
\label{sec:example_deep_learning}
In this section, we set out to empirically study the main questions~\ref{question:efficient}--\ref{question:generalization} of Section~\ref{sec:intro}. For this, we consider two simple but non-trivial MLP netwroks under various settings, namely (a) 1-hidden layer neural networks for image classification, and (b) deep auto-encoder. 
For the following experiments, the empirical performance of the following algorithms, in light of Questions~\ref{question:efficient}--\ref{question:generalization}, are compared\footnote{We excluded ARC variants for reasons alluded to in Section~\ref{sec:conclusions}.}: 
\begin{enumerate}[leftmargin=*,wide=0em, itemsep=-1pt,topsep=-1pt,label=(\arabic*)]
	\item {\it TR Uniform}: Algorithm~\ref{alg:STR_fg} with uniform sub-sampling,
	\item {\it GN Uniform}: sub-sampled variants of Gauss-Newton method with modifications introduced in~\citet{martens2010deep,martens2012training}, and 
	\item {\it SGD with Momentum}: mini-batch SGD with momentum term~\cite{sutskever2013importance} and fixed step-size.
\end{enumerate}

\paragraph{Parameter settings} For Algorithm \ref{alg:STR_fg}, we set $\eta_1 = 10^{-4}$, $\eta_2 = 0.8$, $\gamma_1 = 1.2$, and $\gamma_2 = 2$ (these are some values typical used in literature). The sample size used for Hessian sub-sampling is set to $ 5\%  $ of the total training set, i.e., $ 0.05 n $. For momentum SGD, the mini-batch size is chosen as $ 0.05 n $ and the momentum parameter is set as $ 0.9 $ (which is typically set in literature). For GN method, we use the same parameter settings as in~\citet{martens2010deep}.

\vspace{-1mm}
\subsubsection{One-Hidden Layer Neural Network}
\label{sec:example_1layer}
\begin{figure}[H]
\centering
	\subfigure[tight][Random Initial.: Training Loss]{
		\includegraphics[width=0.31\textwidth]
		{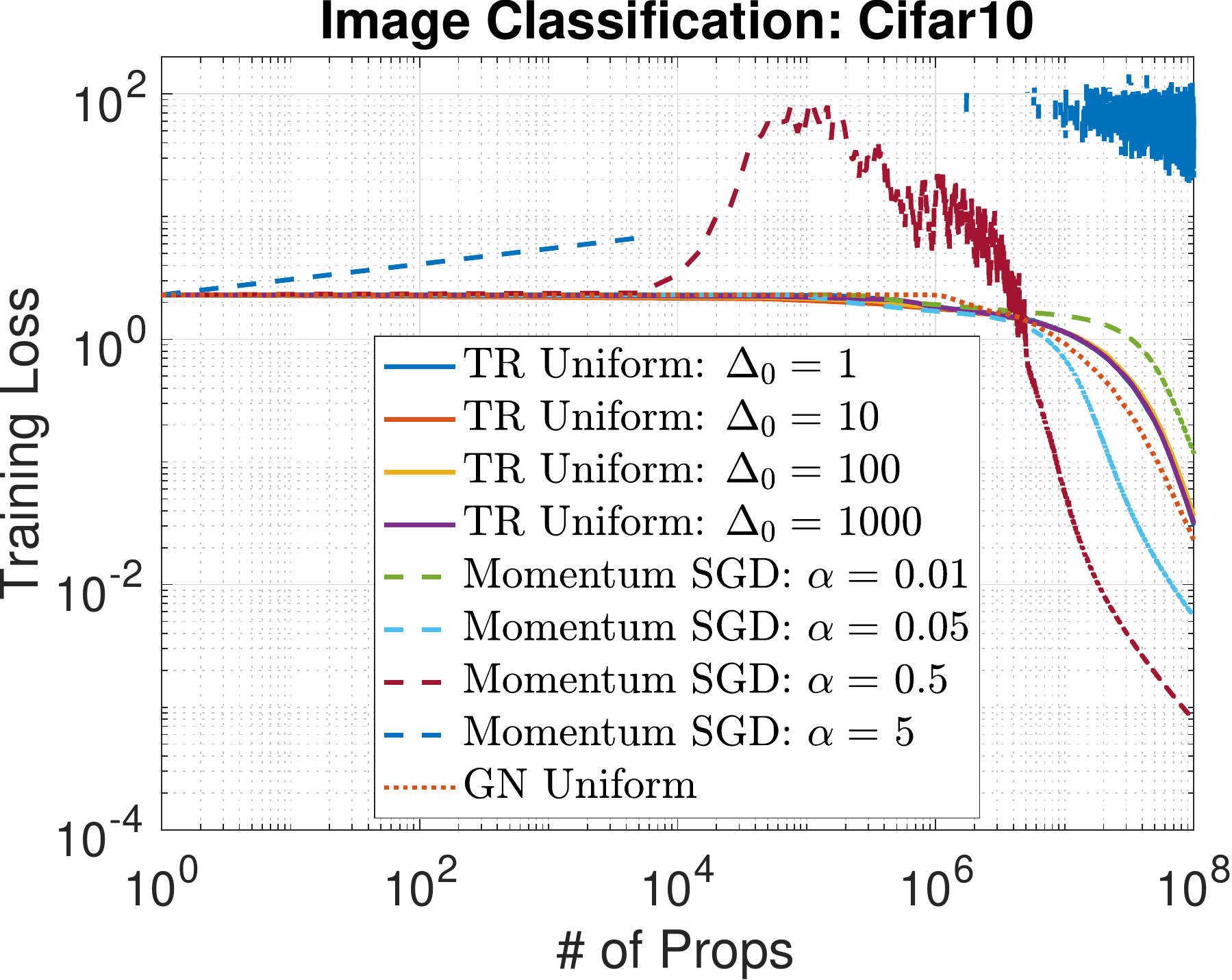}
	}
	\subfigure[Random Initial.: Training Error]{
		\includegraphics[width=0.31\textwidth]
		{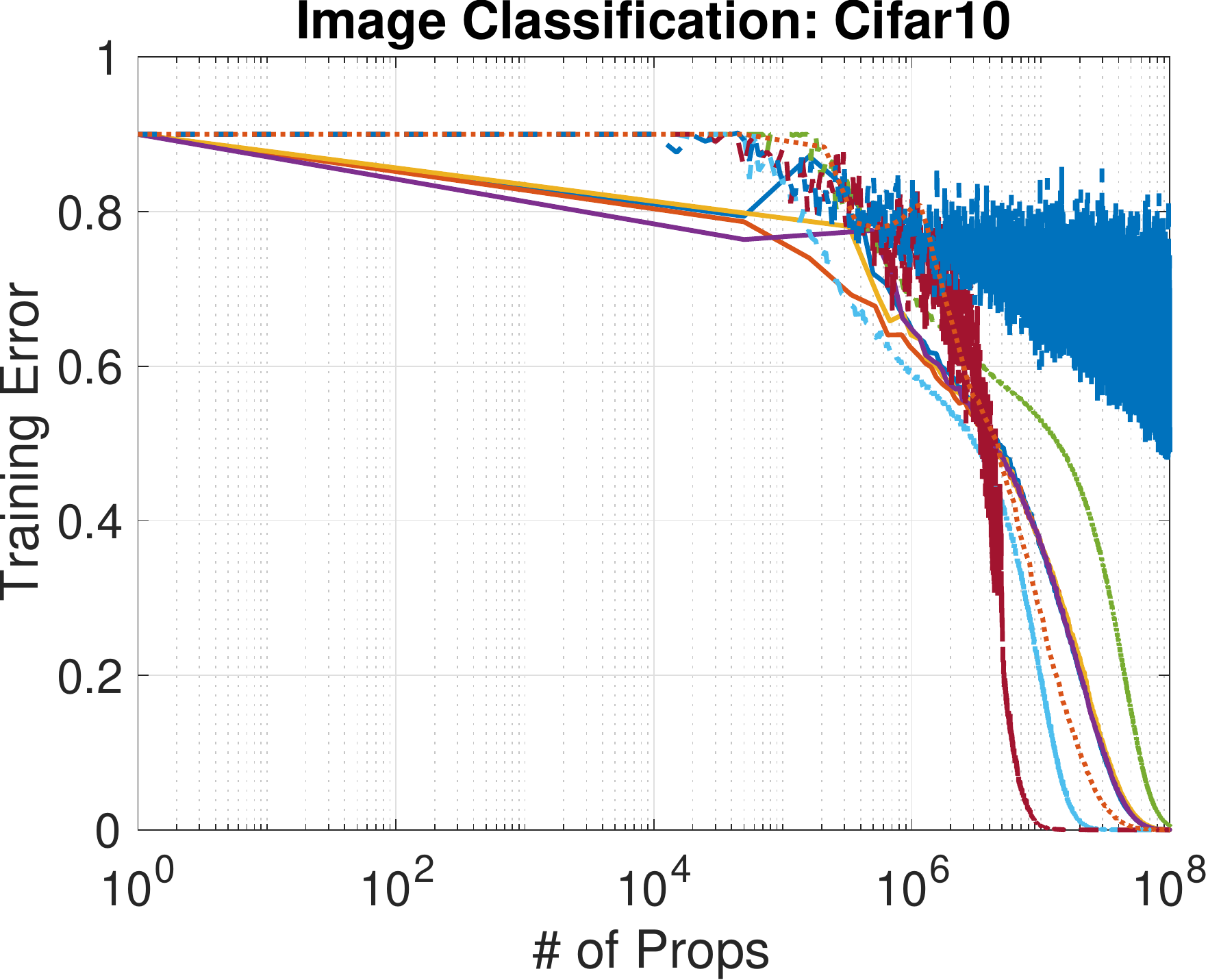}
	}
	\subfigure[Random Initial.: Test Error]{
		\includegraphics[width=0.31\textwidth]
		{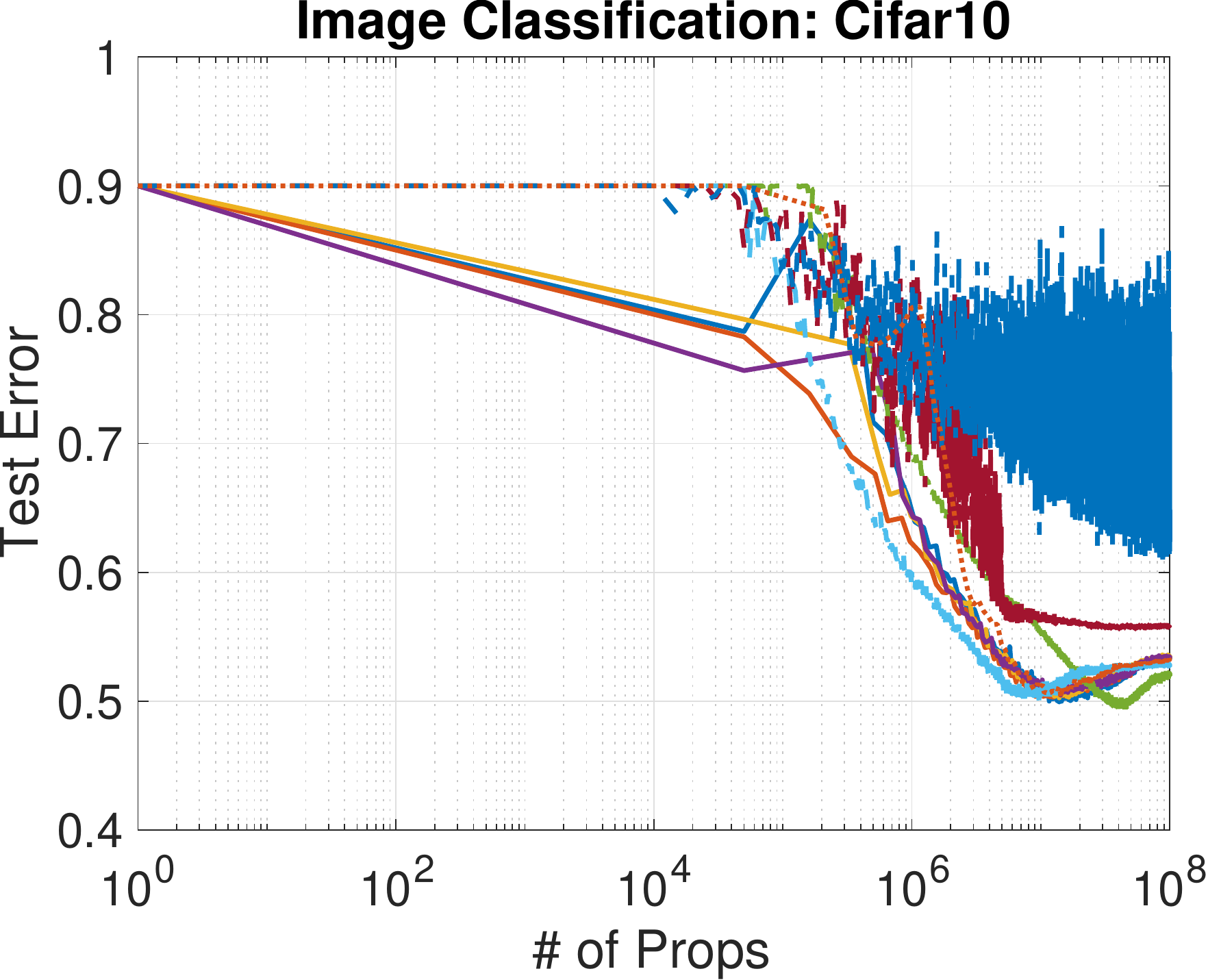}
	}
	\subfigure[Zero Initial.: Training Loss]{
			\includegraphics[width=0.31\textwidth]
			{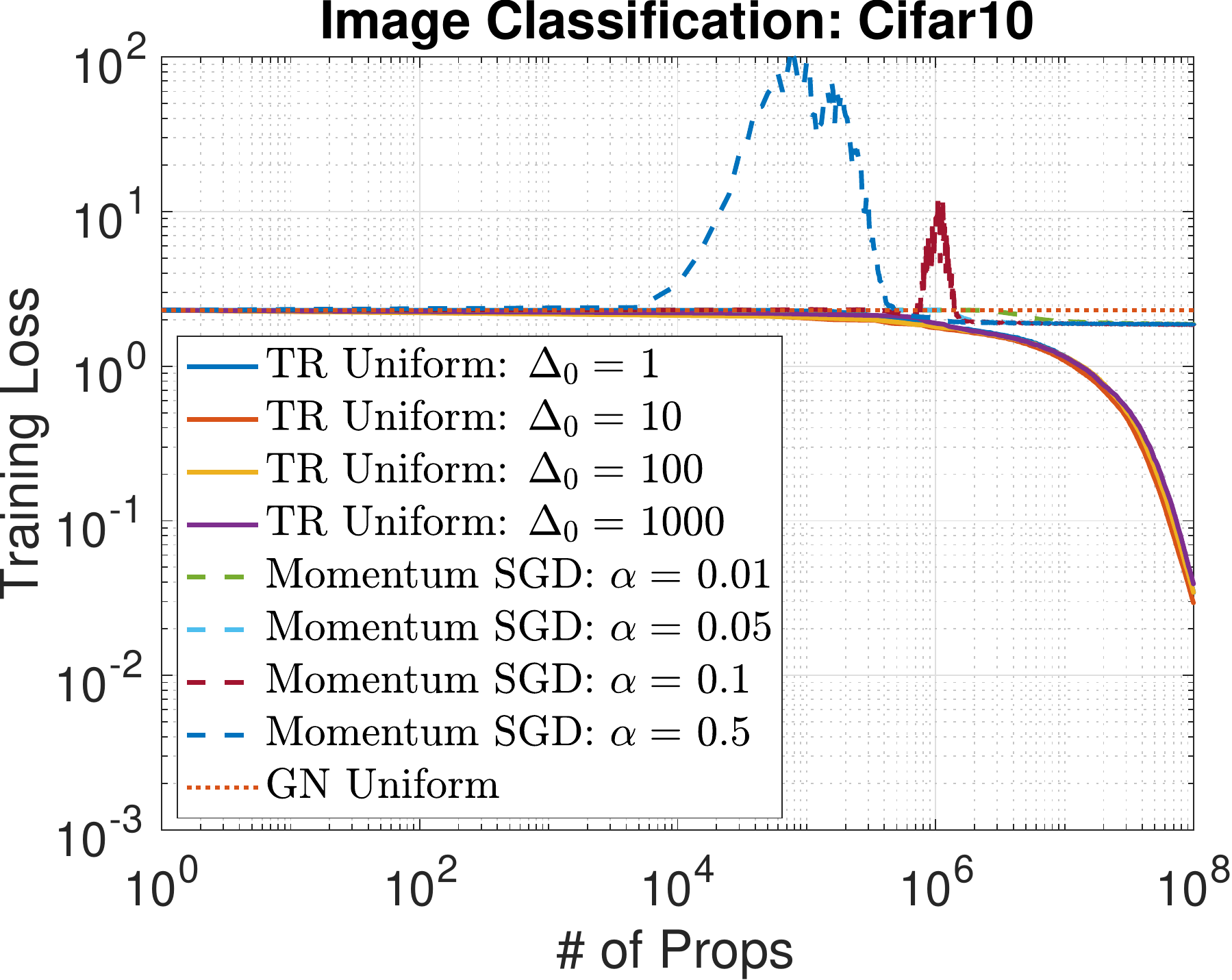}
	}
	\subfigure[Zero Initial.: Training Error]{
			\includegraphics[width=0.31\textwidth]
			{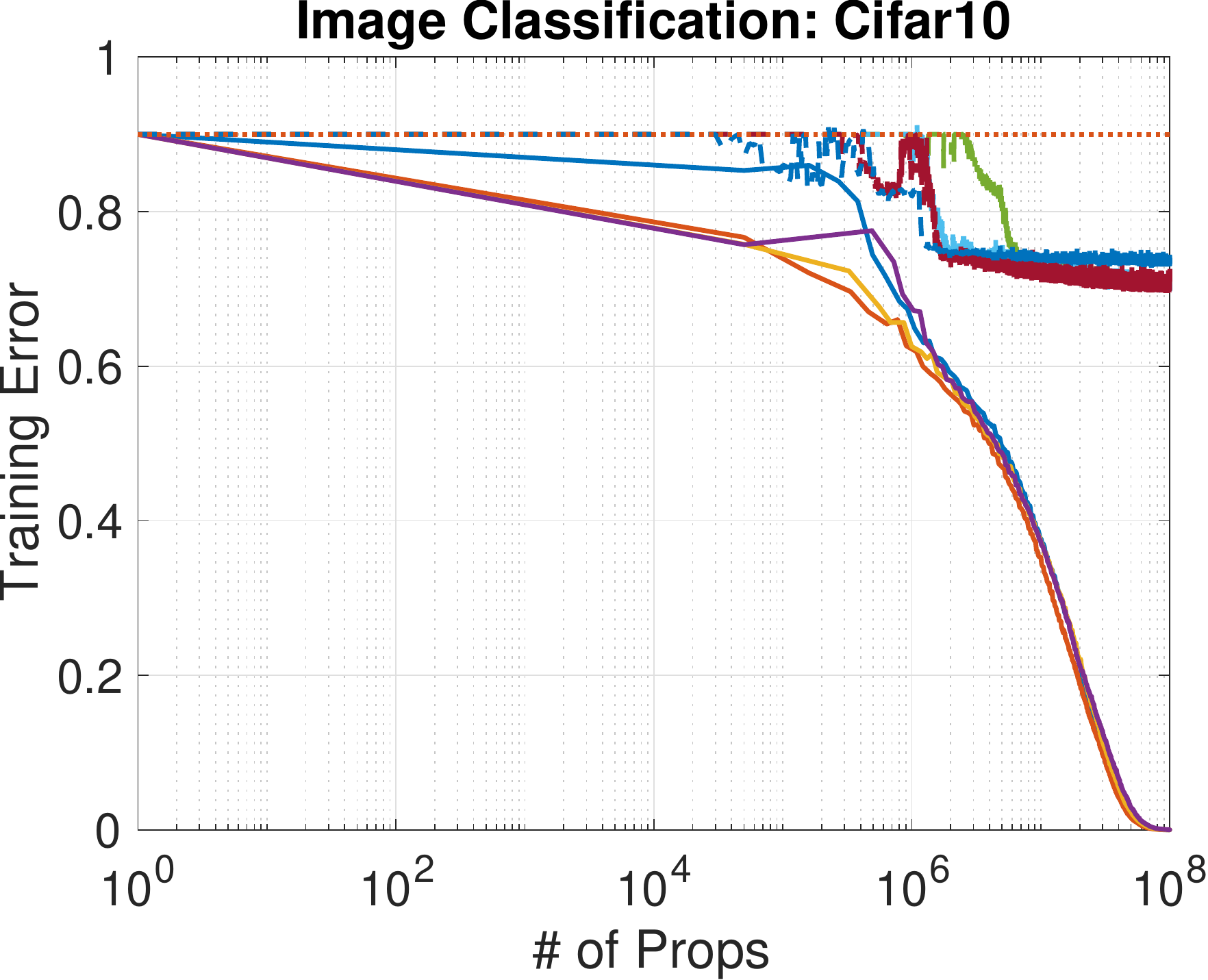}
	}
	\subfigure[Zero Initial.: Test Error]{
			\includegraphics[width=0.31\textwidth]
			{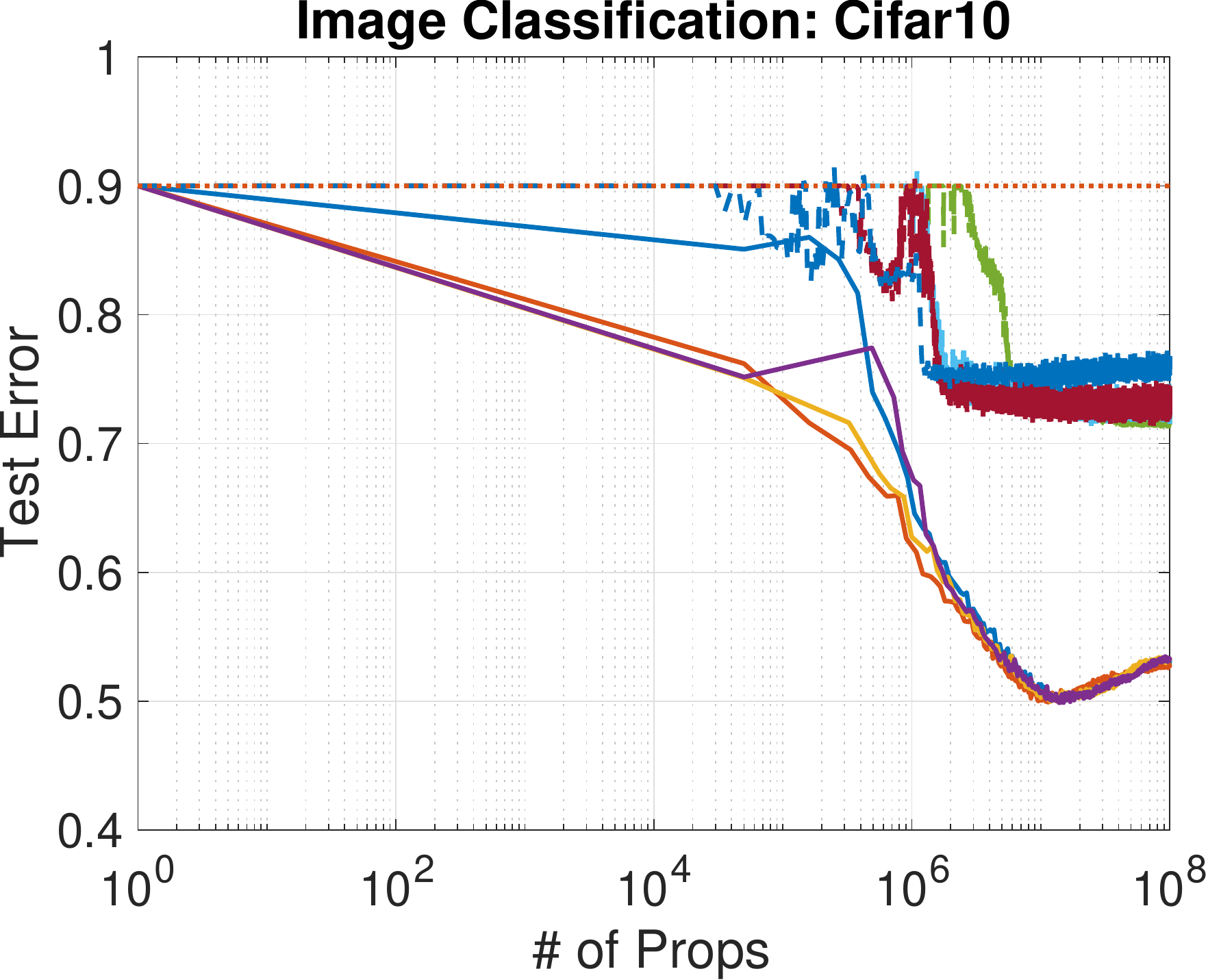}
	}
	\vspace{-3mm}
	\caption{\texttt{Cifar10} dataset on MLP using normalized random (a,b,c) and zero (d,e,f) initialization. $\Delta_{0}$ is the initial trust-region radius of Algorithm~\ref{alg:STR_fg} and $ \alpha $ is the step size for SGD with momentum.}
	\label{fig:cifar}
\end{figure}

\vspace{-1mm}
Here we consider a one-hidden layer neural network for the task of image classification using \texttt{cifar10}~\cite{cifar10} data set. The hidden layer size is $512$, amounting to $ d = 1,578,506 $. 
Initializations are done by setting $ \x_{0} $ to a normalized vector drawn from standard normal distribution and all-zeros vector. 
Figure \ref{fig:cifar} shows training loss, training error and test error for all methods.

In light of \ref{question:efficient}--\ref{question:generalization}, we make the following observations:
\begin{enumerate}[leftmargin=*,wide=0em, itemsep=-1pt,topsep=-1pt,label = {\bfseries Re: Q.\arabic*}]
	\item \textit{(Computational Efficiency)} In Figures \ref{fig:cifar}(a)(b)(c), all algorithms start from a normalized random vector. From Figure \ref{fig:cifar}(a), we observe that in terms of training loss, sub-sampled TR algorithm, though very competitive, can be slightly slower than SGD as long as SGD's step size is appropriately fine-tuned. However, this does not necessarily translate to better test error; see Figure \ref{fig:cifar}(c). In particular, while all TR runs achieve similar test errors, SGD's generalization performance does not mirror its training behavior and appears highly ``chaotic'', i.e., the step-size that achieved the fastest training (red dashed line), generalizes very poorly (compare red and blue dashed  lines).
	
	\item \textit{(Robustness to Hyper-parameters)} From Figures \ref{fig:cifar}(a)(b)(c), we notice the dependence of SGD's performance on the choice of its learning-rate. Well-tuned step size can give both fast convergence of training process and good generalization performance. Too small a step size, however, can lead to slow convergence, while too large a step size, can cause SGD divergence or poor generalization performance. In contrast, Algorithm~\ref{alg:STR_fg} with drastically different initial trust-region radii exhibit comparable performances; similar phenomenon are seen in Figures~\ref{fig:cifar}(d)(e)(f).
	
	\item \textit{(Escaping Saddle Point)} In Figure \ref{fig:cifar}(d)(e)(f), all algorithms start from the origin. We can clearly see that SGD and GN easily get trapped at/near saddle points and/or flat regions (it is easy to check the gradients there are extremely small) and can barely make any progress. In contrast, sub-sampled TR, which effectively utilizes the Hessian, seamlessly escapes these regions and makes continued progress.

	\item \textit{(Generalization Performances)} In all initialization schemes, as shown in Figures~\ref{fig:cifar}, sub-sampled TR obtains competitive, if not better, generalization performance.
\end{enumerate}

\vspace{-1mm}
\subsubsection{Deep Auto-Encoder}
\label{sec:example_encoder}
\vspace{-1mm}

\begin{figure}[!htb]
	\centering
	\includegraphics[width = 0.4\textwidth]{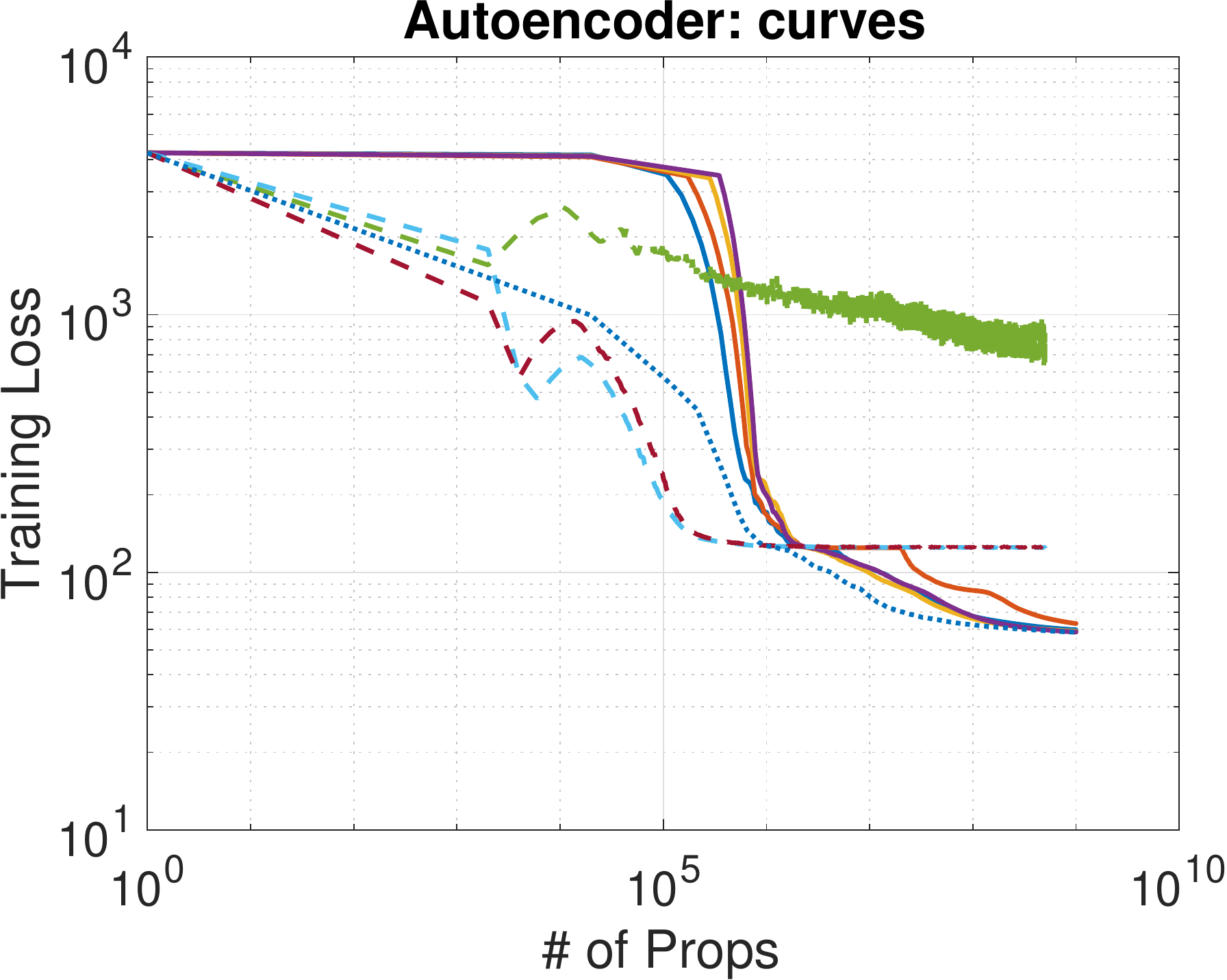}
	\includegraphics[width = 0.4\textwidth]{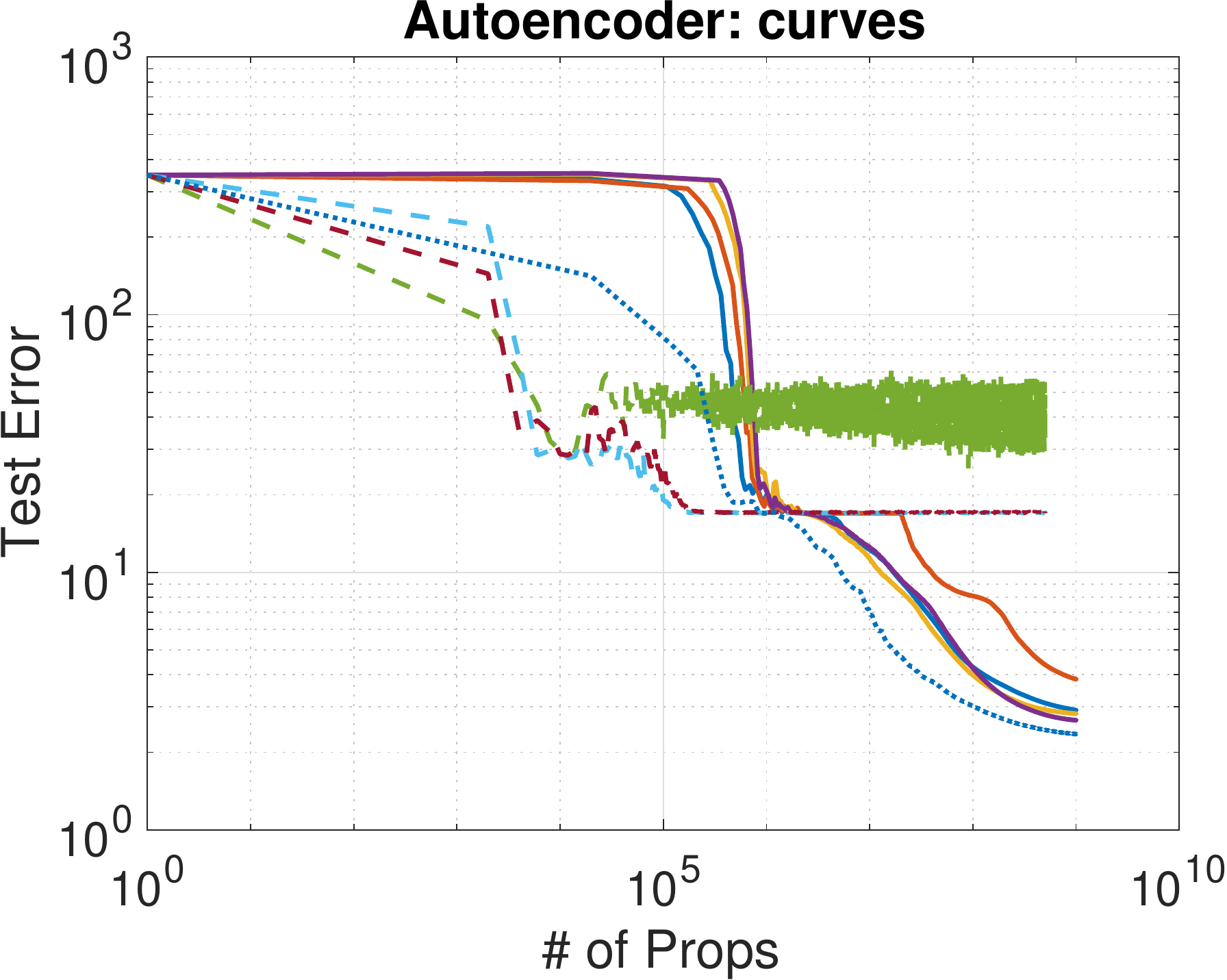}
		\subfigure[Random Initialization]{
			\includegraphics[width = 0.65\textwidth]{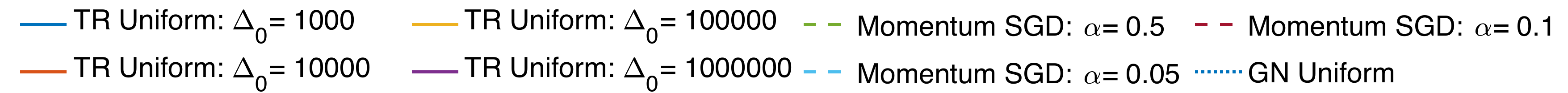}
			
	}
	
	\includegraphics[width = 0.4\textwidth]{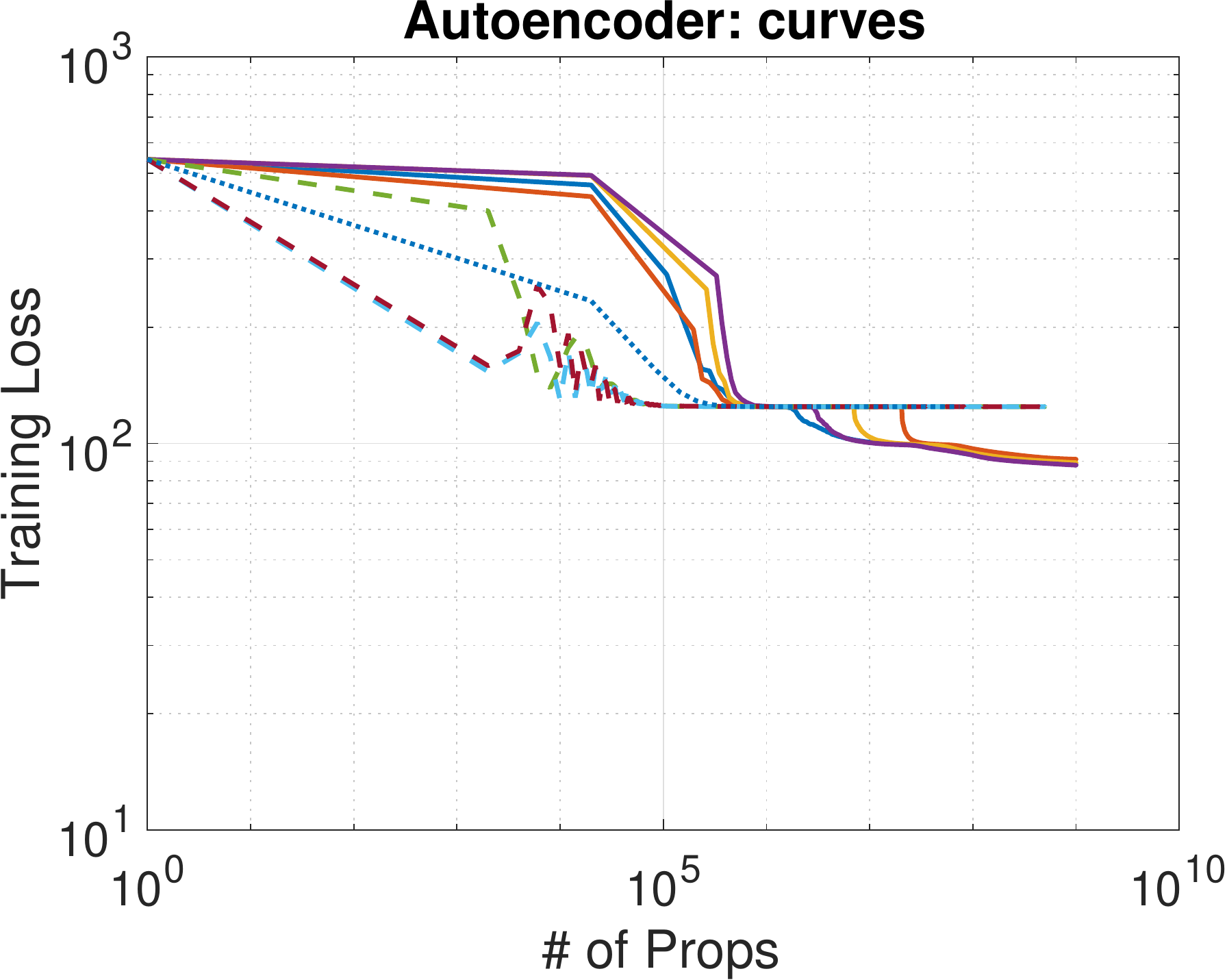}
	\includegraphics[width = 0.4\textwidth]{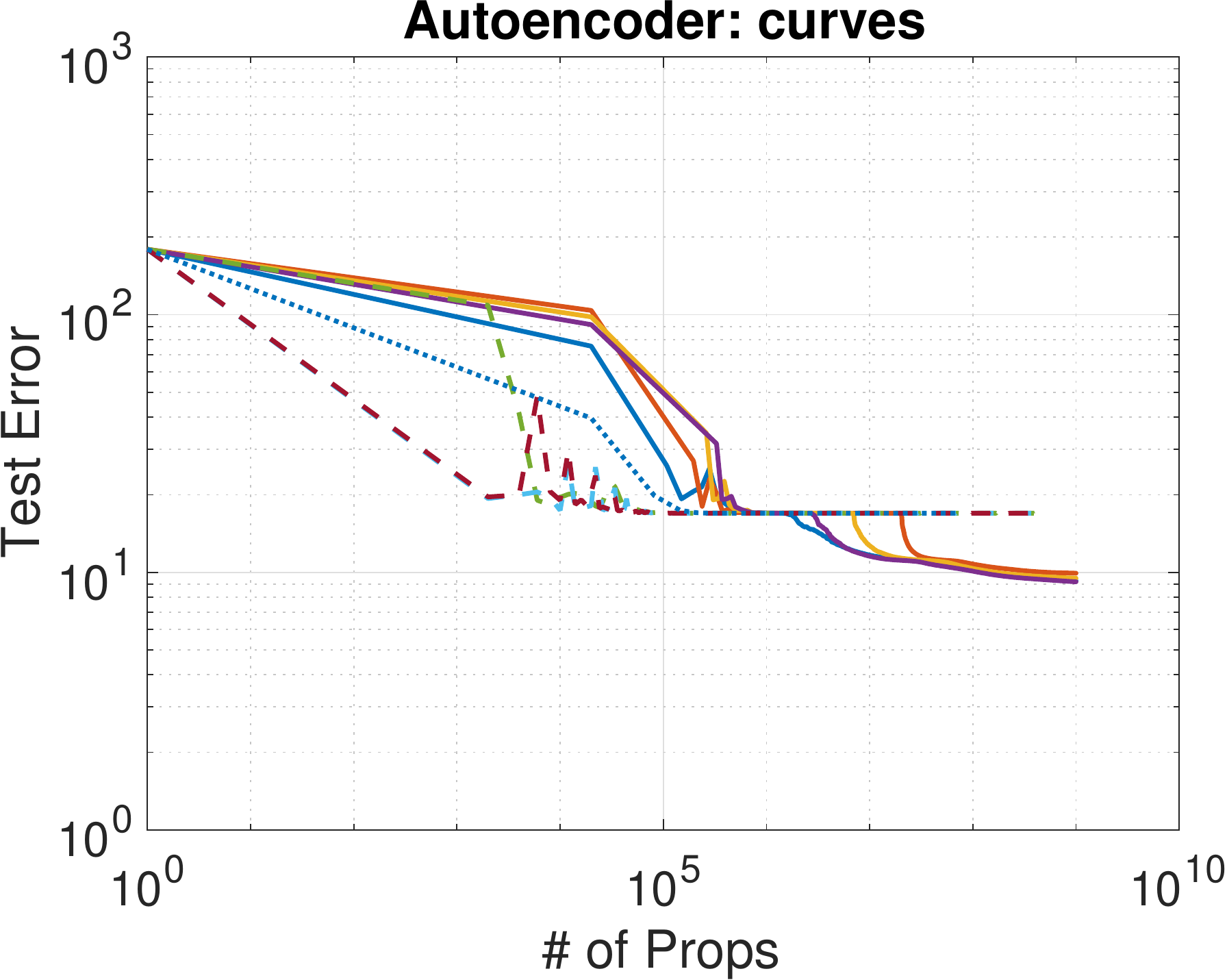}
		\subfigure[Zero Initialization]{
			\includegraphics[width = 0.65\textwidth]{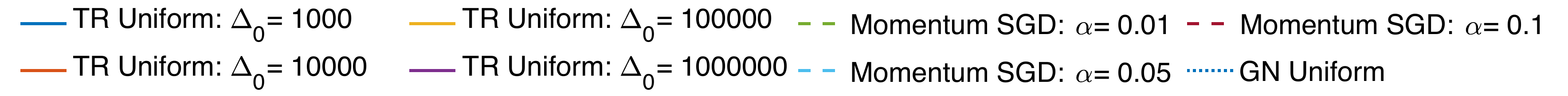}
	}
	\vspace{-3mm}
	\caption{Deep autoencoder on \texttt{curves} dataset using random (a) and zero (b) initial point.  $\Delta_{0}$ is the initial trust-region radius of Algorithm~\ref{alg:STR_fg} and $ \alpha $ is the step size for SGD with momentum.}
	\label{fig:curves}
\end{figure}

\begin{figure}[!htb]
	\centering
	        \includegraphics[width = 0.4\textwidth]{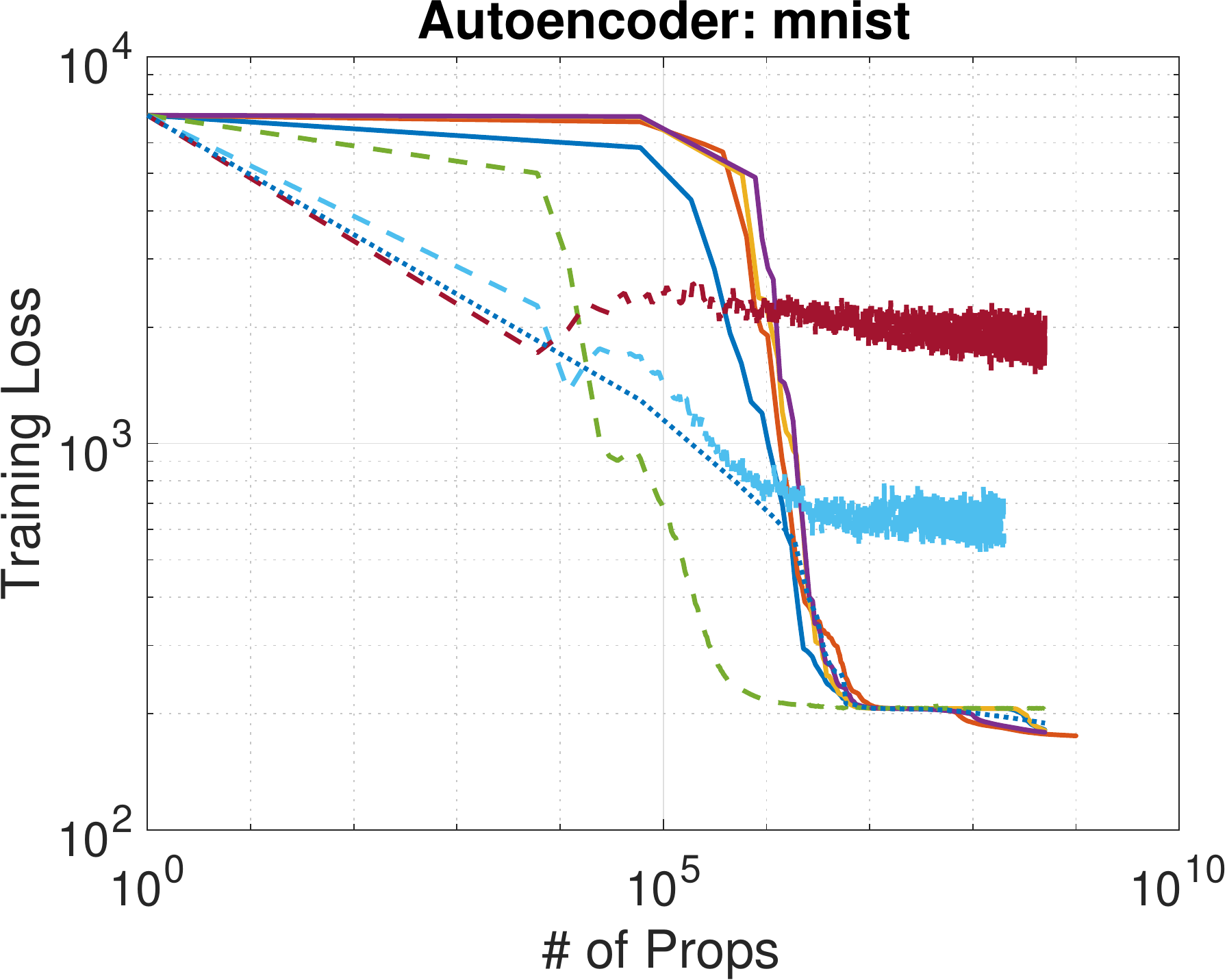}
			\includegraphics[width = 0.4\textwidth]{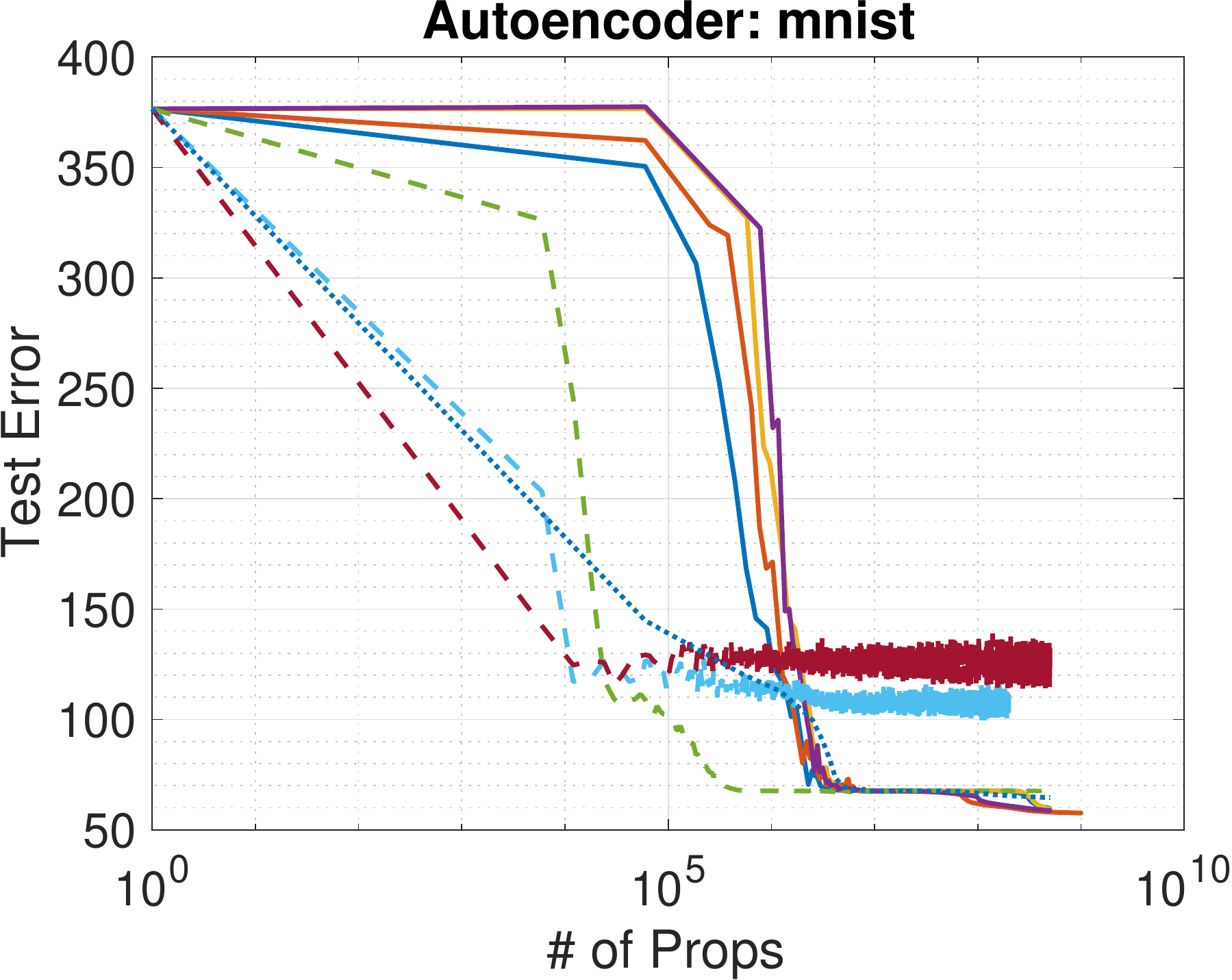}
		\subfigure[Random Initialization]{

			\includegraphics[width = 0.65\textwidth]{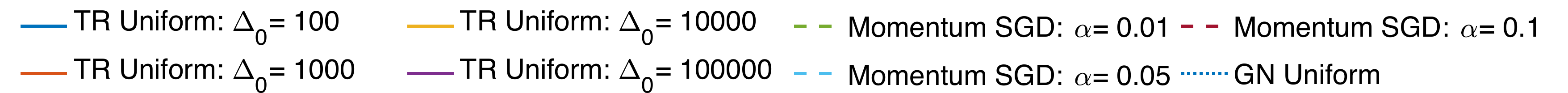}
	}
			\includegraphics[width = 0.4\textwidth]{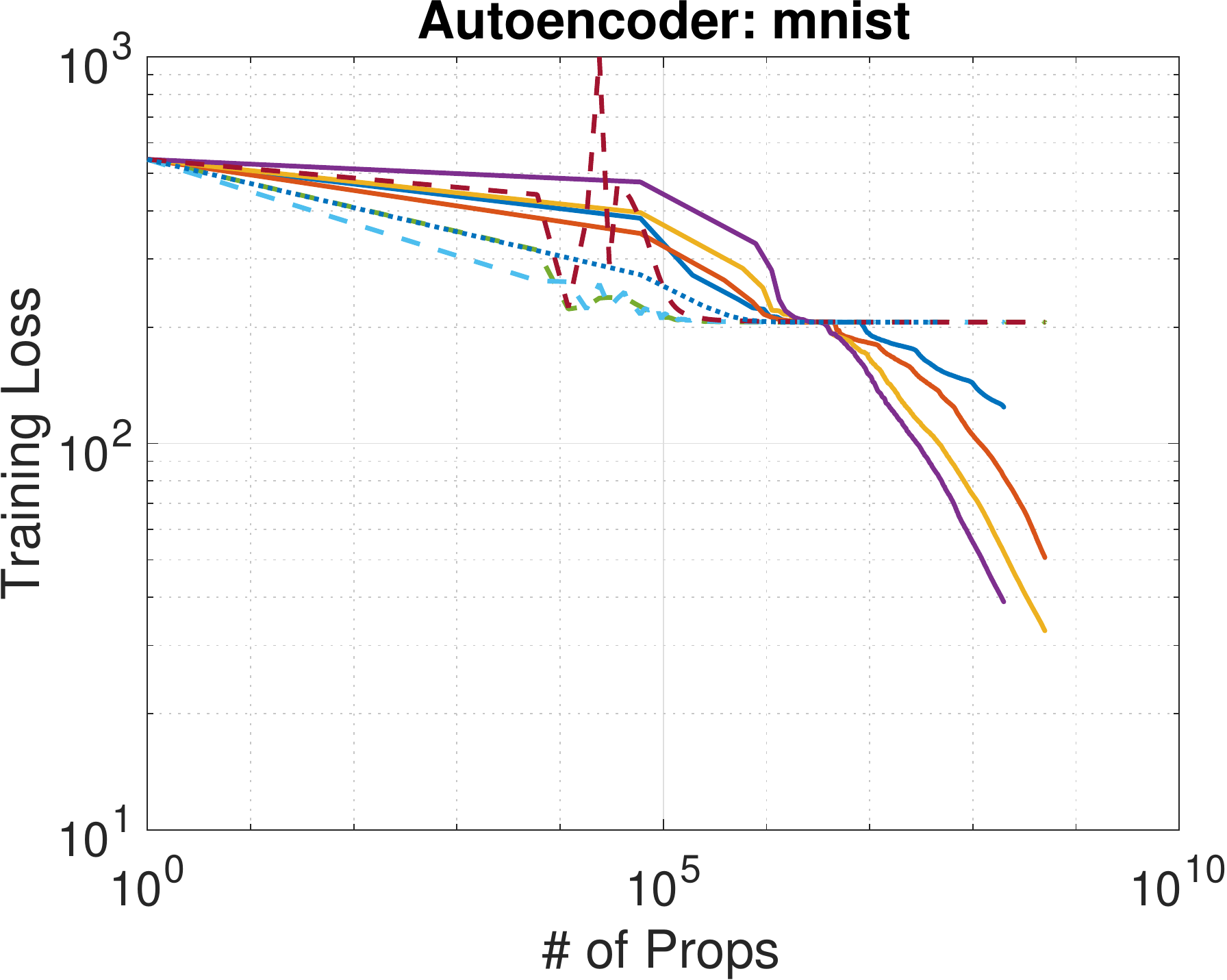}
			\includegraphics[width = 0.4\textwidth]{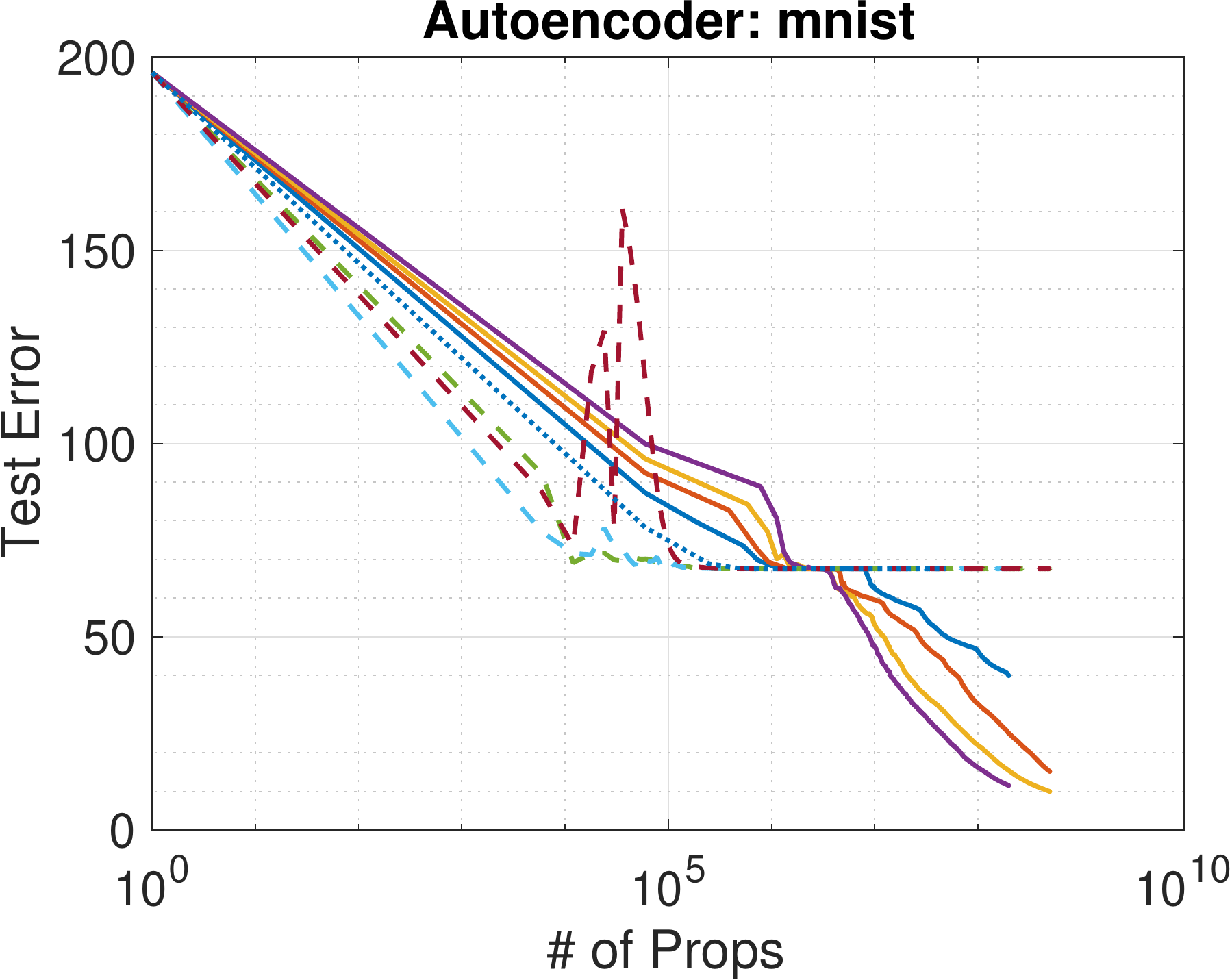}
		\subfigure[zero Initialization]{
		\centering

			\includegraphics[width = 0.65\textwidth]{figs/mnist/zeros0/mnist_zeros_0_legend.pdf}
	}
	\vspace{-3mm}
  	\caption{Deep autoencoder on \texttt{mnist} dataset using random (a) and zero (b) initial point.  $\Delta_{0}$ is the initial trust-region radius of Algorithm~\ref{alg:STR_fg} and $ \alpha $ is the step size for SGD with momentum.}
	\label{fig:mnist}
\end{figure}

Here, we consider the deep auto-encoder problem~\cite{goodfellow2016deep} and use the same model architectures as well as loss functions as in~\citet{martens2010deep}. The dataset and network architectures are given in Table \ref{tab: data2}.
The experiments in Figures~\ref{fig:curves} and \ref{fig:mnist} are each done with initialization to a vector drawn from standard normal distribution as well as the all-zeros vector.

\begin{table}[!htbp]
\centering
\small
\caption{Datasets for deep auto-encoder experiment.}\label{tab: data2}
\begin{tabular}{cccc}
\hline
Dataset & Size &  Encoder Network  ~(\# parameters)\\
\texttt{curves} & $20000$  &$784$-$400$-$200$-$100$-$50$-$25$-$6$ ~($842,340$)\\
\texttt{mnist} & $60000$  &$784$-$1000$-$500$-$250$-$30$ ~($2,837,314$)\\
\hline
\end{tabular}
\end{table}

Although deep auto-encoder networks here are much more complex than 1-hidden layer network of Section \ref{sec:example_1layer}, but we can still make very similar observations as they relate to \ref{question:efficient}--\ref{question:generalization}. 
In particular, on both datasets, Algorithm \ref{alg:STR_fg} converges comparably as fast, or faster than, other methods in terms of number of propagations. More importantly, it exhibits great robustness to its hyper-parameter, the initial trust-region radius, in contrast to SGD with momentum, which is heavily dependent on the choice of step size. Since these are rather complex networks, the optimization landscape is riddled with saddle and/or very flat regions. In this light, SGD and GN algorithms can both, rather easily, get trapped at/near these regions, which is quite contrary to the behavior of Algorithm~\ref{alg:STR_fg} (Figures~\ref{fig:curves}(b) and \ref{fig:mnist}(b)). In terms of test error, Algorithm~\ref{alg:STR_fg} is competitive to other methods on some experiments (Figures~\ref{fig:curves}(a) and \ref{fig:mnist}(a)) and clearly outperforms on others (Figures~\ref{fig:curves}(b) and \ref{fig:mnist}(b)).

\begin{figure}[htb]
	\centering
	\subfigure[Normalized Random Initialization]{
		\includegraphics[width = 0.4\textwidth]{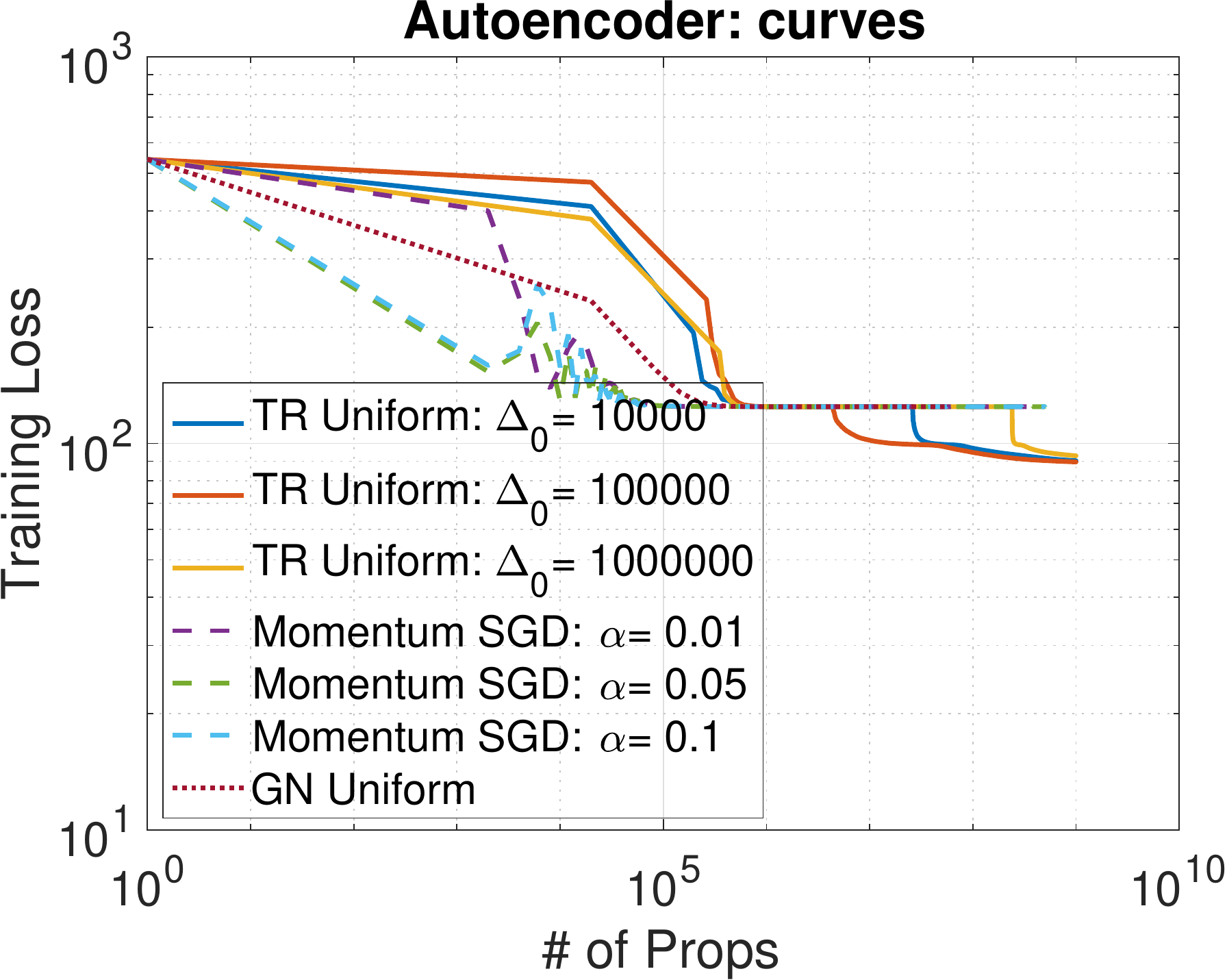}
		\includegraphics[width = 0.4\textwidth]{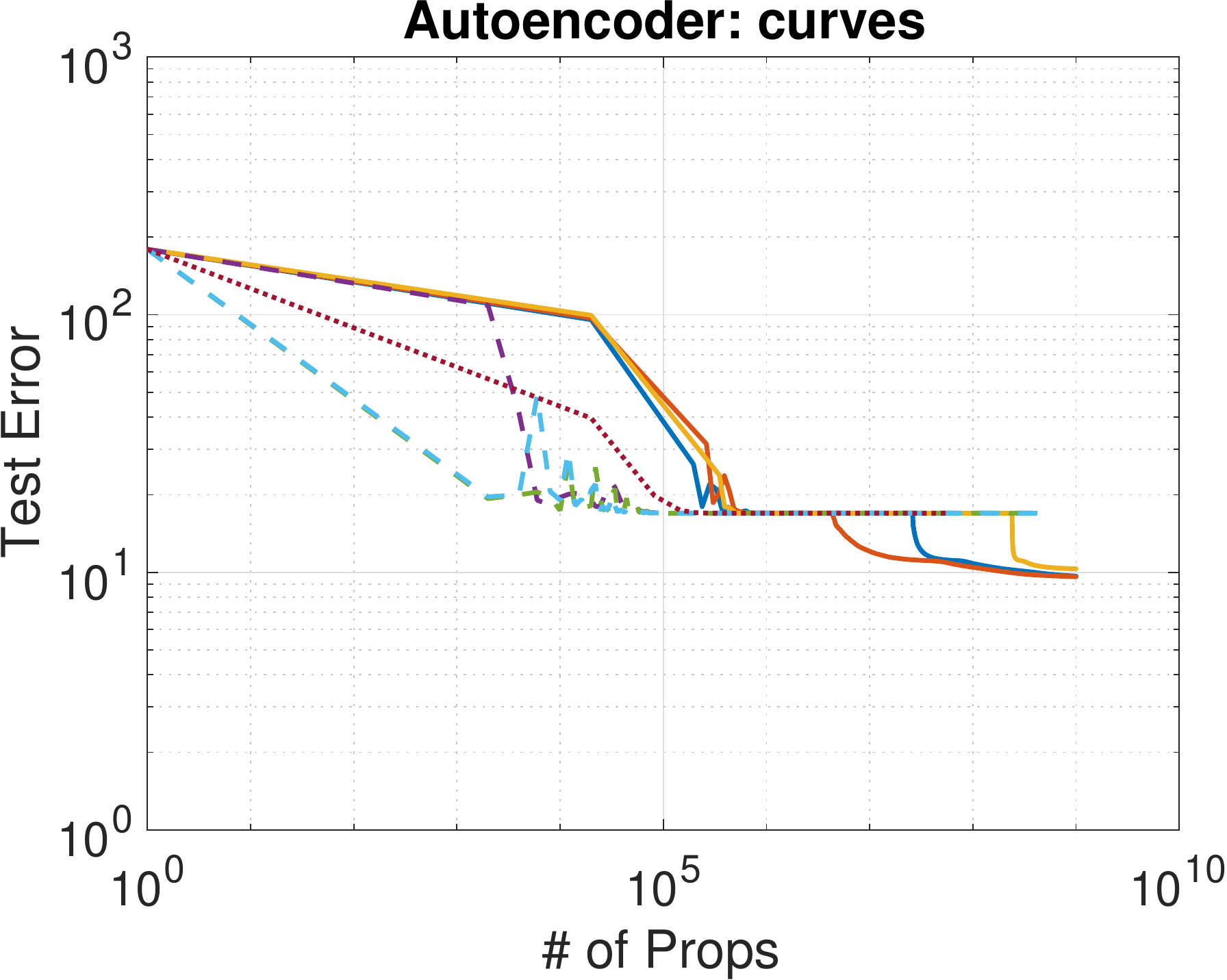}
	}
	\subfigure[Scaled Random Initialization]{
		\includegraphics[width = 0.4\textwidth]{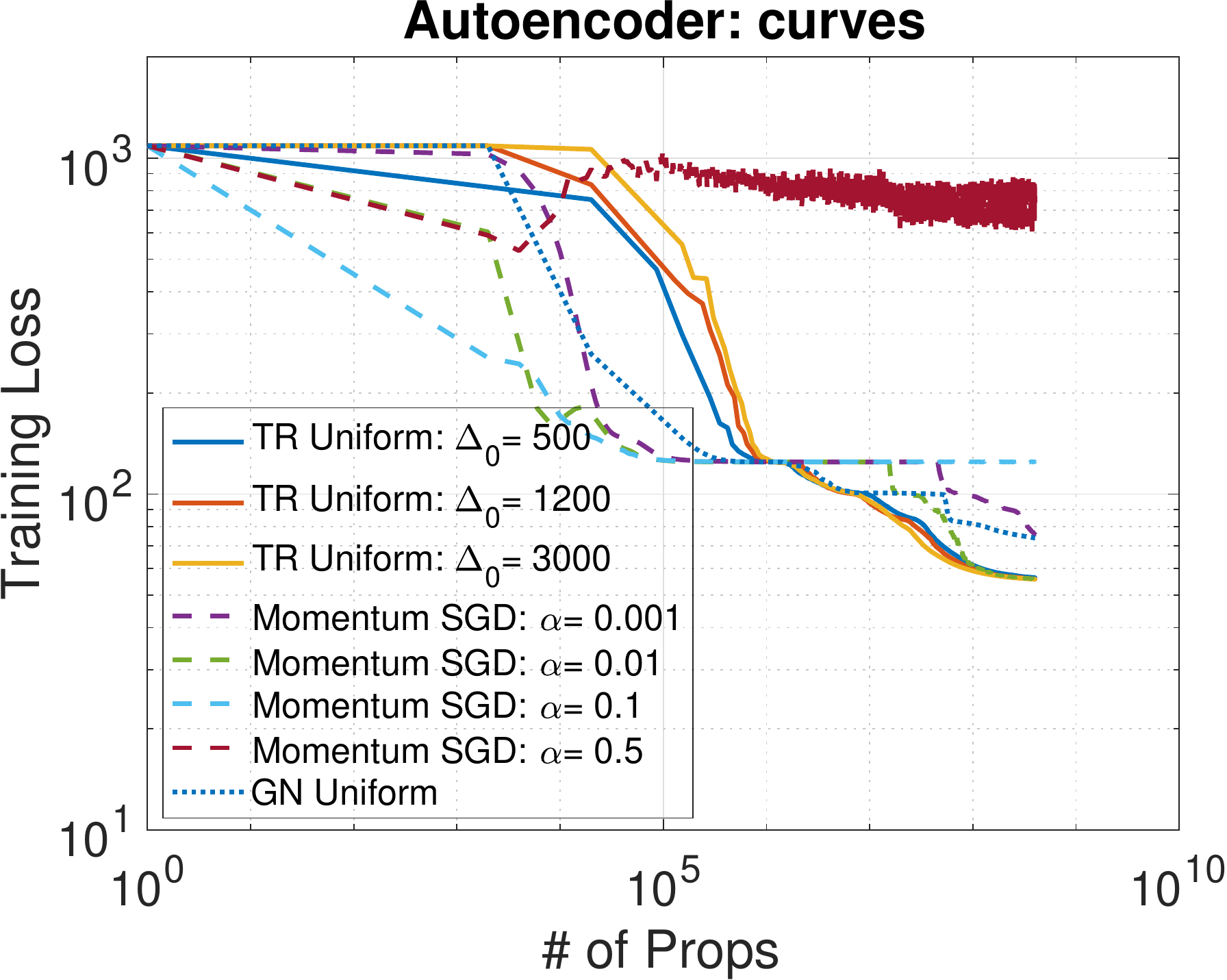}
		\includegraphics[width = 0.4\textwidth]{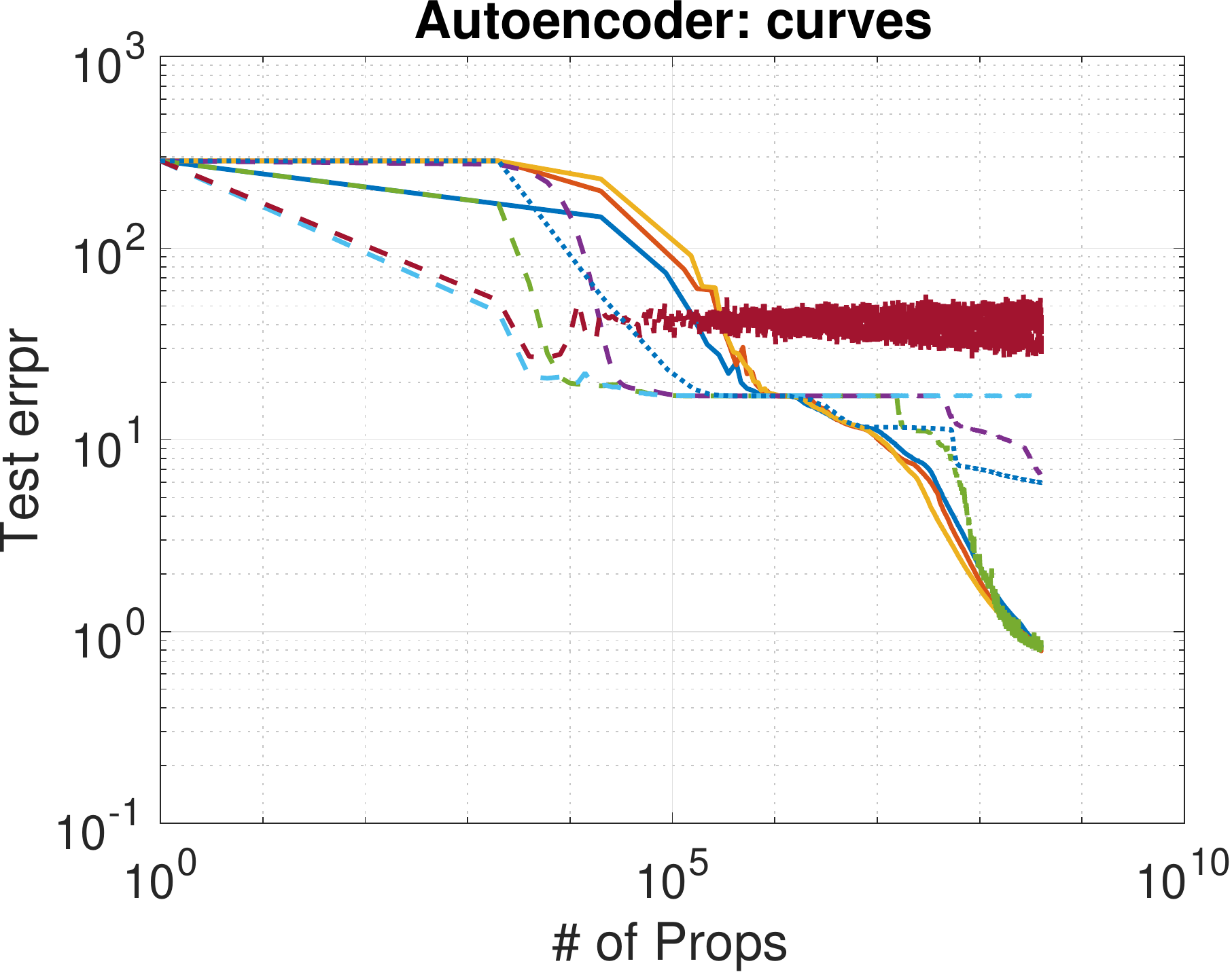}
	}
	\vspace{-3mm}
	\caption{Deep autoencoder on \texttt{curves} dataset using (a) normalized and (b) scaled random initial points.  $\Delta_{0}$ is the initial trust-region radius of Algorithm~\ref{alg:STR_fg} and $ \alpha $ is the step size for SGD.}
	\label{fig:curves2}
\end{figure}

We further examine two additional initialization strategies: normalized random initial point (Figure \ref{fig:curves2}(a)) as well as a random vector scaled by $ 0.25 $ (Figure  \ref{fig:curves2}(b)). 
Similar observations as earlier can also be made here. In particular, unlike Section~\ref{sec:example_1layer}, for this problem normalized random initial point (Figure \ref{fig:curves2}(a)) seems to paint a different picture, i.e., SGD with momentum as well as GN all get trapped at high training levels while Algorithm~\ref{alg:STR_fg} makes continued progress. It is worth mentioning that among all the initialization schemes we considered here, only under the particular scaled random initialization SGD can obtain desirable performance. This further demonstrates the versatility of Algorithm~\ref{alg:STR_fg}.
\vspace{-2mm}
\subsection{Non-Linear Least Squares}
\label{sec:example_nls}

We now turn to study Questions \ref{question:sampling} and \ref{question:comparison}. For this, we consider the class of NLS problems, and focus on the task of binary classification with square loss as a concrete instance. Since logistic loss, which is the ``standard'' loss used in this task, leads to a convex problem, we use square loss to obtain a non-convex objective. More specifically, suppose we are given training data $\left\{\x_i, y_i\right\}_{i=1}^n$, where $\x_i \in \bbR^d$ and $y_i\in\{0,1\}$ are, respectively, $ i^{\text{th}} $ feature vector and the corresponding label. Consider minimizing the empirical risk problem $ \min_{\w \in \bbR^d} \sum_{i=1}^n \left(y_i - \phi\big(\lin{\x_i,\w}\big)\right)^2/n$, 
where $\phi(z)$ is the sigmoid function, i.e., $\phi(z) = {1}/{(1 + e^{-z})}$.
Since~ this is an example of \eqref{eq:obj_sum_ERM}, we can apply both uniform and non-uniform sampling schemes.

Table \ref{tab:data1} summarizes the real data sets used for the experiments of this section. All datasets are from \texttt{LIBSVM} library~\cite{libsvm}. 
The following algorithms are compared (exact Hessian refers to $ \H_t = \nabla^{2} F(\x_{t}) $): 
\begin{enumerate}[leftmargin=*,wide=0em, itemsep=-1pt,topsep=-1pt,label=(\arabic*)]
\item {\it TR Full/Uniform/Non-Uniform}: Algorithm \ref{alg:STR_fg} with full or uniform/non-uniform estimation of Hessian, respectively,
\item {\it ARC Full/Uniform/Non-Uniform}: Algorithm~\ref{alg:SCR_fg} with full or uniform/non-uniform estimation of Hessian, respectively,
\item {\it GN Full/Uniform/Non-Uniform}: GN with full or uniform/non-uniform estimation of Hessian, respectively,
\item {\it LBFGS-100}: the standard L-BFGS method~\citep{liu1989lbfgs} with history size $100$ and using line-search.
\end{enumerate}
\vspace{-4mm}
\begin{table}[htbp]
	\caption{Datasets used in binary linear classification.}
	\label{tab:data1}
	\centering
	\small
	\begin{tabular}{cccc}
		\toprule
		\sc Data & Training Size ($n$) & \# Features ($d$) & Test Size  \\ 
		\midrule
		{\tt covertype} & $464,810$ & $54$ & $116,202$
		\\
		{\tt mnist} & $60,000$& $784$ &$10,000$
		\\
		\bottomrule
	\end{tabular}
\end{table}

\begin{figure}[htb]
\centering
\includegraphics[scale=0.35]{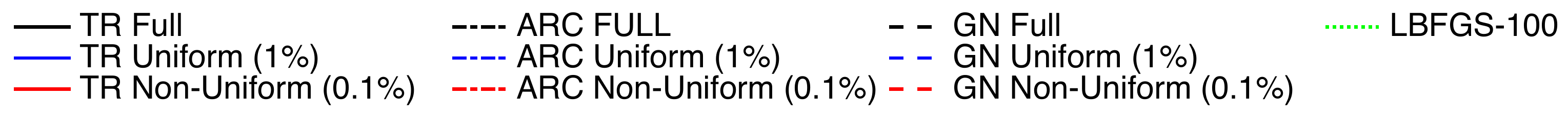}

\includegraphics[width=0.244\textwidth]{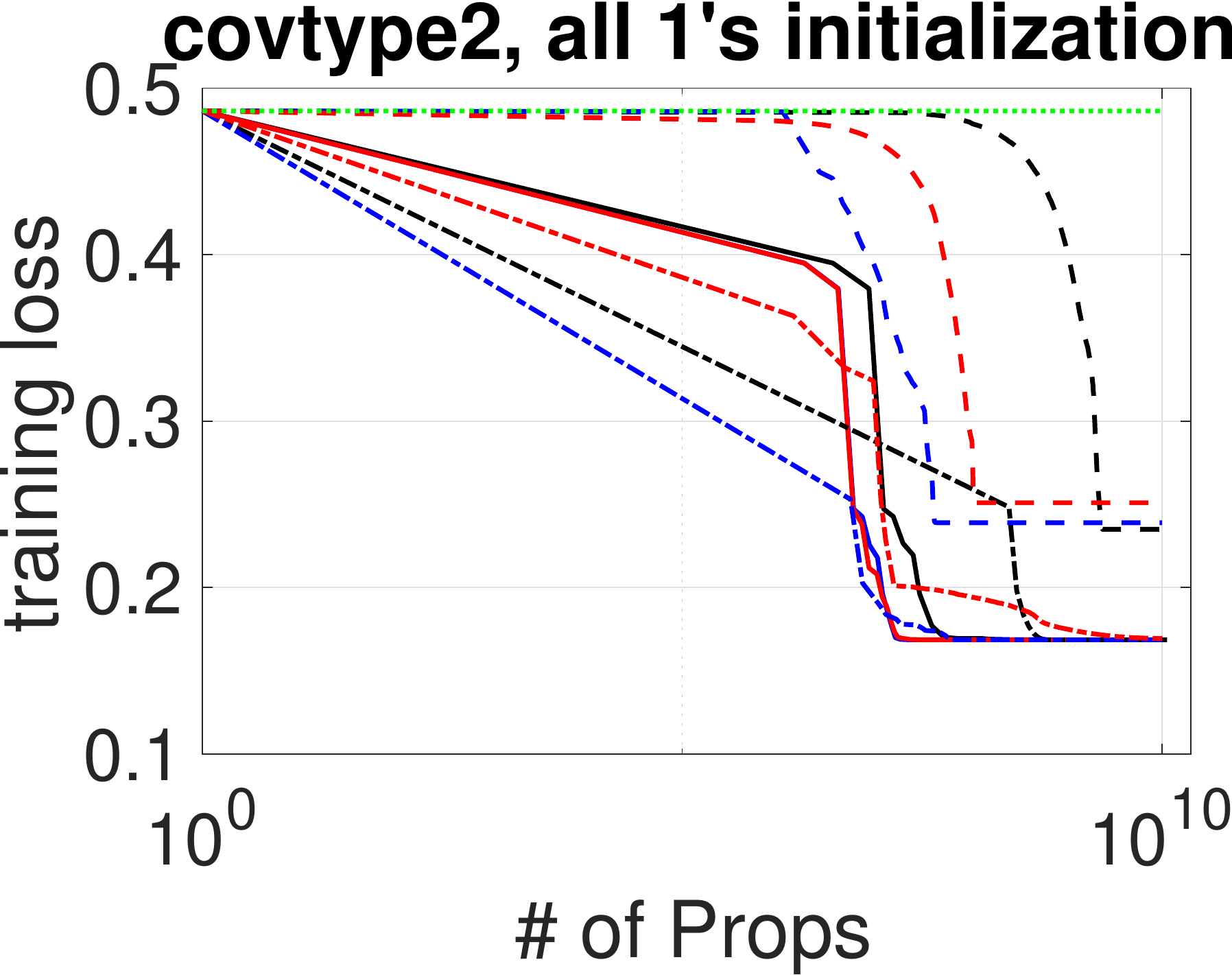}
\includegraphics[width=0.244\textwidth]{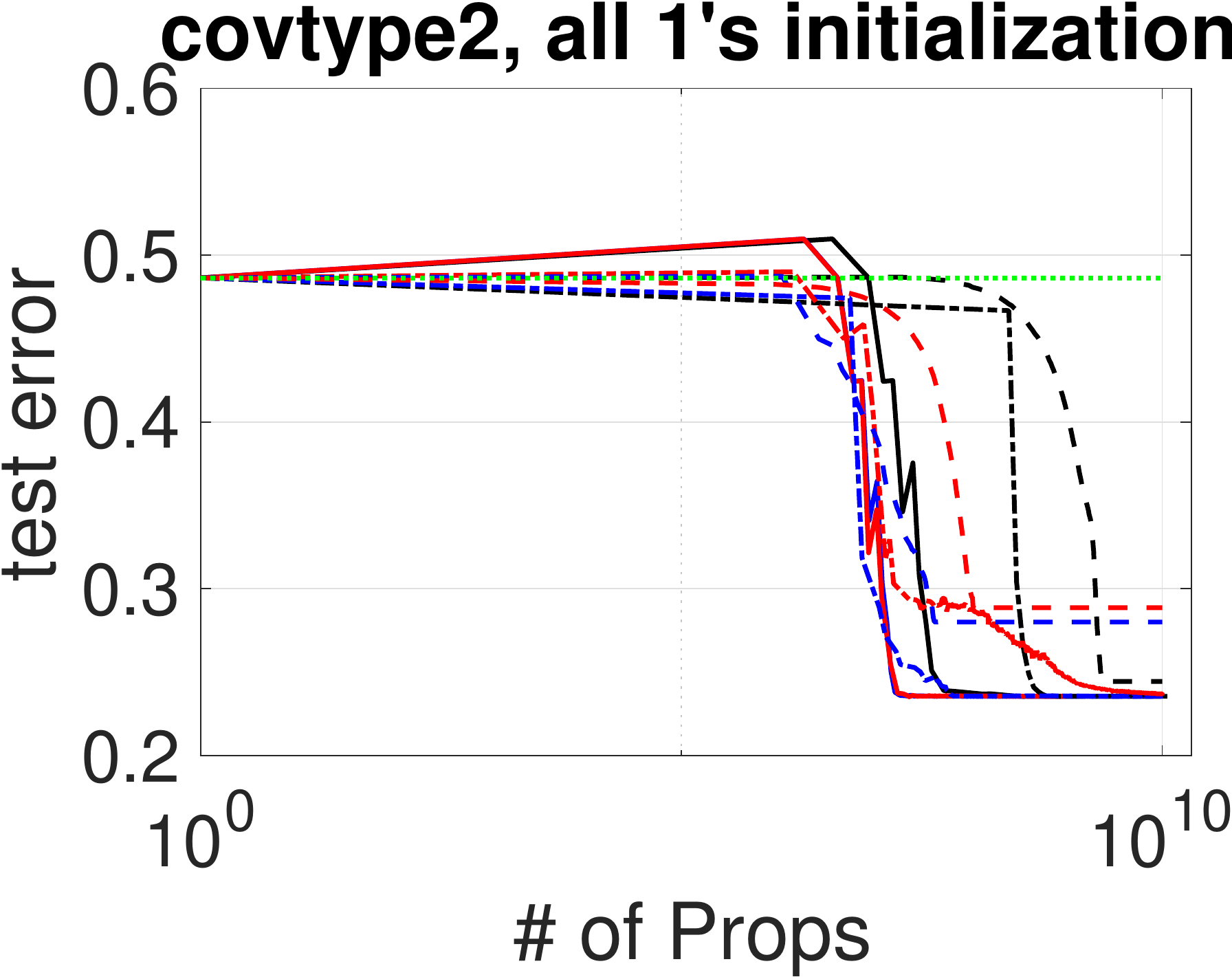}
\includegraphics[width=0.244\textwidth]{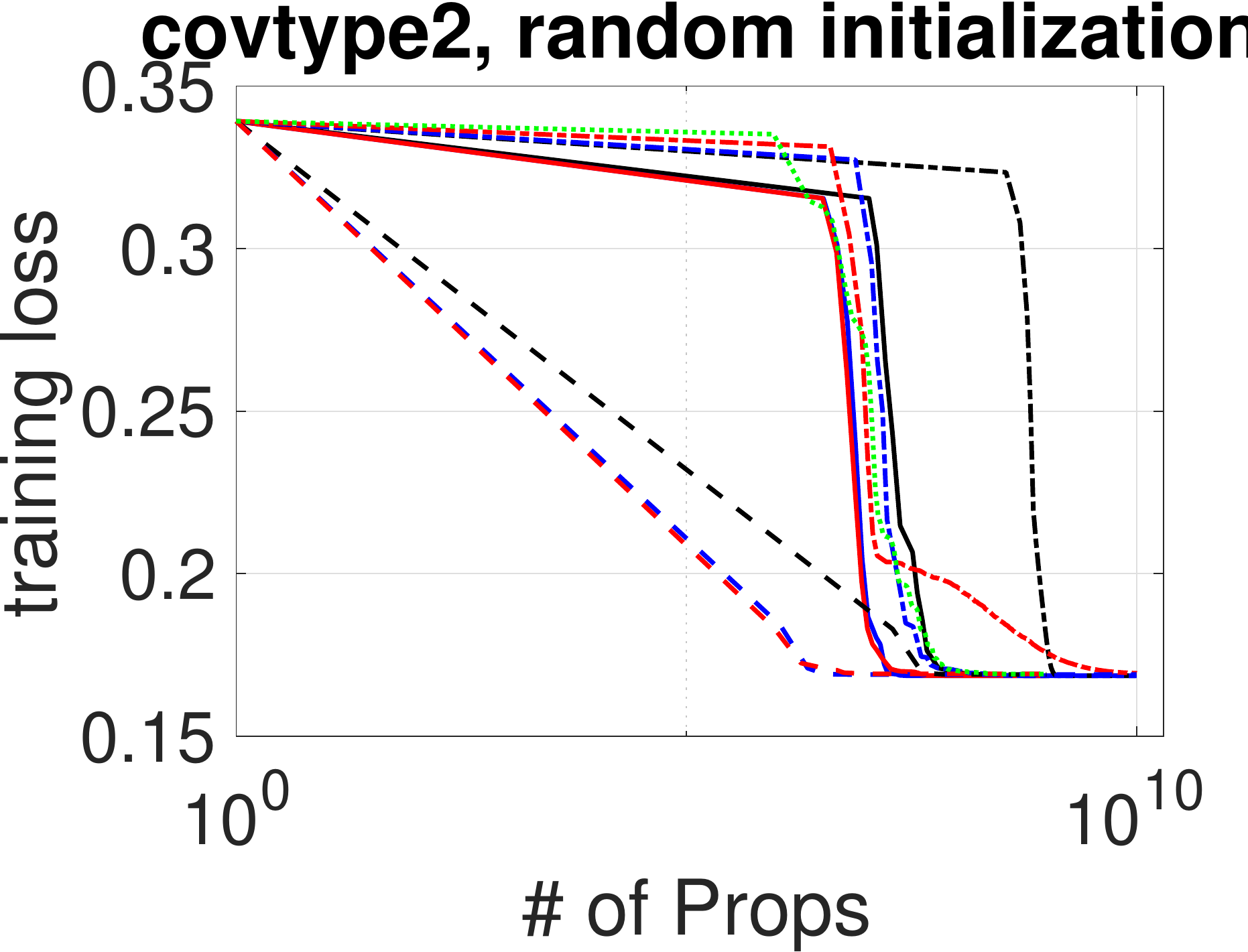}
\includegraphics[width=0.244\textwidth]{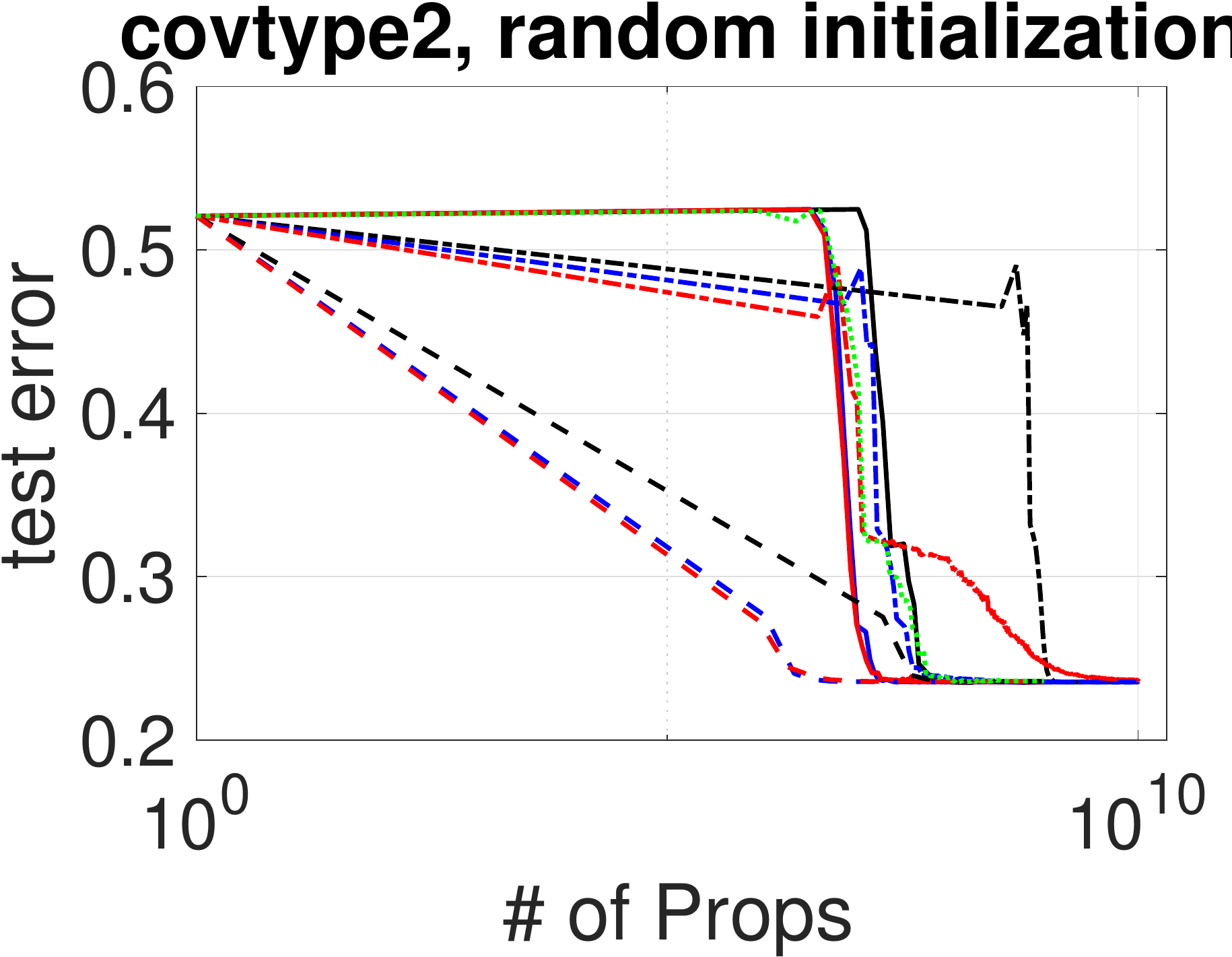}

\includegraphics[width=0.243\textwidth]{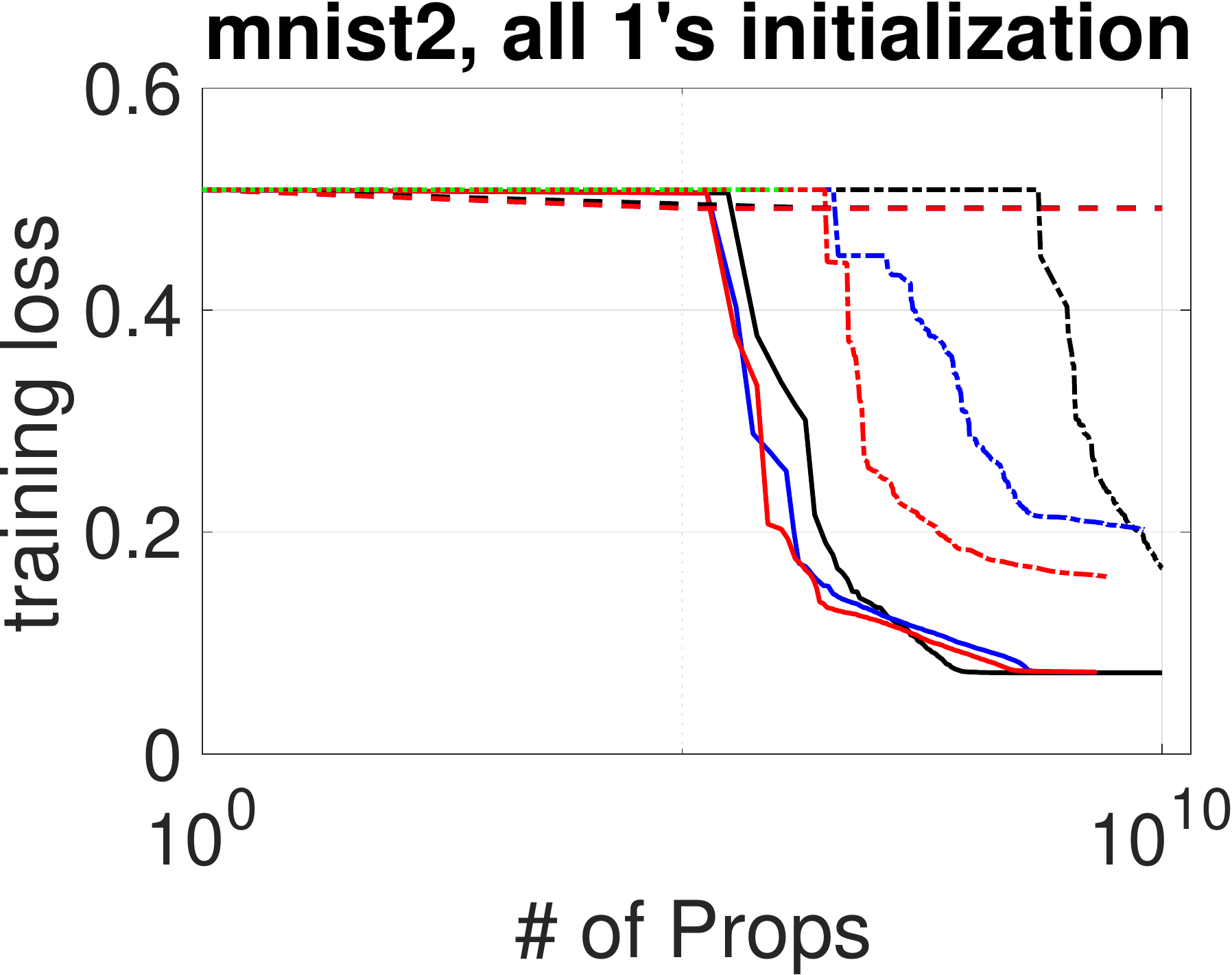}
\includegraphics[width=0.243\textwidth]{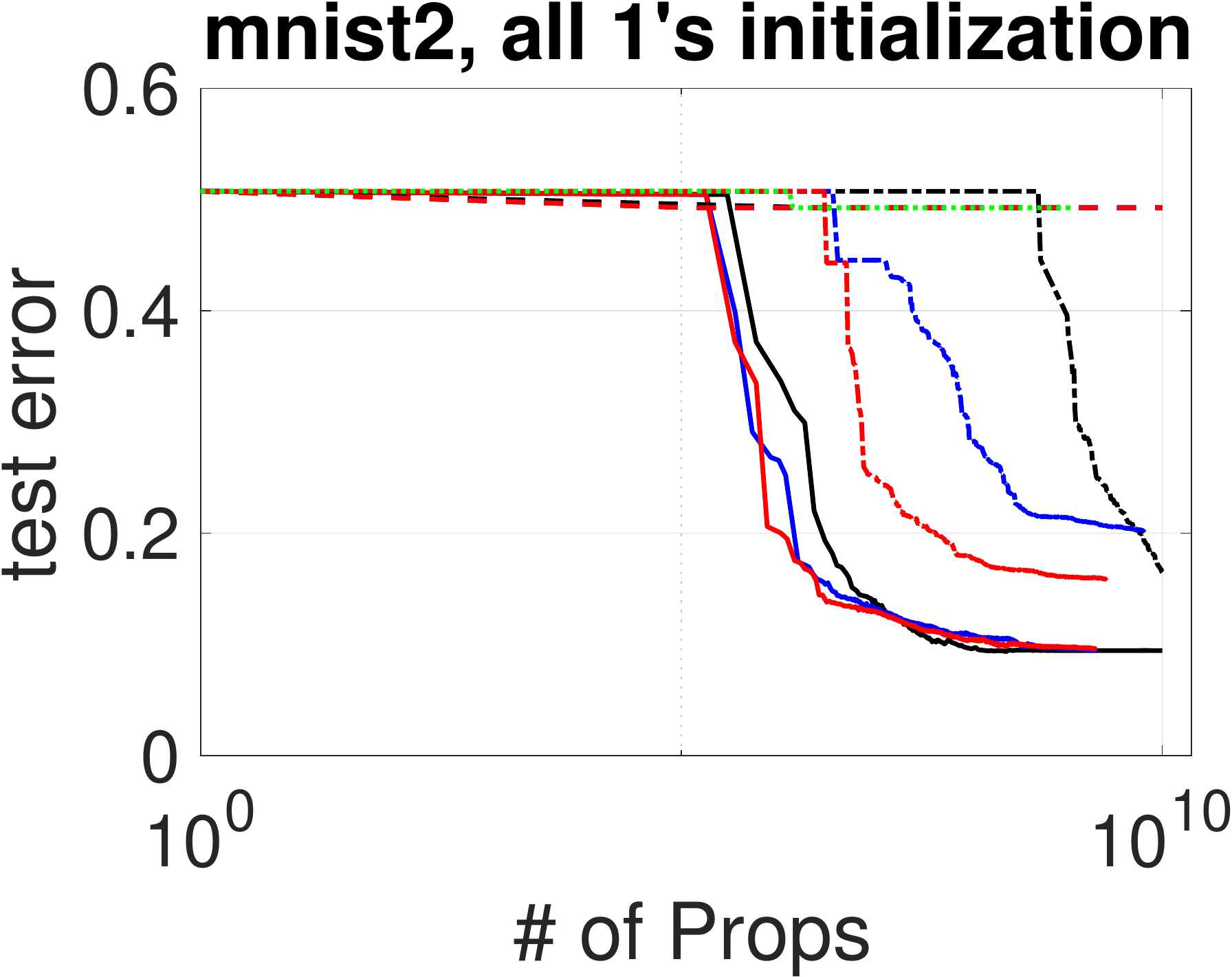}
\includegraphics[width=0.243\textwidth]{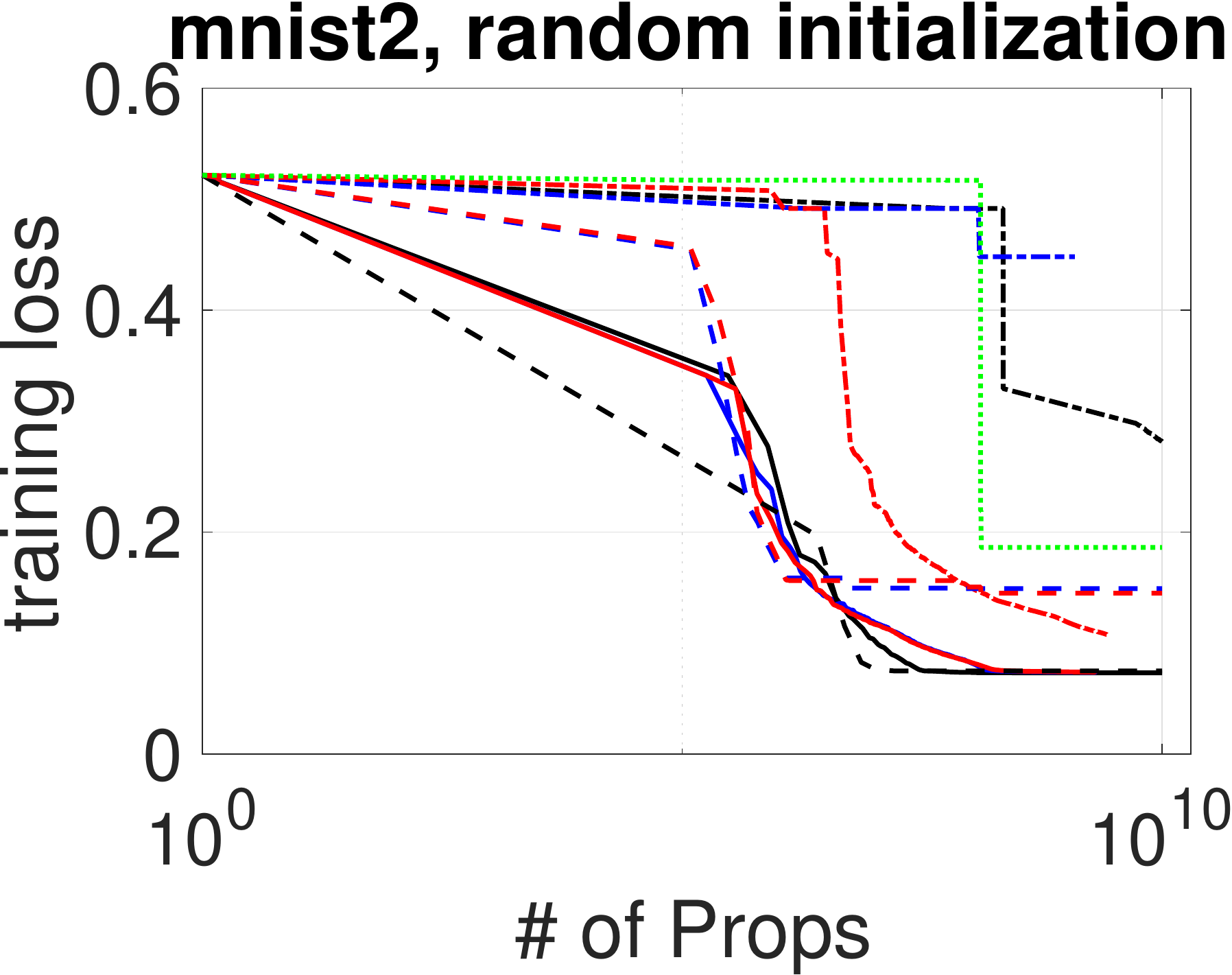}
\includegraphics[width=0.243\textwidth]{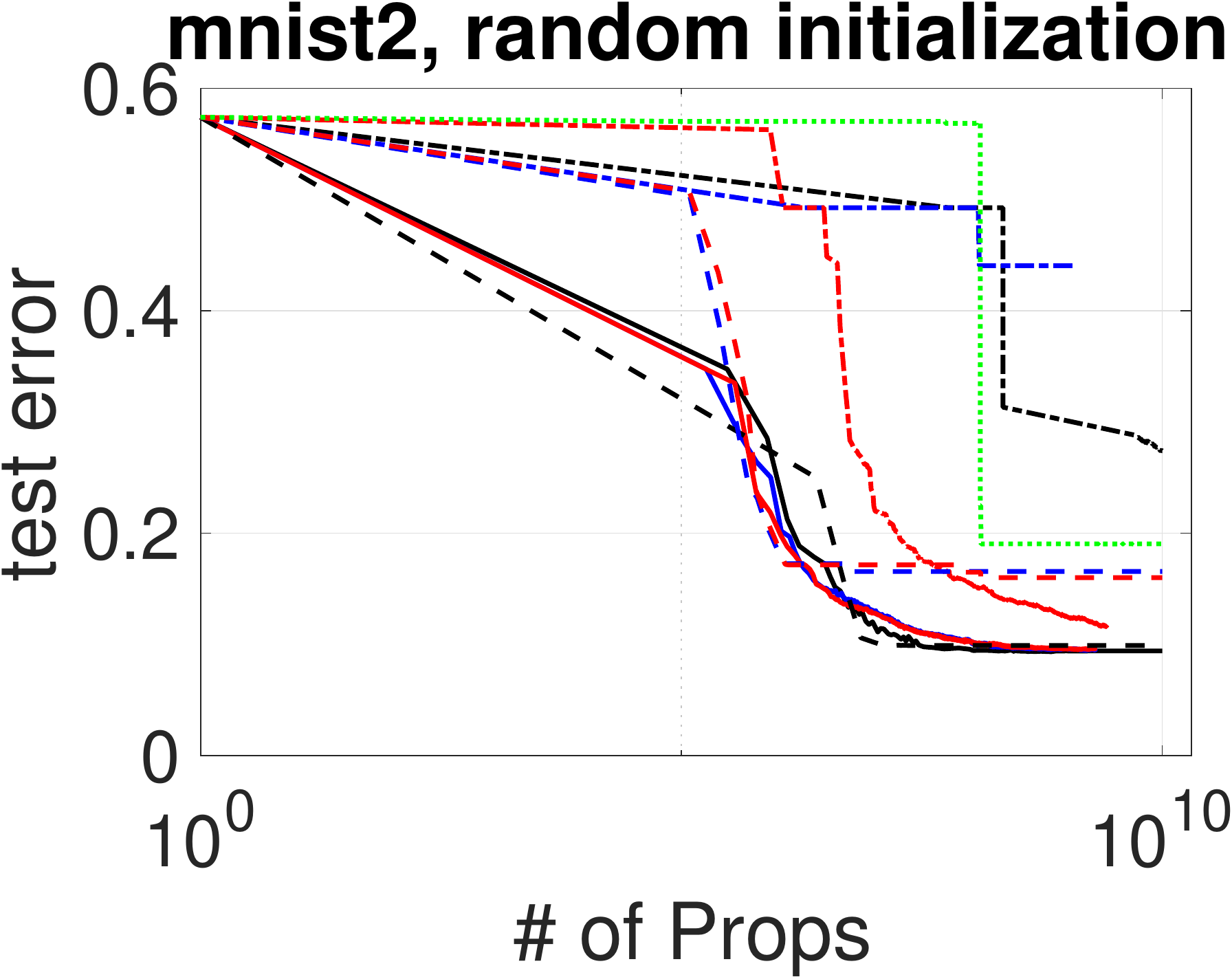}
\vspace{-3mm}
\caption{Binary Classification on Different Datasets with all $1$'s initialization (left two columns) and randon initialization (right two columns). Each row corresponds to a dataset. 
The x-axis is the number of propagations in log-scale.}
\label{fig:blc}
\end{figure}

\vspace{-2mm}
\paragraph{Parameter settings} For both Algorithms \ref{alg:STR_fg} and \ref{alg:SCR_fg}, we set $\eta_1, \eta_2, \gamma_1$ and $\gamma_2$, the same as Section~\ref{sec:example_deep_learning}. The sampling ratios, i.e., $ |\mathcal{S}|/n$, for uniform and non-uniform sampling are set to $ 1\% $ and $ 0.1\% $, respectively. For all datasets,  we set $ \Delta_{0} = 10$ for Algorithm~\ref{alg:STR_fg} and $ \sigma_{0} = 10^{-4}$ for Algorithm~\ref{alg:SCR_fg}. For GN, we do not regularize Hessian as in~\citet{martens2010deep}. Figure \ref{fig:blc} gathers all comparison results of this section.

\begin{enumerate}[leftmargin=*,wide=0em, itemsep=-1pt,topsep=-1pt,label = {\bfseries Re: Q.\arabic*}, resume]
\item \textit{(Benefits of Sub-Sampling)} Overall, both Algorithms \ref{alg:STR_fg} and \ref{alg:SCR_fg} compare well with the classical TR and ARC methods. In fact, sub-sampling can, at times, help increase the efficiency, e.g., TR variants for {\tt covtype2}. However, with too small a sample, the performance can hurt; ARC Full vs.\ ARC Non-Uniform and ARC Uniform for {\tt covtype2} and {\tt mnist2}, respectively. 
The benefits of non-uniform sampling over uniform  alternative are far more pronounced in the performance of Algorithm~\ref{alg:SCR_fg} than Algorithm~\ref{alg:STR_fg}. This can be attributed mainly to their respective sub-problem solvers in terms of total number of performed Hessian-vector products . In particular, CG-Steihaug used for the sub-problem \eqref{eq:STR_subp} of Algorithm~\ref{alg:STR_fg} typically terminates in a handful of iterations whereas the generalized Lanczos method for solving the sub-problem~\eqref{eq:SARC_subp} of Algorithm~\ref{alg:SCR_fg} usually exhausts the allotted 250 iterations.

\item \textit{(Comparison Among Second-Order Methods)}  One can observe the consistent poor performances of L-BFGS-100 (green dotted lines) methods on all datasets, in particular with all $1$'s initialization vector. This is rather expected as contrary to popular belief, BFGS is not quite a ``full-fledged'' second-order method. Indeed, BFGS merely employs first-order information, i.e. gradients, to approximate the curvature, and starting from all $1$'s vector, L-BFGD cannot capture enough curvature information to navigate its way out of this region effectively. Gauss-Newton (dash lines), which has been specifically designed to effectively solve NLS problems, performs very well with random initialization. Starting from all $1$'s vector, however, where the gradient is very small, GN performs poorly. This is also expected because GN, similar in spirit to BFGS, does not fully utilize the Hessian information. In particular, in exchange for obtaining a positive  definite approximation matrix, GN completely ignores the information from \emph{negative curvature}, which is critical for allowing to escape from regions with small gradient. We can also observe that ARC is consistently no better than TR. This is an empirical evidence that the optimal worst-case complexity of ARC \citep{cartis2011adaptiveII,xuNonconvexTheoretical2017}, though theoretically highly interesting, might be hard to observe in many practical settings.
\end{enumerate}

\vspace{-1mm}
\section{Conclusion}
\label{sec:conclusions}
In this paper, we aimed at painting a more complete picture of the practical advantages of Newton-type algorithms in general, and sub-sampled variants of trust-region and adaptive cubic regularization methods in particular, as compared with first-order alternatives. In the context of (deep) multi-layer perceptron networks as well as non-linear least squares, two simple, yet illustrative, non-convex machine learning applications,  by making the following observations, we empirically attempted to make a case for the application of such second-order methods for machine learning.
\begin{enumerate}[leftmargin=*,wide=0em, itemsep=-1pt,topsep=-1pt, label = {\bfseries A.\arabic*}]
    \item({\it Computational Efficiency}) The randomized sub-sampling approaches described here, and studied in detail in~\citet{xuNonconvexTheoretical2017}, can effectively make Newton-type methods \emph{computationally efficient} enough to be competitive with popular first-order methods, widely used in machine learning, e.g., SGD with momentum. This is indeed due to the amortized combination of (1) low per-iteration cost offered by randomized sampling, and (2) small number of overall iterations due to the application of curvature information.
	\item({\it Robustness to Hyper-parameters}) In contrast to first-order algorithms whose performance is greatly affected by the choice of hyper-parameters, most notably \emph{step-size}, the performance of the proposed Newton-type methods exhibit great  \emph{robustness} to such \emph{parameter tuning}.
	\item({\it Escaping Saddle Point}) A greatly beneficial advantage of employing Hessian information is that it allows for such Newton-type algorithms, unlike many first-order alternatives, to seamlessly escape regions near \emph{saddle points}.
	\item({\it Generalization Performance}) Second-order methods prove beneficial for the downstream machine learning objective of obtaining \emph{good generalization error}. In particular, one can obtain very good levels of prediction accuracy only after a few iterations of such methods. This is highly beneficial in, say, distributed settings where the communication across the network is the main computational bottleneck. 
	\item({\it Benefits of Sub-sampling})  On several real datasets, we validated the effectiveness of sub-sampled Newton-type methods, as compared with classical versions, in speeding up computations. Also, the advantages of non-uniform sampling over the oblivious uniform alternative was verified. 
	\item({\it Comparison Among Second-Order Methods}) There are clear advantages in using (sub-sampled) TR and ARC methods over other second-order alternatives, e.g., L-BFGS and GN, in terms of effective exploitation of curvature.
\end{enumerate}

For the examples of Section~\ref{sec:example_deep_learning}, despite the best of our efforts, we were unable to obtain the expected performance of Algorithm~\ref{alg:SCR_fg} using a variety of implementations, e.g., our own hand-written code as well as some existing packages such as GALAHAD~\cite{gould2003galahad}. 
We believe that this is tightly connected to the choice of the sub-problem solver in all these implementations. 
Since we were unable to pinpoint the source of the problem, we did not include Algorithm~\ref{alg:SCR_fg} in examples of Section~\ref{sec:example_deep_learning}.

Finally, we acknowledge that, although we presented various experiments, the empirical study of methods such as the ones considered here, takes more than a single ``proof-of-concept'' paper, and the results presented here should be viewed as merely a glimpse into their various properties. 

\section*{Acknowledgment}
We would like to acknowledge ARO, DARPA, Cray, and NSF for providing partial support of this work. FR gratefully acknowledges the support of the Australian Research Council through a Discovery Early Career Researcher Award (DE180100923). We also sincerely thank Profs.\ Dominique Orban, Nicholas I.M. Gould, and Coralia Cartis for kindly helping us with the code for adaptive cubic regularization as well as setting up GALAHAD package. We also greatly appreciate Dr.\ Felix Lenders' help with the installation of the trlib package. We would like to thank Dr.\ Amir Gholaminejad for valuable comments on our empirical evaluations and suggestions for improving them.

\printbibliography

\newpage
\onecolumn
\appendix
\section{Image Classification with \texttt{Cifar10}}\label{sec:cifarapp}

Here, we evaluate the sensitivity of various algorithms, when initialized with different random seeds. We do this by consider similar set up as in \cref{sec:example_1layer}. 
We present the 10 different runs of Algorithm~\ref{alg:STR_fg} and SGD with momentum with fixed the configuration. Overall, we observed that algorithms did not show much sensitivity with respect to random seed.

\begin{figure}[H]
	\centering
	\subfigure[9 runs of SGD]{
		\includegraphics[width=0.28\textwidth]{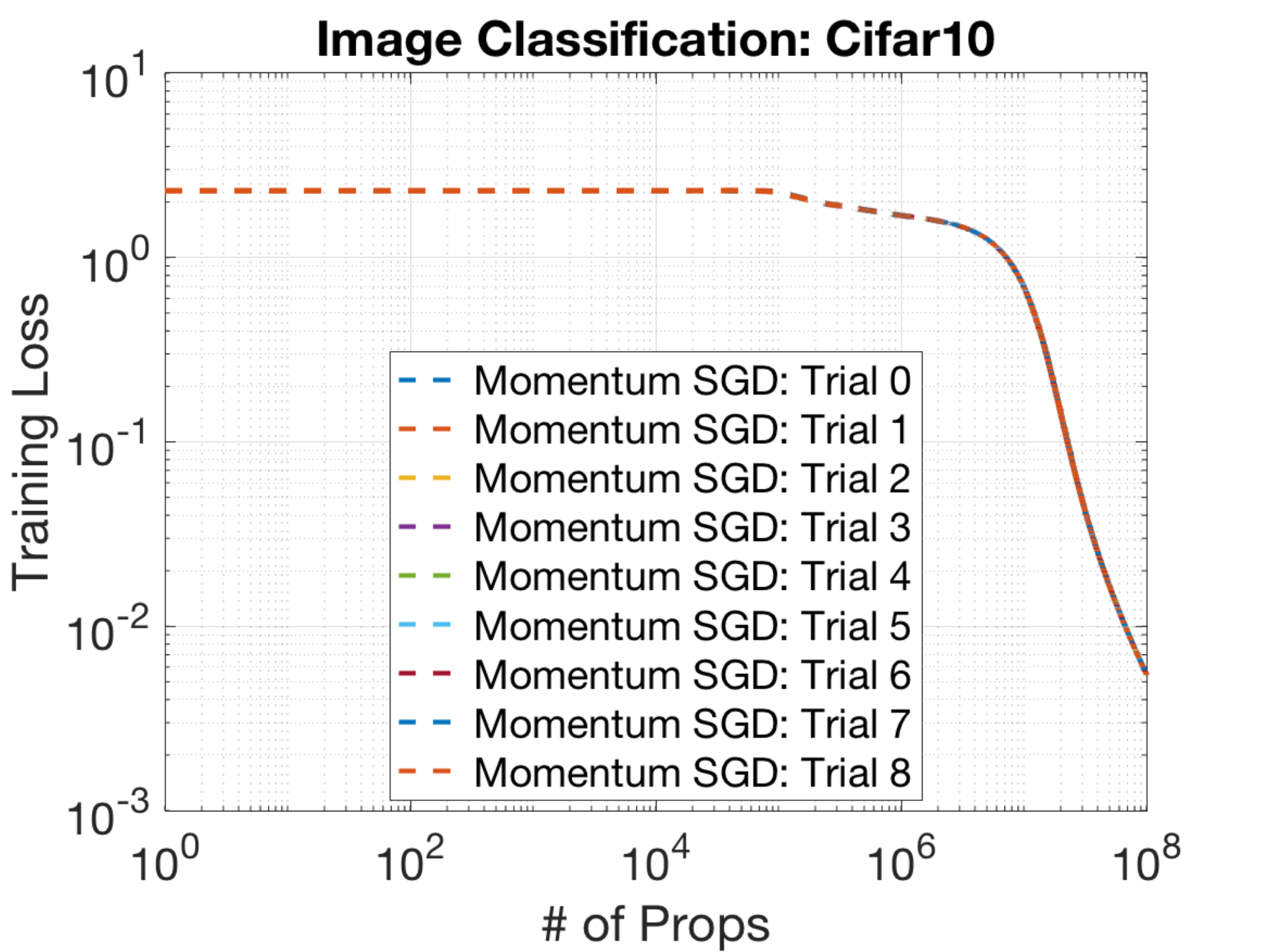}
		\includegraphics[width=0.28\textwidth]{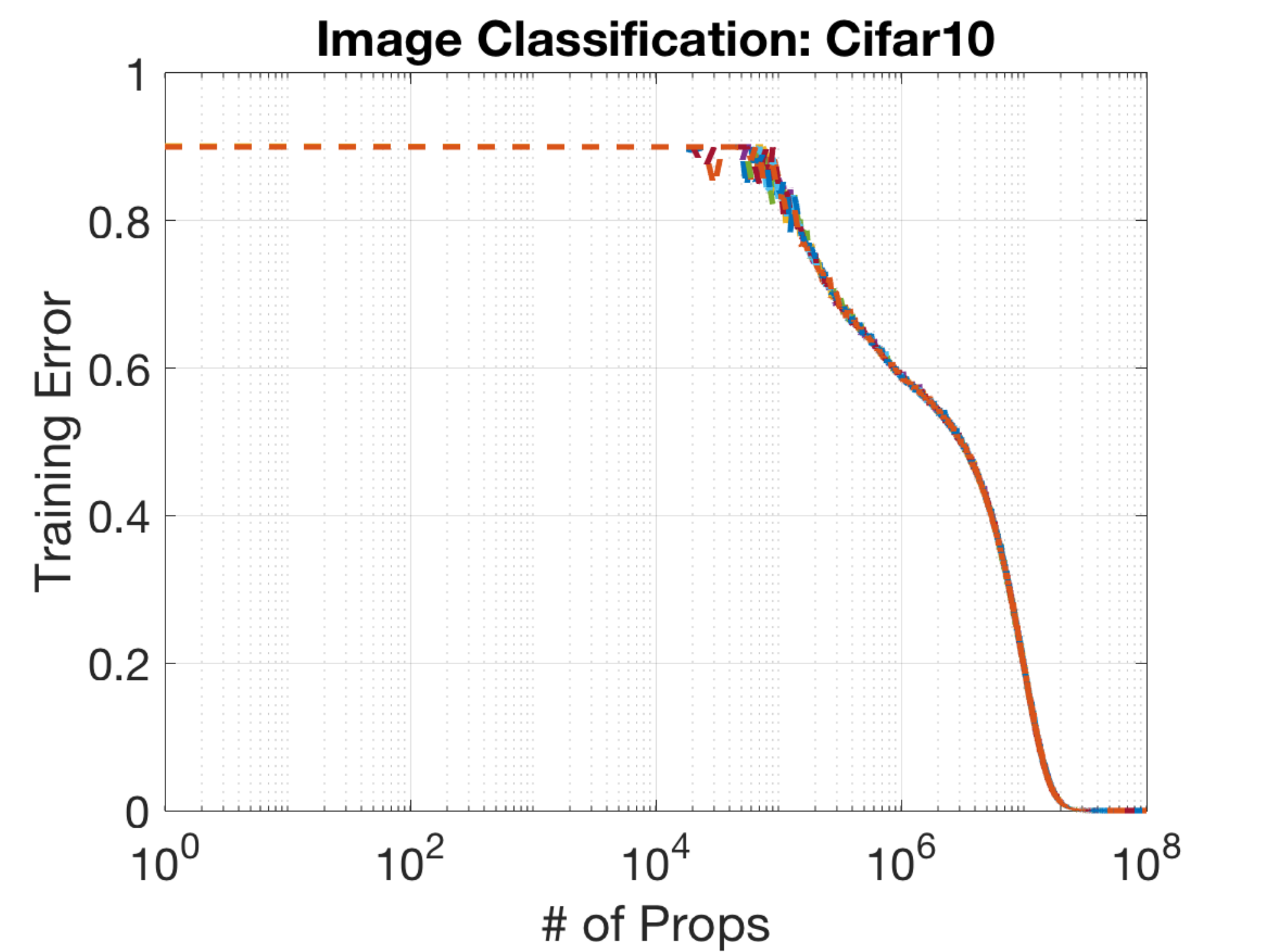}
		\includegraphics[width=0.28\textwidth]{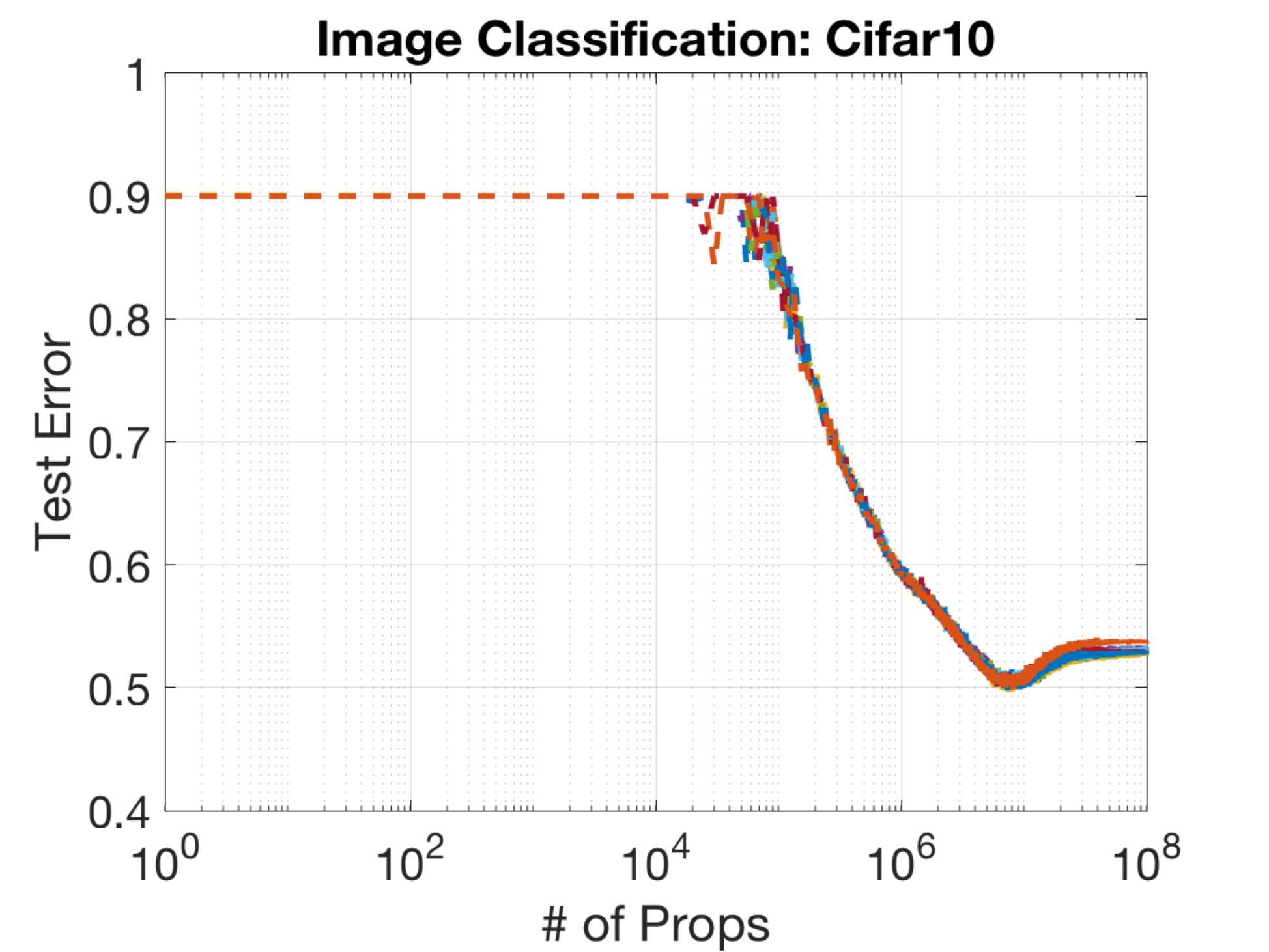}
	}
	\subfigure[9 runs of Sub-sampled TR]{
		\includegraphics[width=0.27\textwidth]{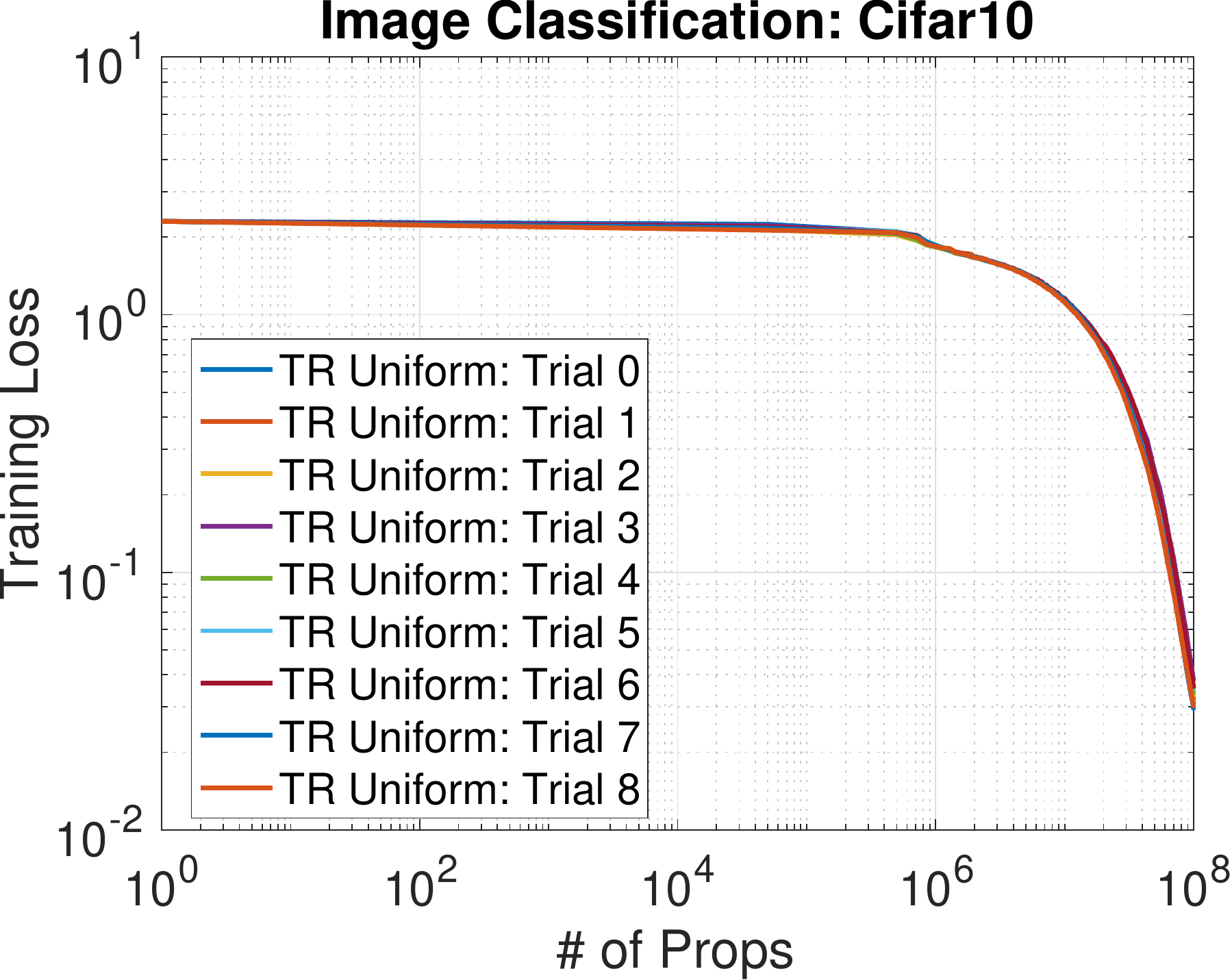}
		\includegraphics[width=0.27\textwidth]{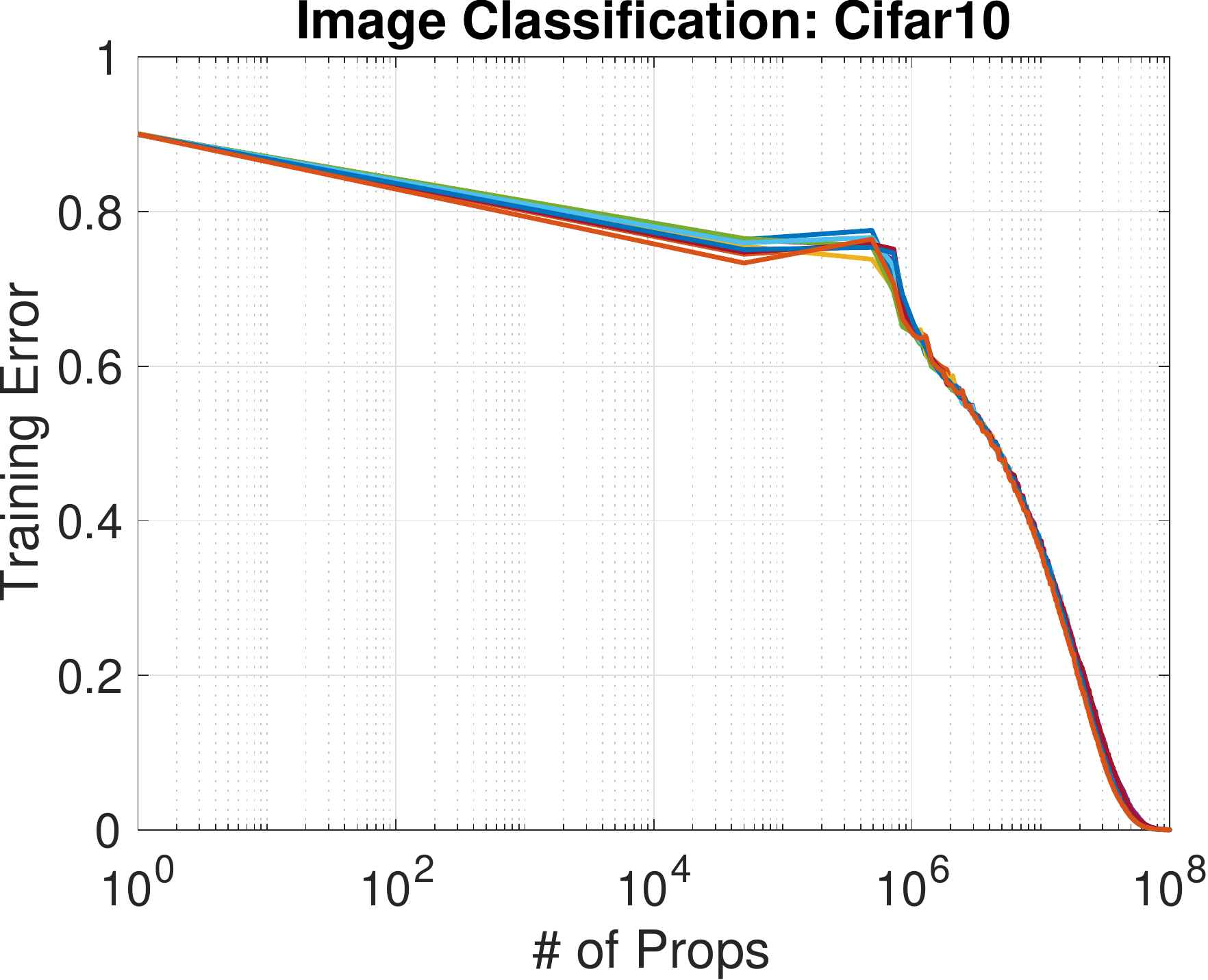}
		\includegraphics[width=0.27\textwidth]{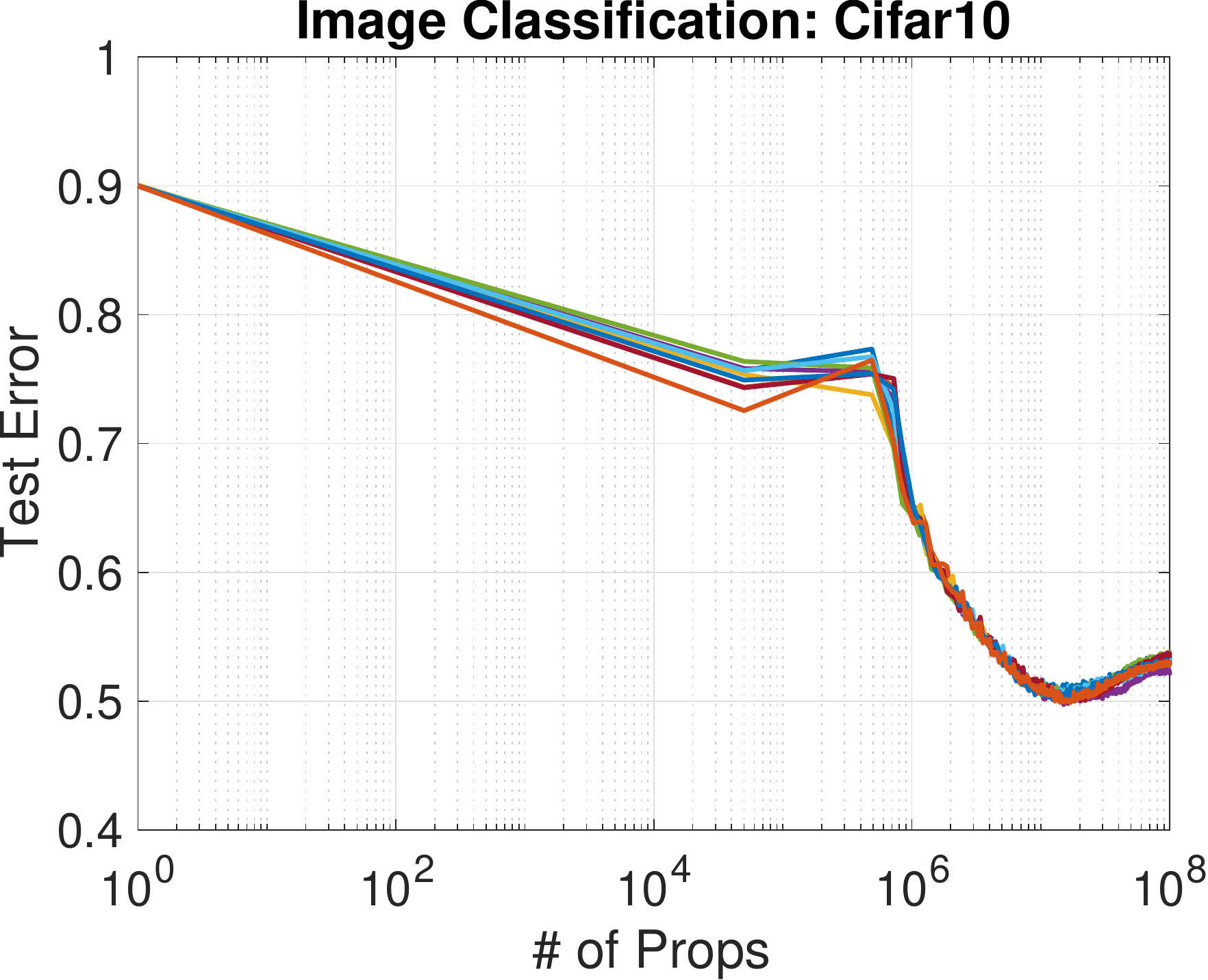}
	}
	\caption{Multiple runs of SGD and Algorithm~\ref{alg:STR_fg} with different random seeds. All runs start with normalized random initialization.}
	\label{fig:cifar3}
\end{figure}

\begin{figure}[H]
	\centering
	\subfigure[9 runs of SGD]{
		\includegraphics[width=0.28\textwidth]{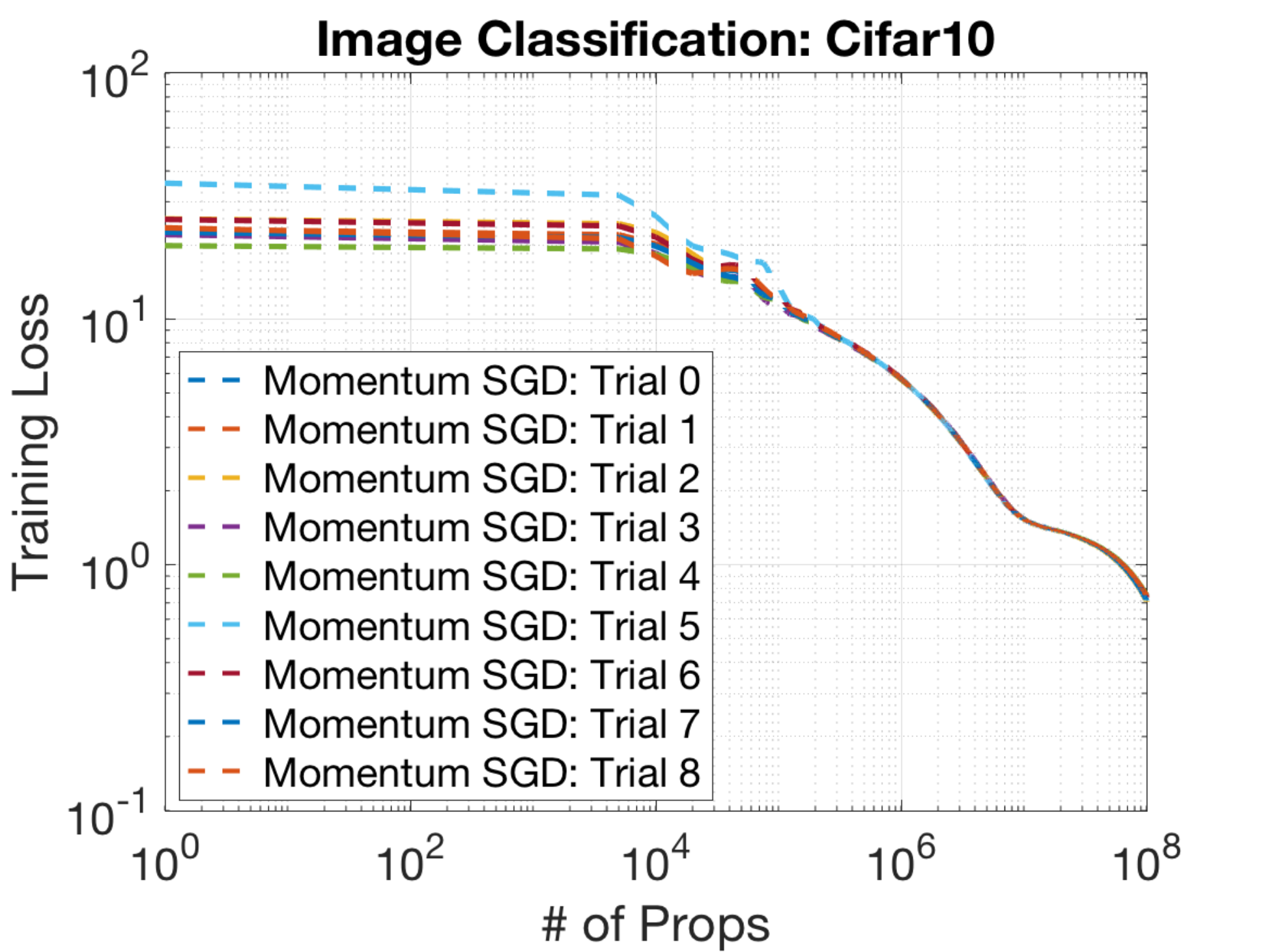}
		\includegraphics[width=0.28\textwidth]{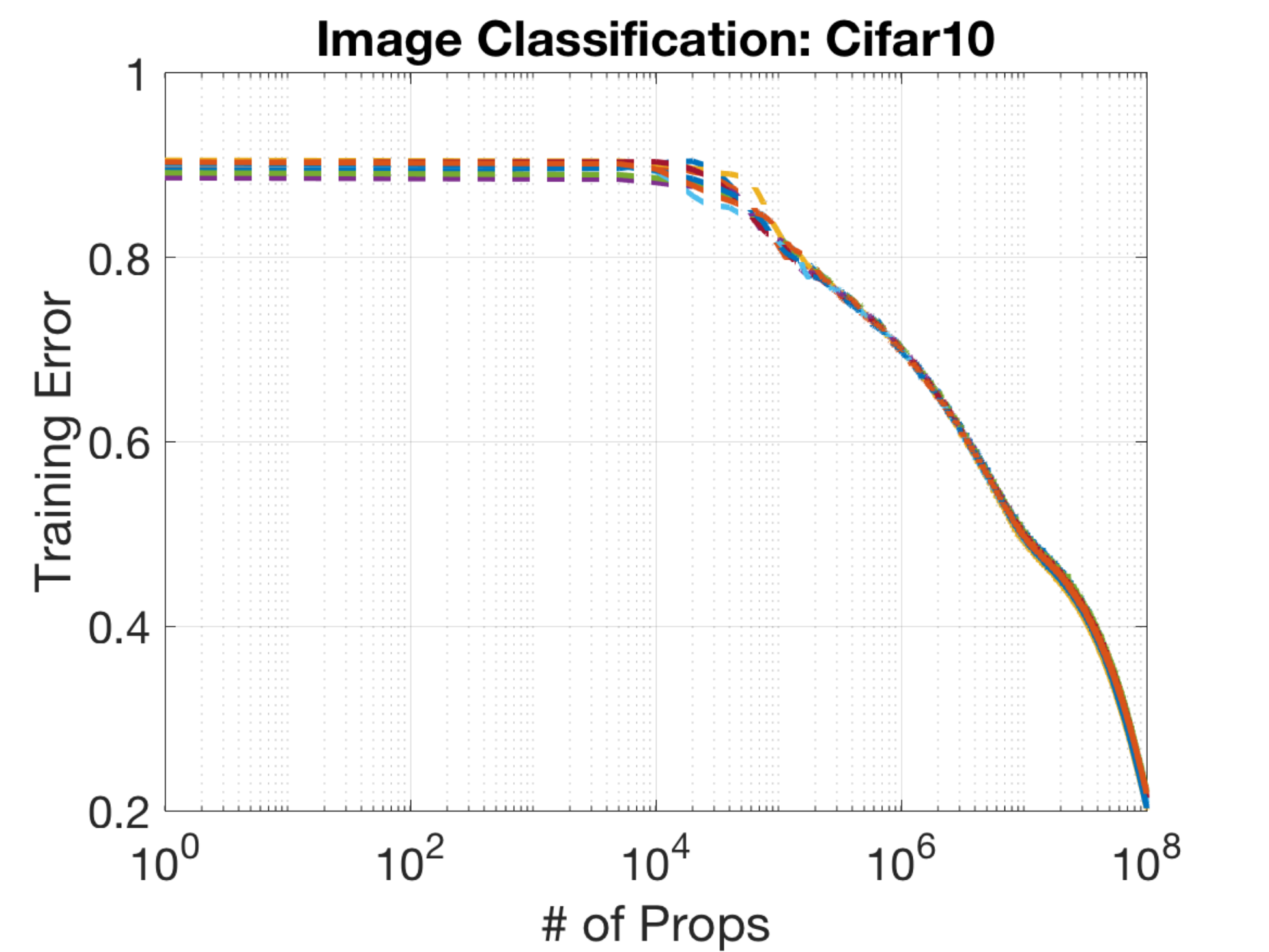}
		\includegraphics[width=0.28\textwidth]{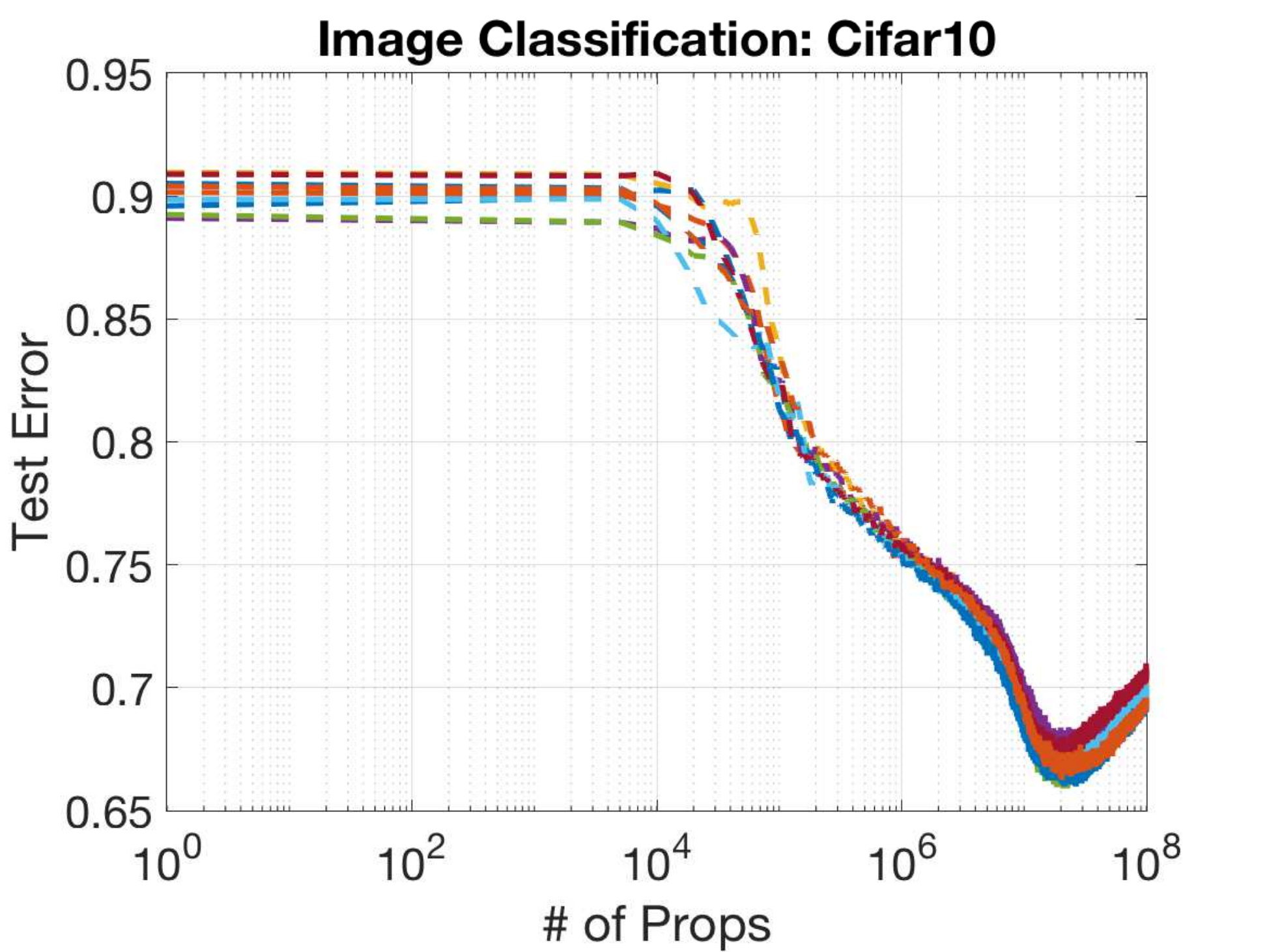}
	}
	\subfigure[9 runs of Sub-sampled TR]{
		\includegraphics[width=0.27\textwidth]{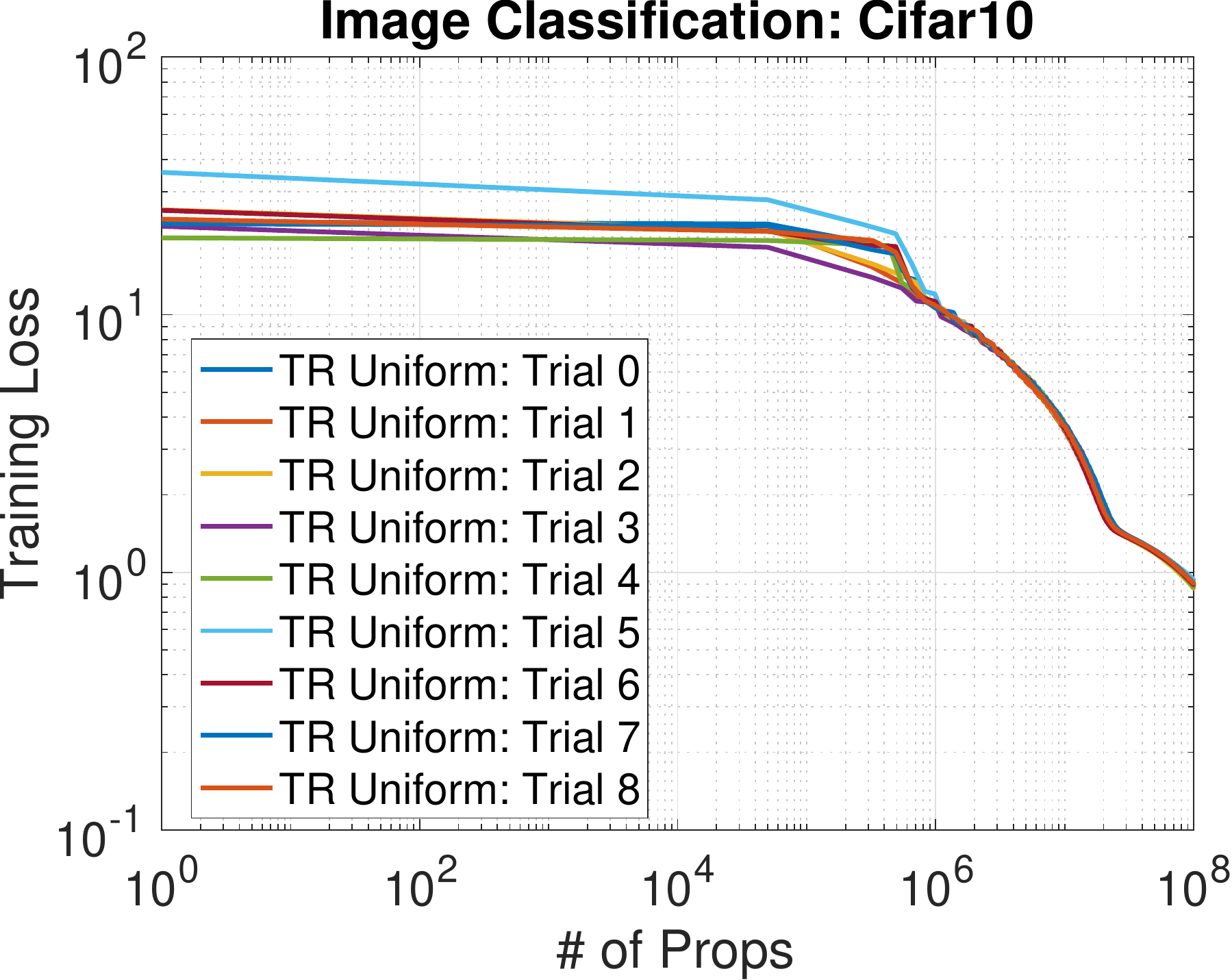}
		\includegraphics[width=0.27\textwidth]{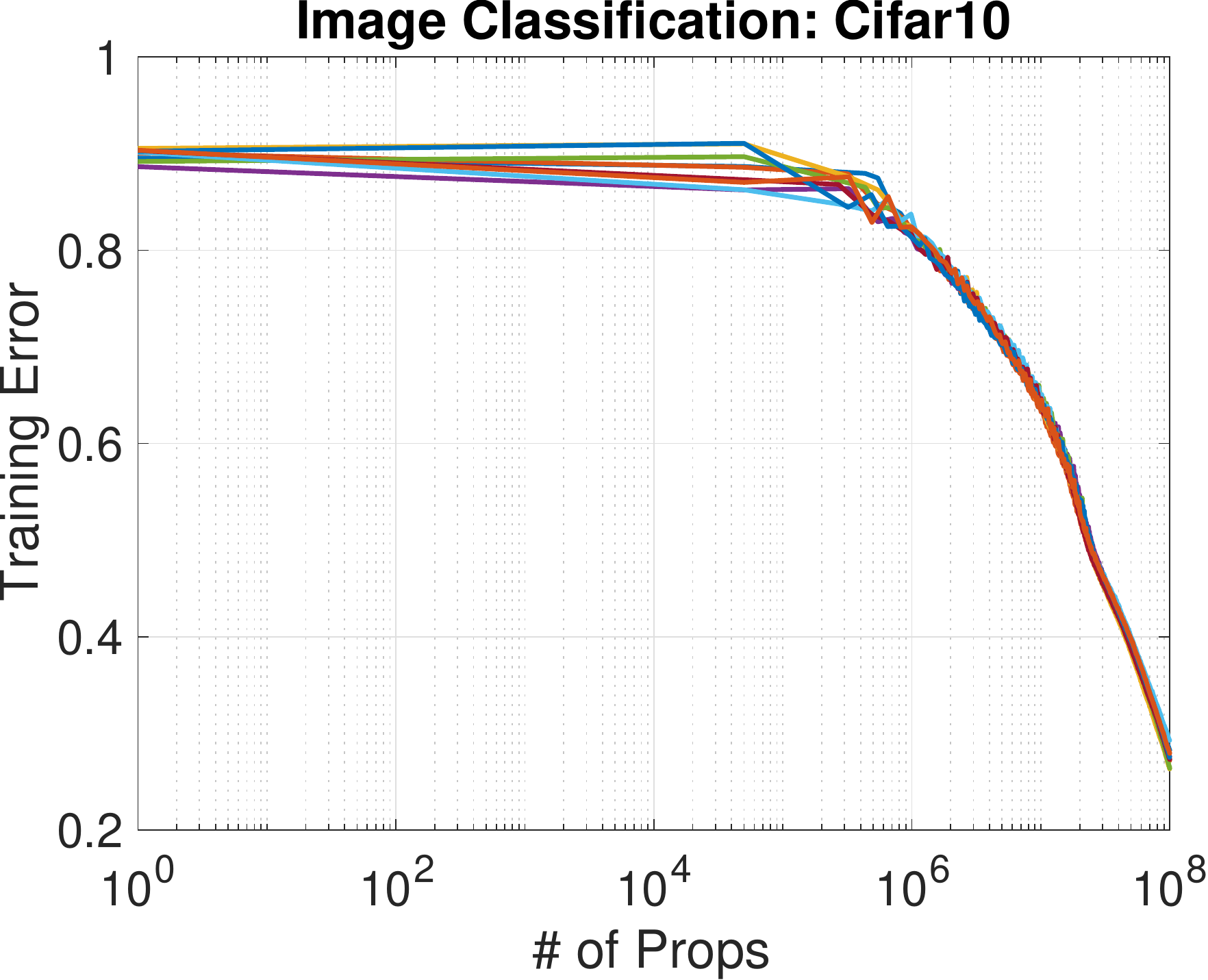}
		\includegraphics[width=0.27\textwidth]{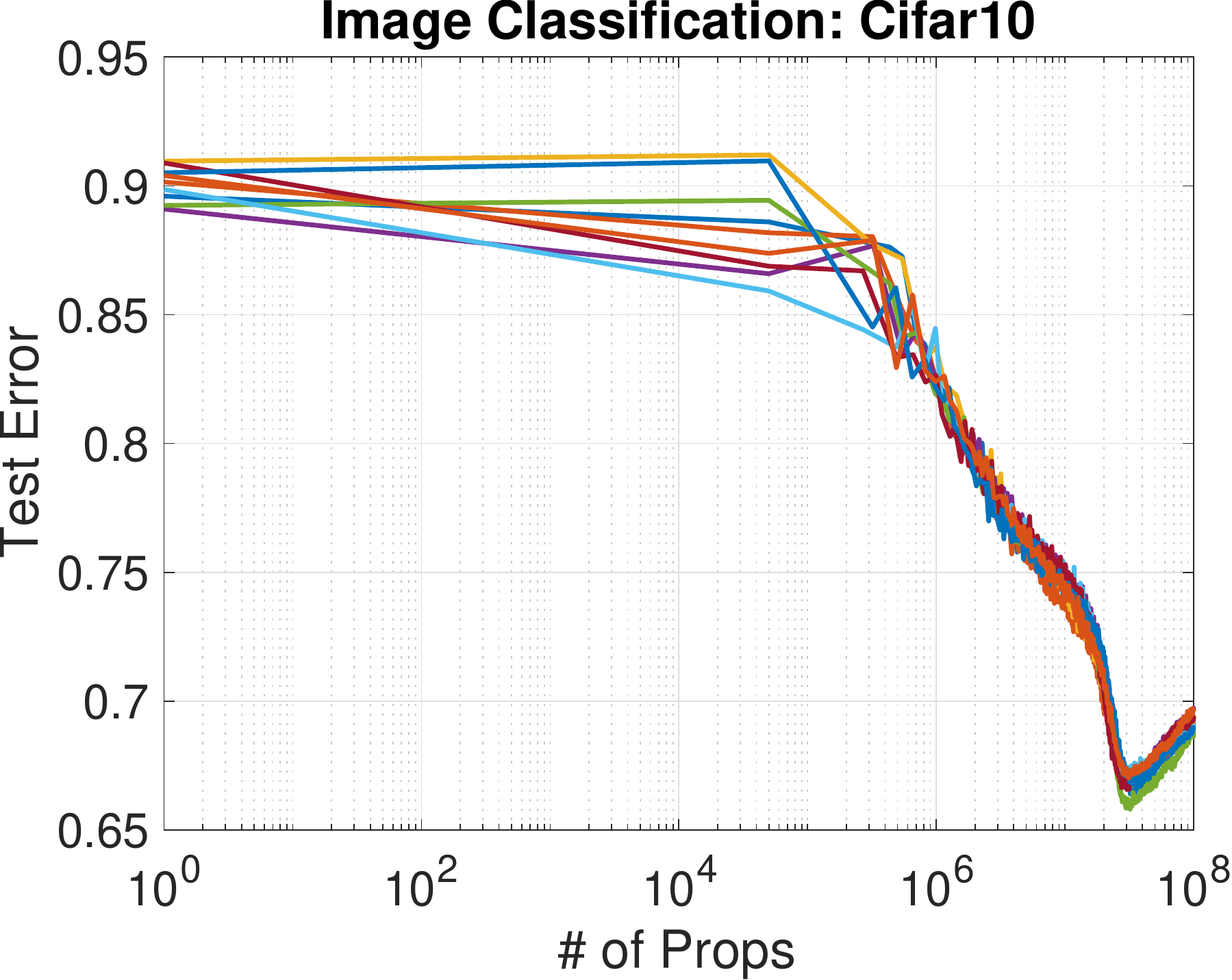}
	}
	\caption{Multiple runs of SGD and Algorithm~\ref{alg:STR_fg} with different random seed. All the runs start with random initialization.}
	\label{fig:cifar4}
\end{figure}

\section{Non-Linear Least Squares}
In this section, we provide more results on different datasets as the experiment in \cref{sec:example_nls}. We use the exact same setup as before and the datasets we considered here are gathered in \cref{tab:data3}.
In this set of experiments, we consider the sub-sampling ratio $1\%$ and $5\%$ of the training data for all sub-sampling methods. 
\begin{table}[htbp]
	\caption{Datasets used in binary linear classification. \texttt{mnist2} is taking even digits as label 1 and odd digits as label 0. \texttt{mnist-28} only consists of digit 2 and 8 from \texttt{mnist}.}
	\label{tab:data3}
	\centering
	\small
	\begin{tabular}{cccc}
		\toprule
		\sc Data & Training Size ($n$) & \# Features ($d$) & Test Size  \\ 
		\midrule
		{\tt a9a} & $32,561$ & $123$ & $16,281$
		\\
		{\tt mnist-28} & $13,007$  & $784$ & $2,163$
		\\
		{\tt mnist2} & $60,000$& $784$ &$10,000$
		\\
		{\tt ijcnn1} & $49,990$ & $22$ & $16,281$ 
		\\
		\bottomrule
	\end{tabular}
\end{table}
\begin{figure}[htbp]
\centering
\includegraphics[width=0.6\textwidth]{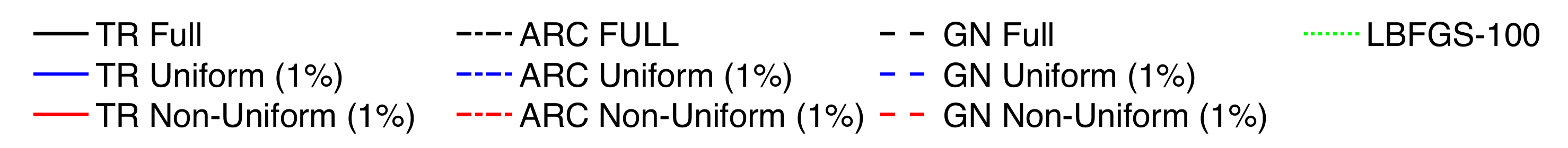}\vspace{-2mm}
\subfigure[use $1\%$ of training data]{
\includegraphics[width=0.24\textwidth]{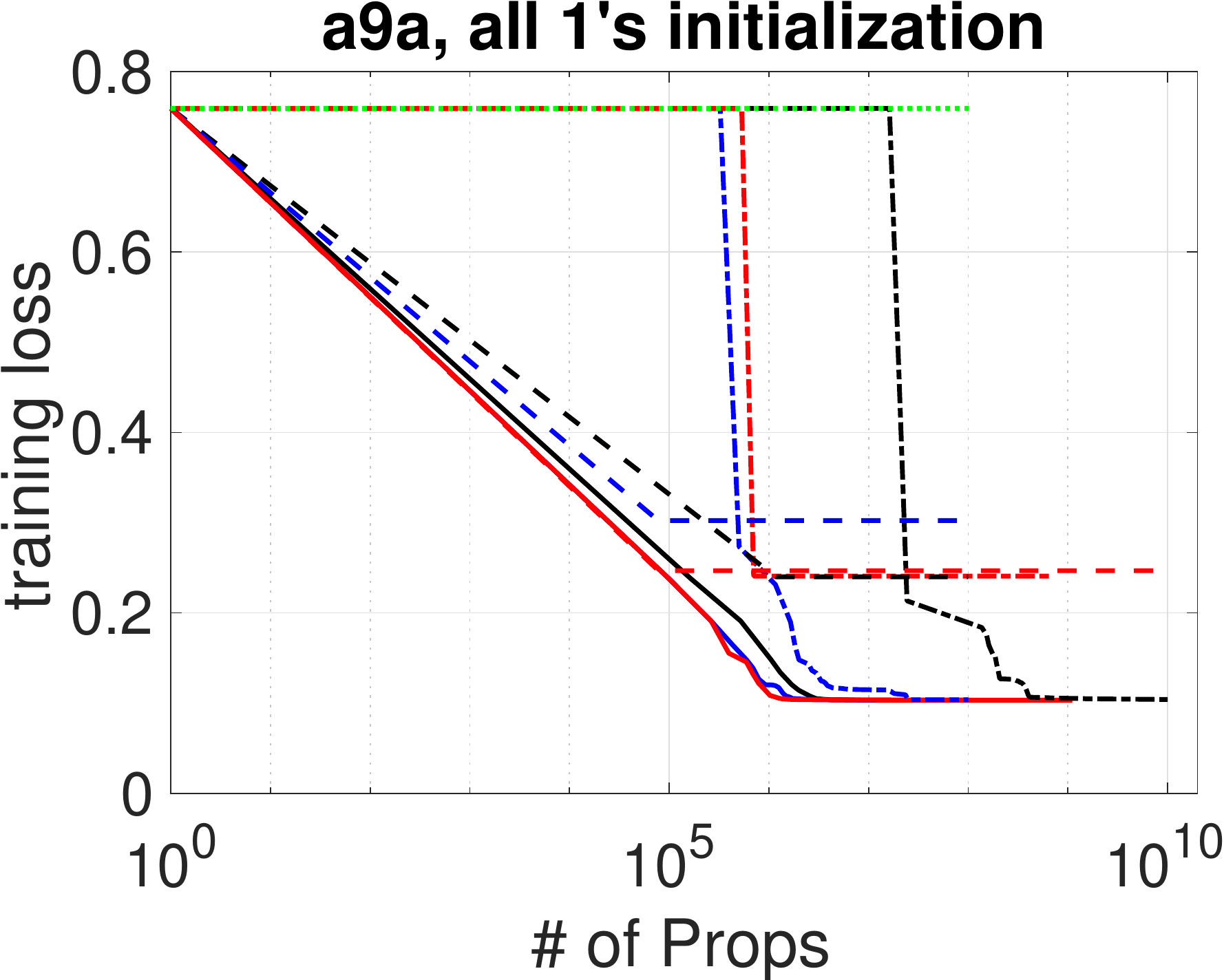}
\includegraphics[width=0.24\textwidth]{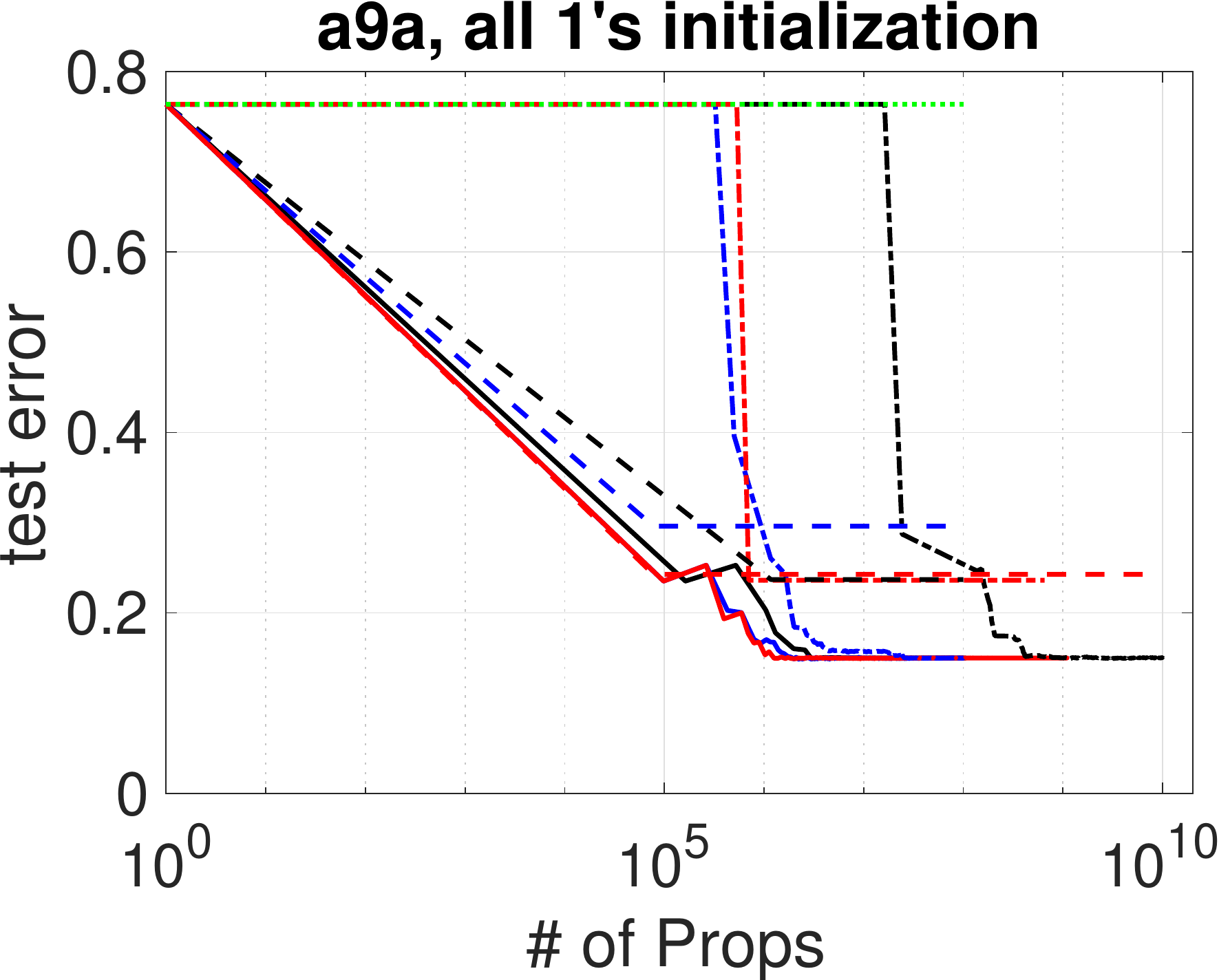}
\includegraphics[width=0.24\textwidth]{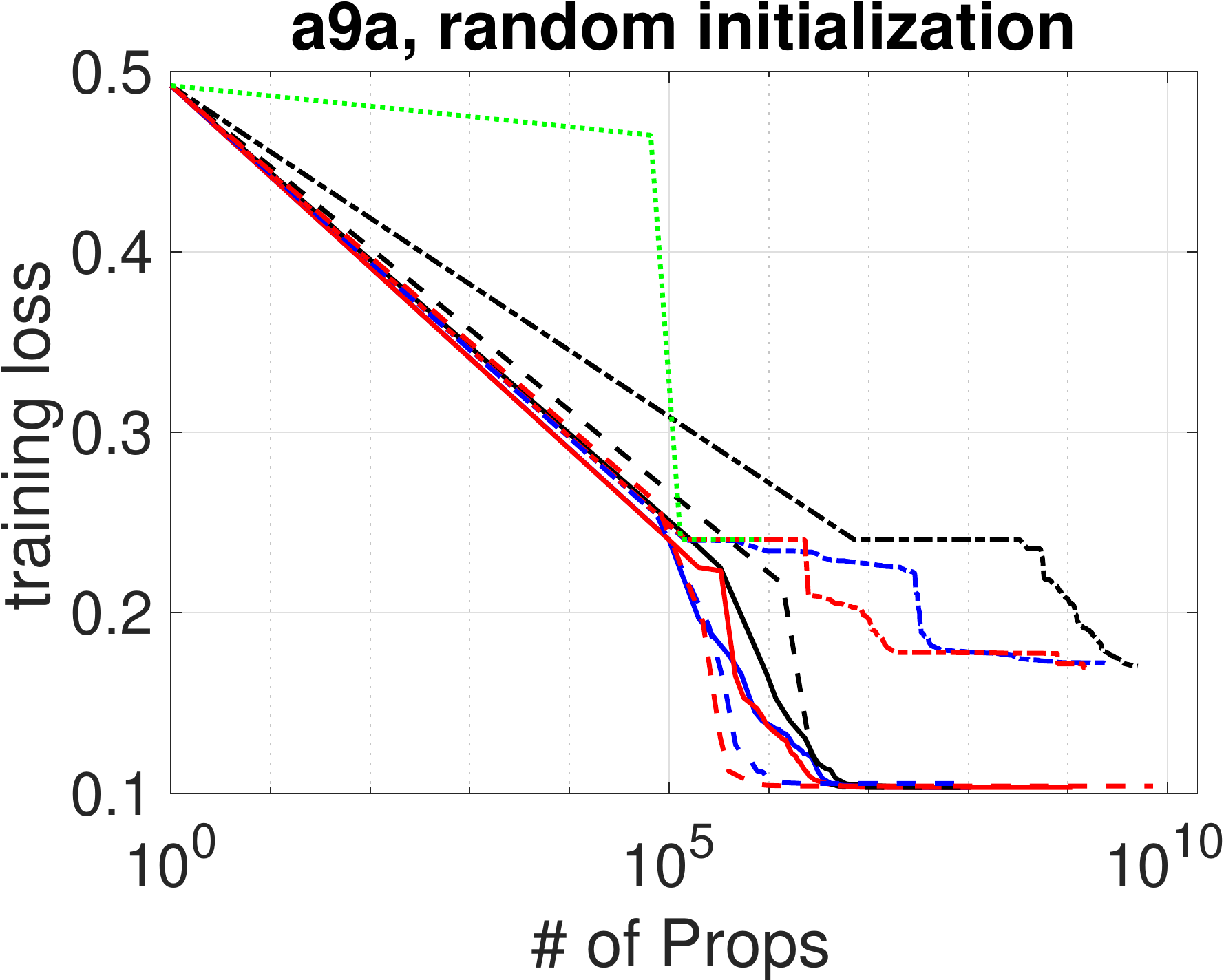}
\includegraphics[width=0.24\textwidth]{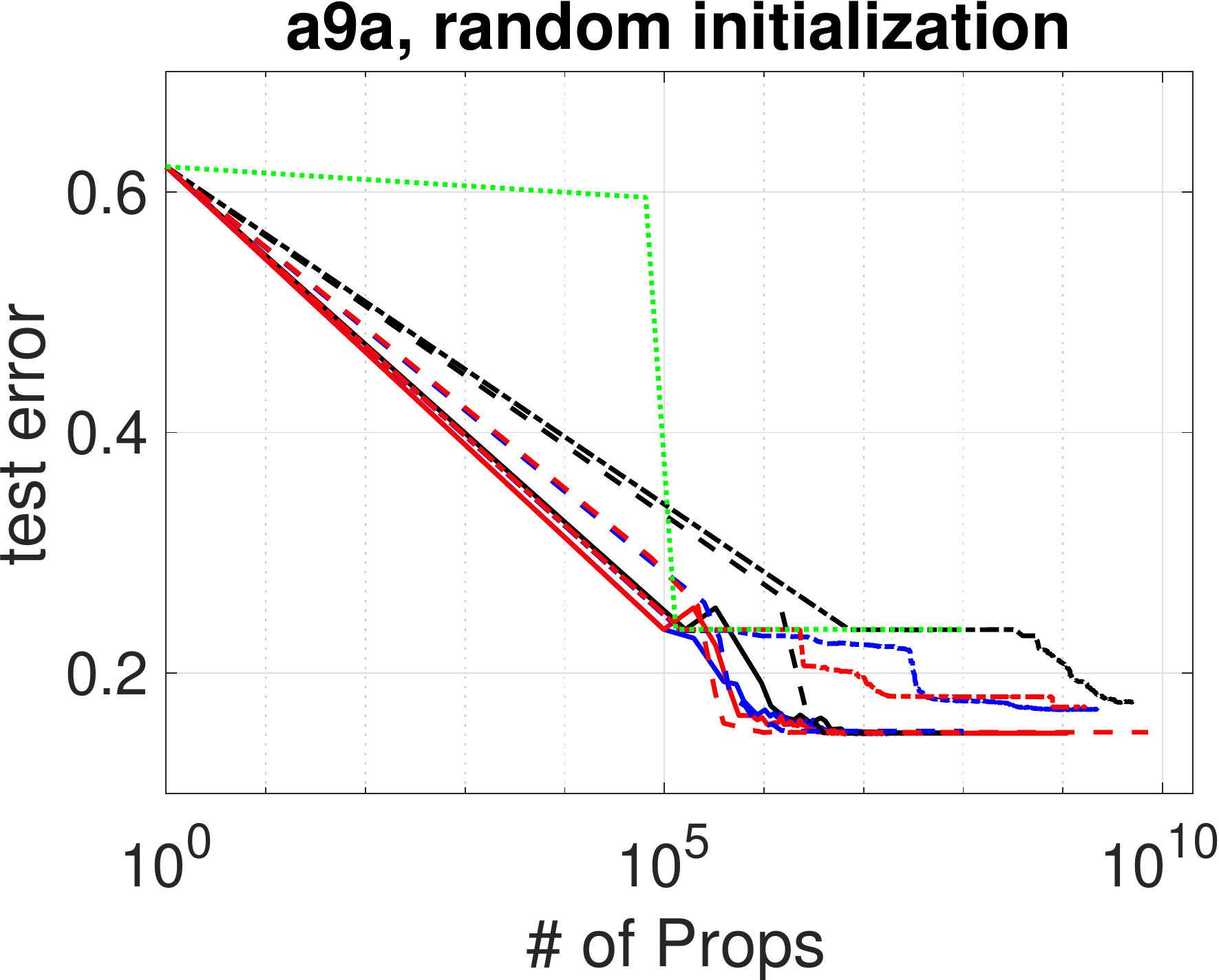}
}
\label{fig:a9a-1}
\includegraphics[width=0.6\textwidth]{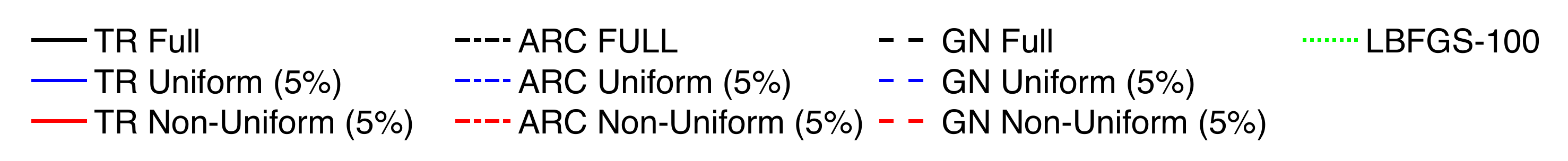}\vspace{-2mm}
\subfigure[use $5\%$ of training data]{
\includegraphics[width=0.24\textwidth]{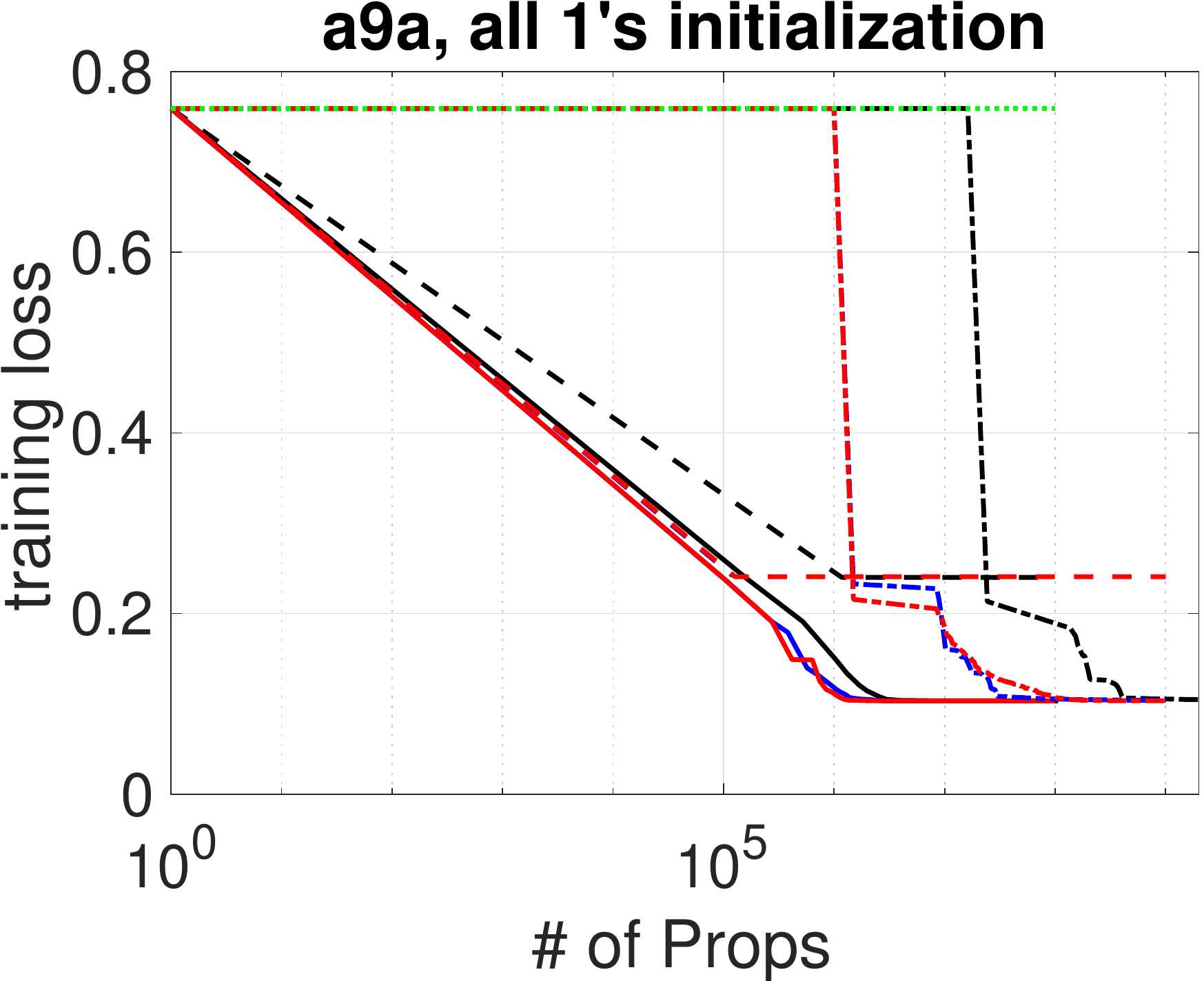}
\includegraphics[width=0.24\textwidth]{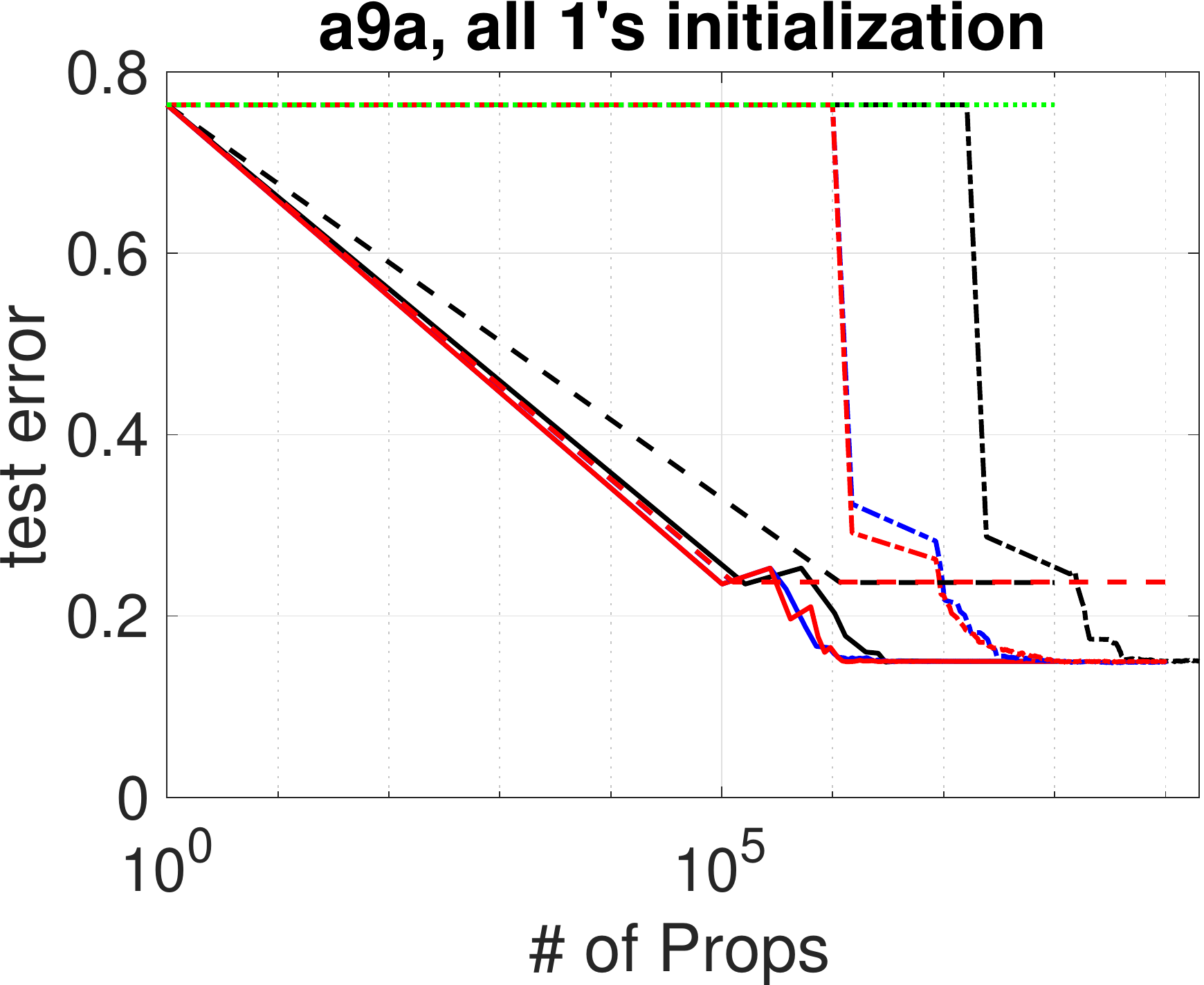}
\includegraphics[width=0.24\textwidth]{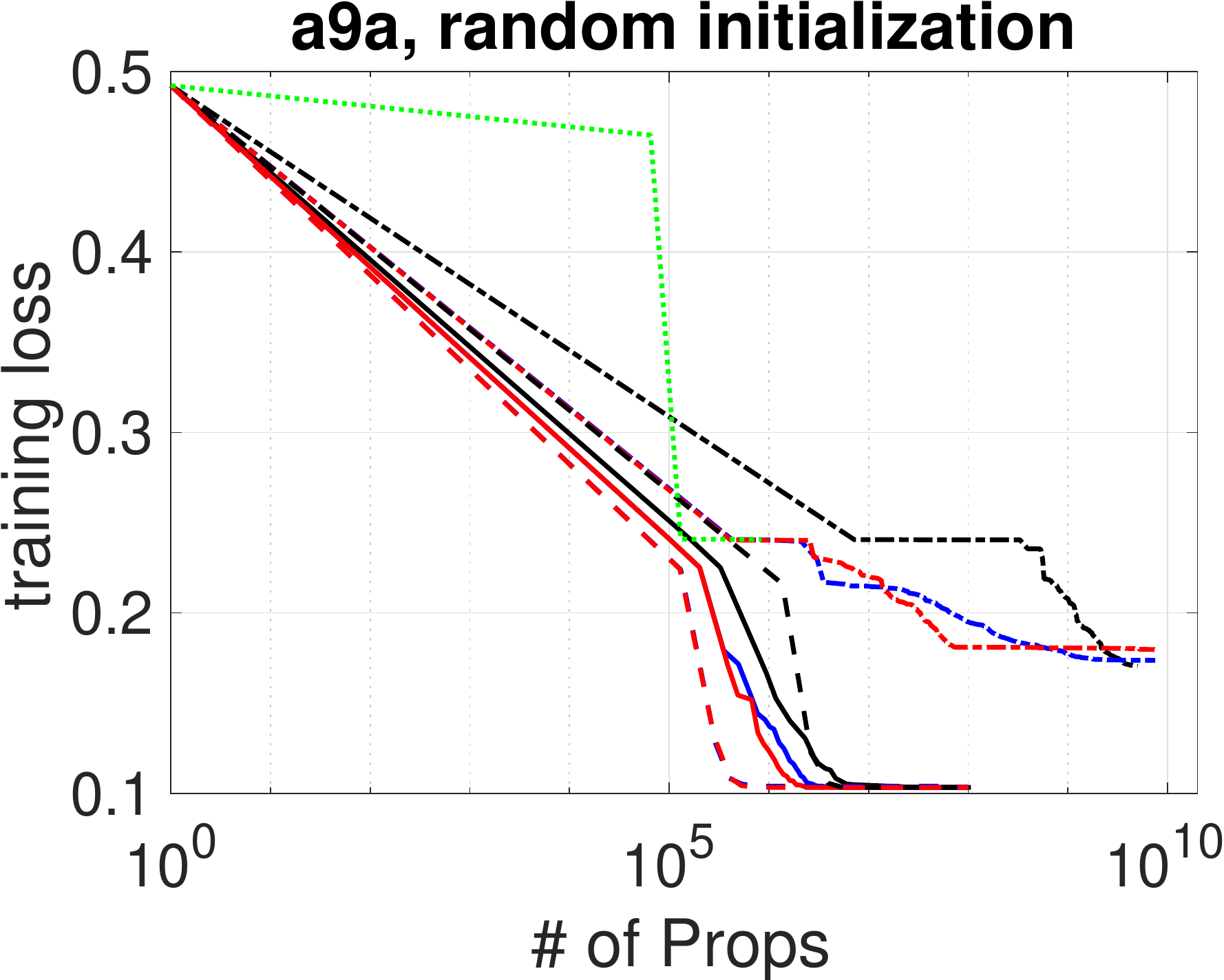}
\includegraphics[width=0.24\textwidth]{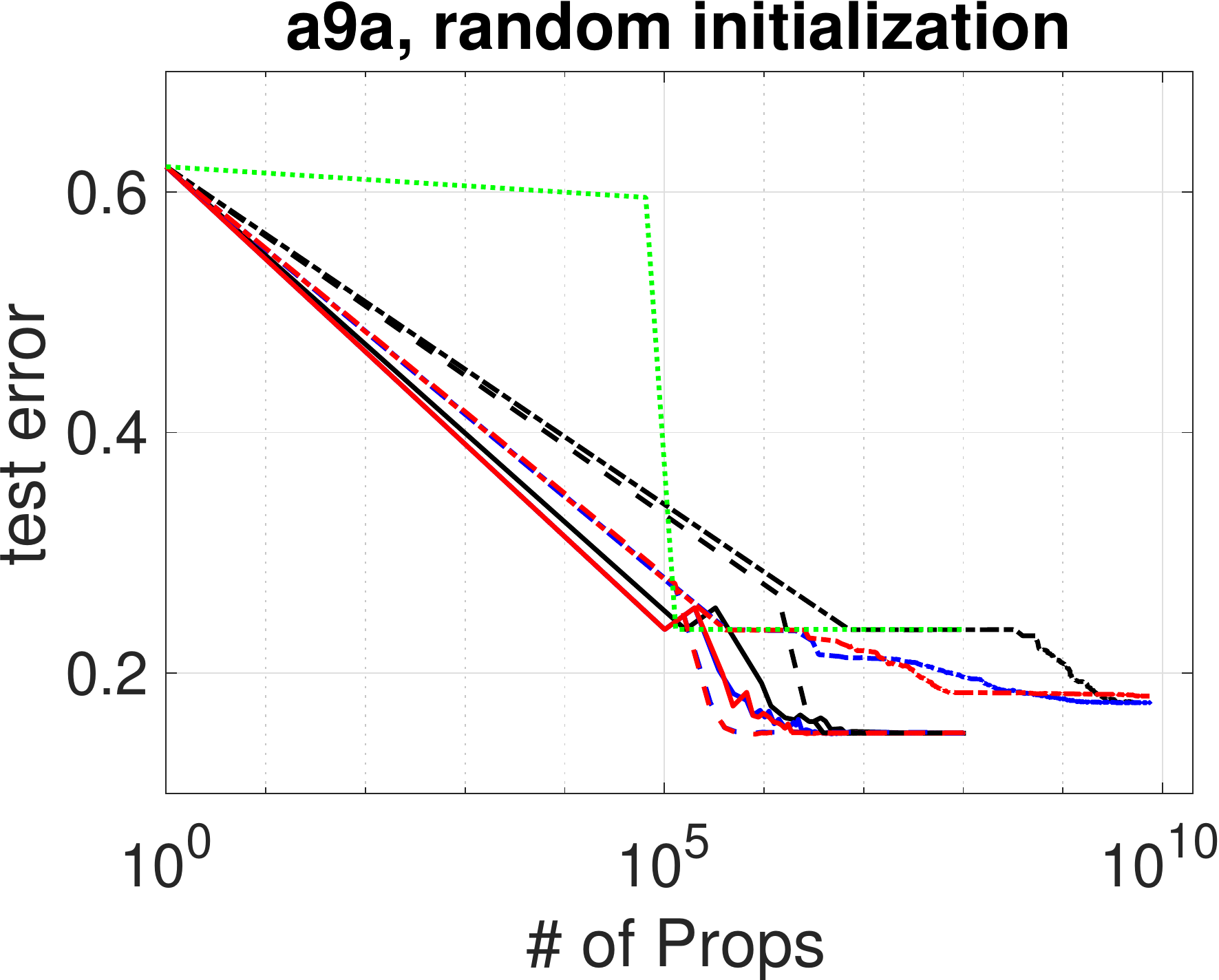}
}
\caption{Binary Classification on \texttt{a9a} with different subsampling sizes. We observe that for both sampling ratios, \emph{TR Uniform} and \emph{TR Non-Uniform} perform very similar in all cases and they converge faster than \emph{TR Full}. Similarly, \emph{ARC Uniform} and \emph{ARC Non-Uniform} outperform then \emph{ARC Full}. Increasing the sample size from $1\%$ to $5\%$ does not seem to help, suggesting that $1\%$ might in fact be a sufficient sample size. Further, GN performs well with random initialization (unlike all 1's vector), which is consistent with our previous observation.}
\label{fig:a9a}
\end{figure}

\begin{figure}[htbp]
\centering
\includegraphics[width=0.6\textwidth]{figs/blc_legend_1}\vspace{-2mm}
\subfigure[use $1\%$ of training data]{
\includegraphics[width=0.24\textwidth]{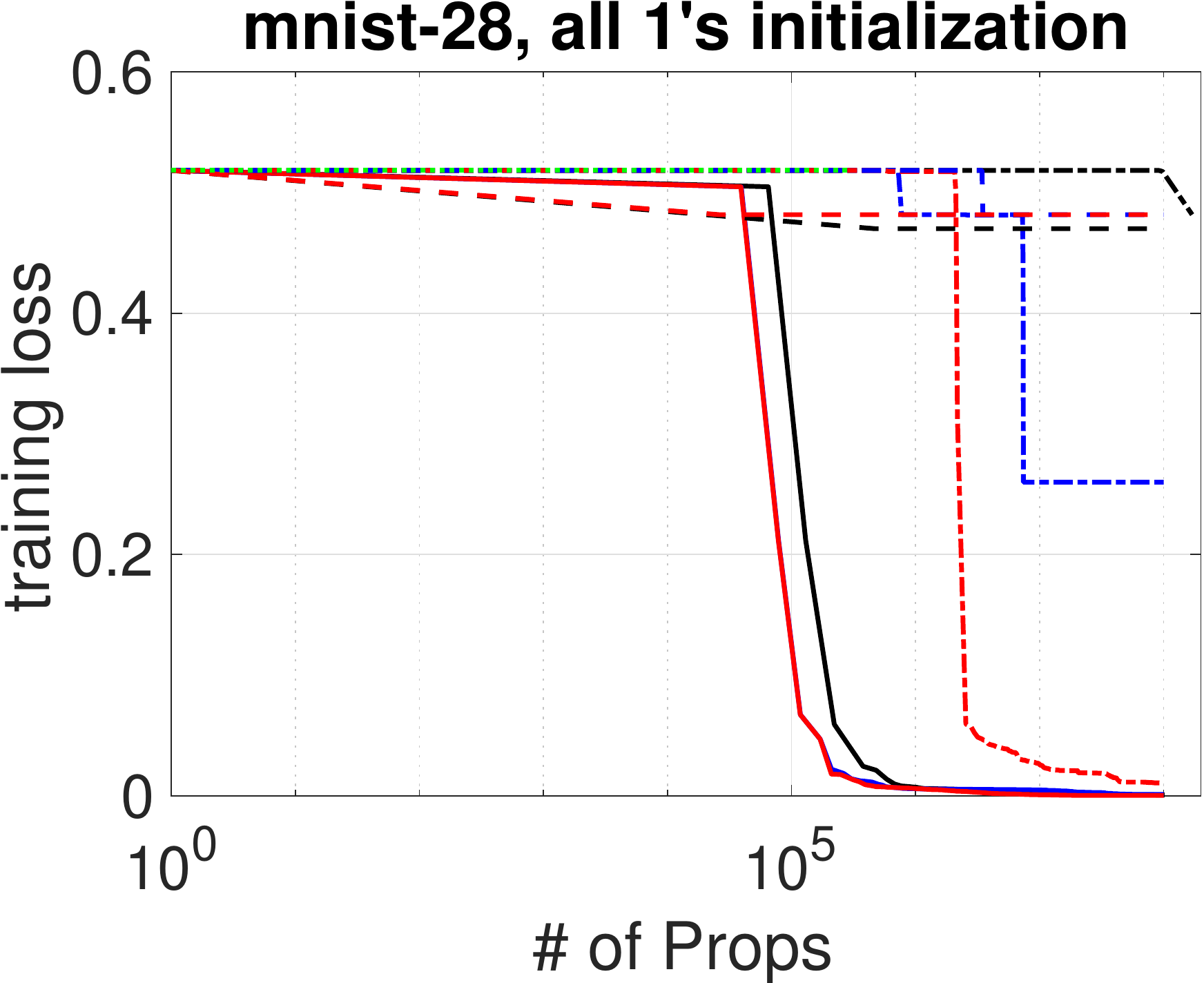}
\includegraphics[width=0.24\textwidth]{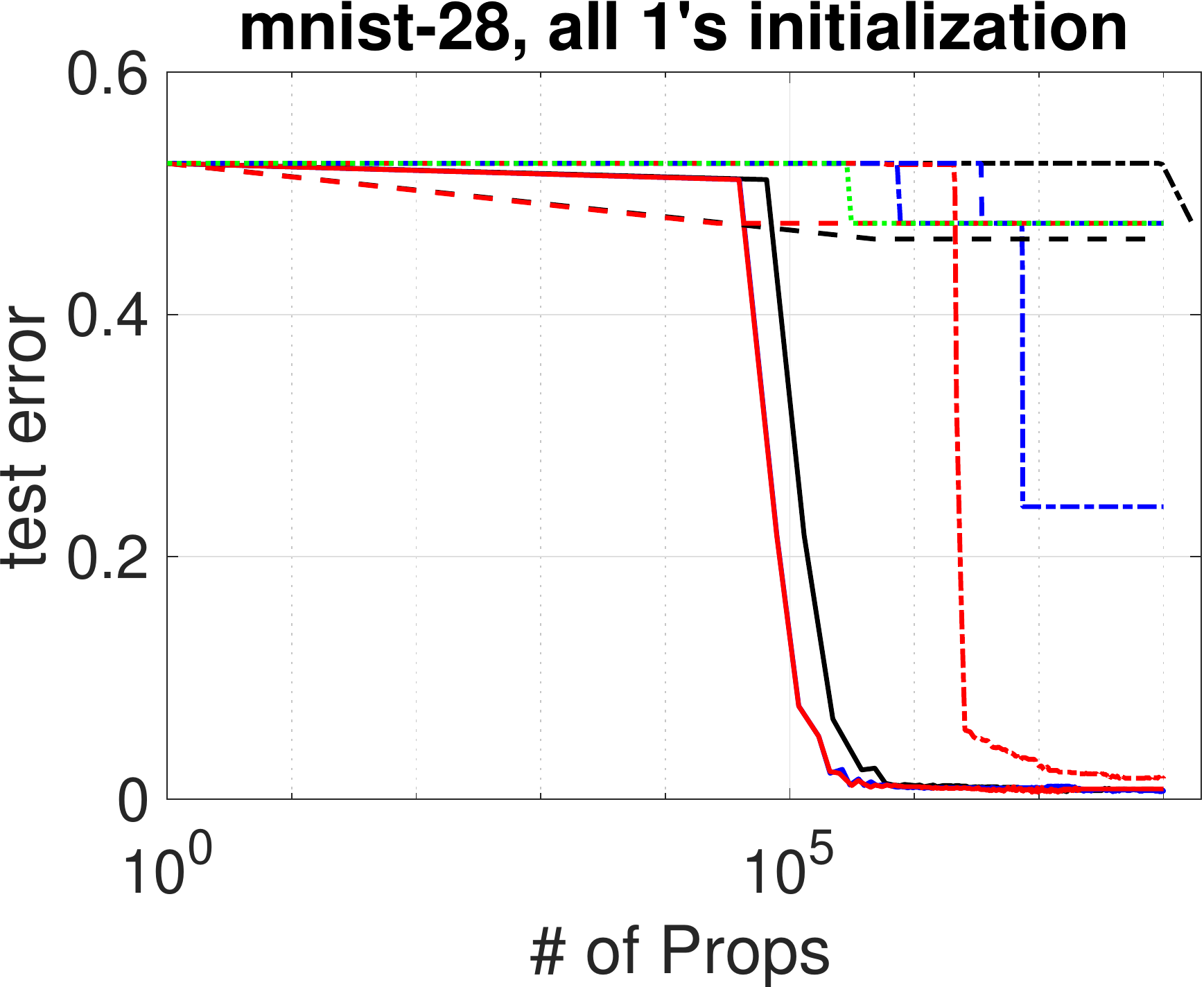}
\includegraphics[width=0.24\textwidth]{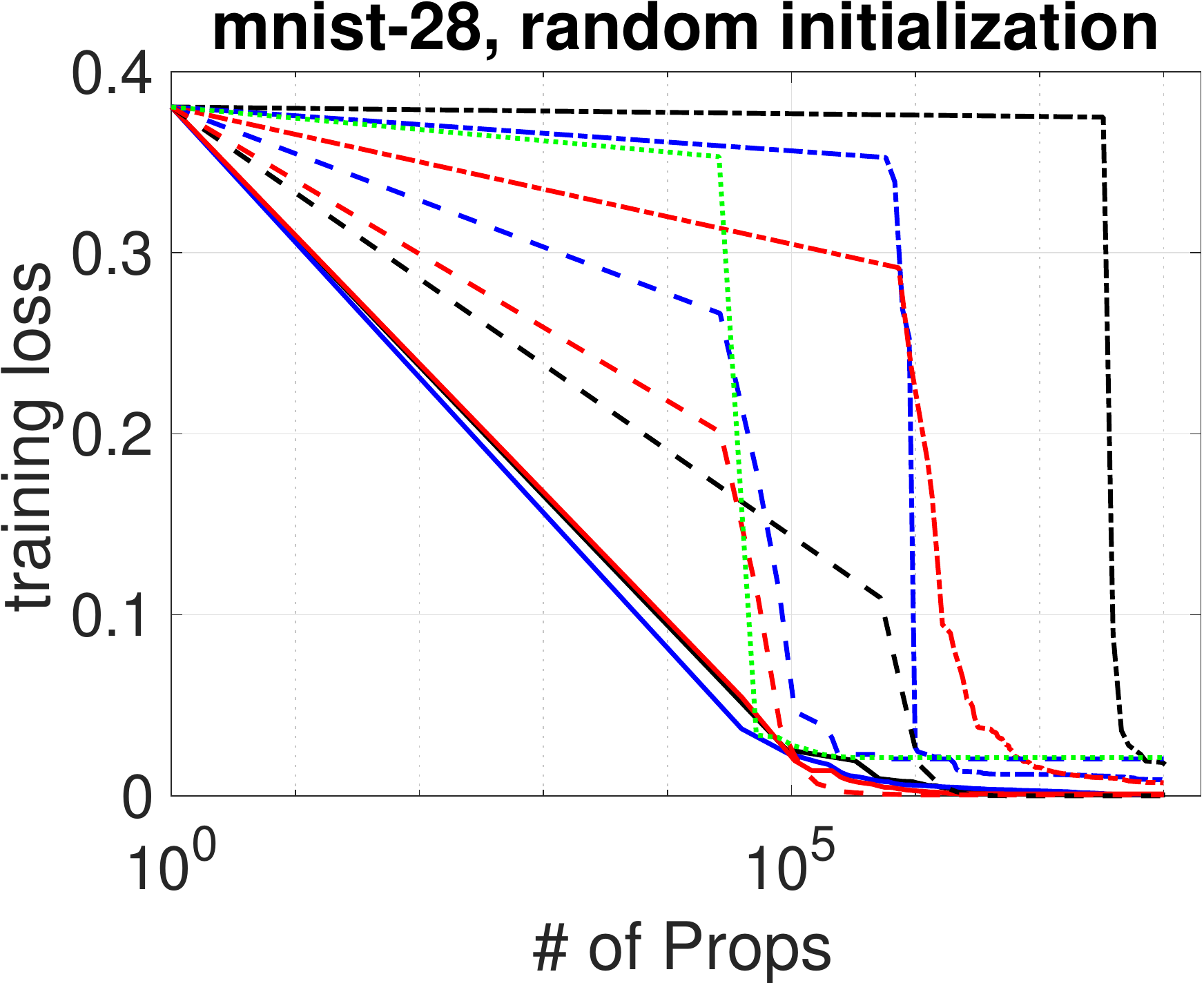}
\includegraphics[width=0.24\textwidth]{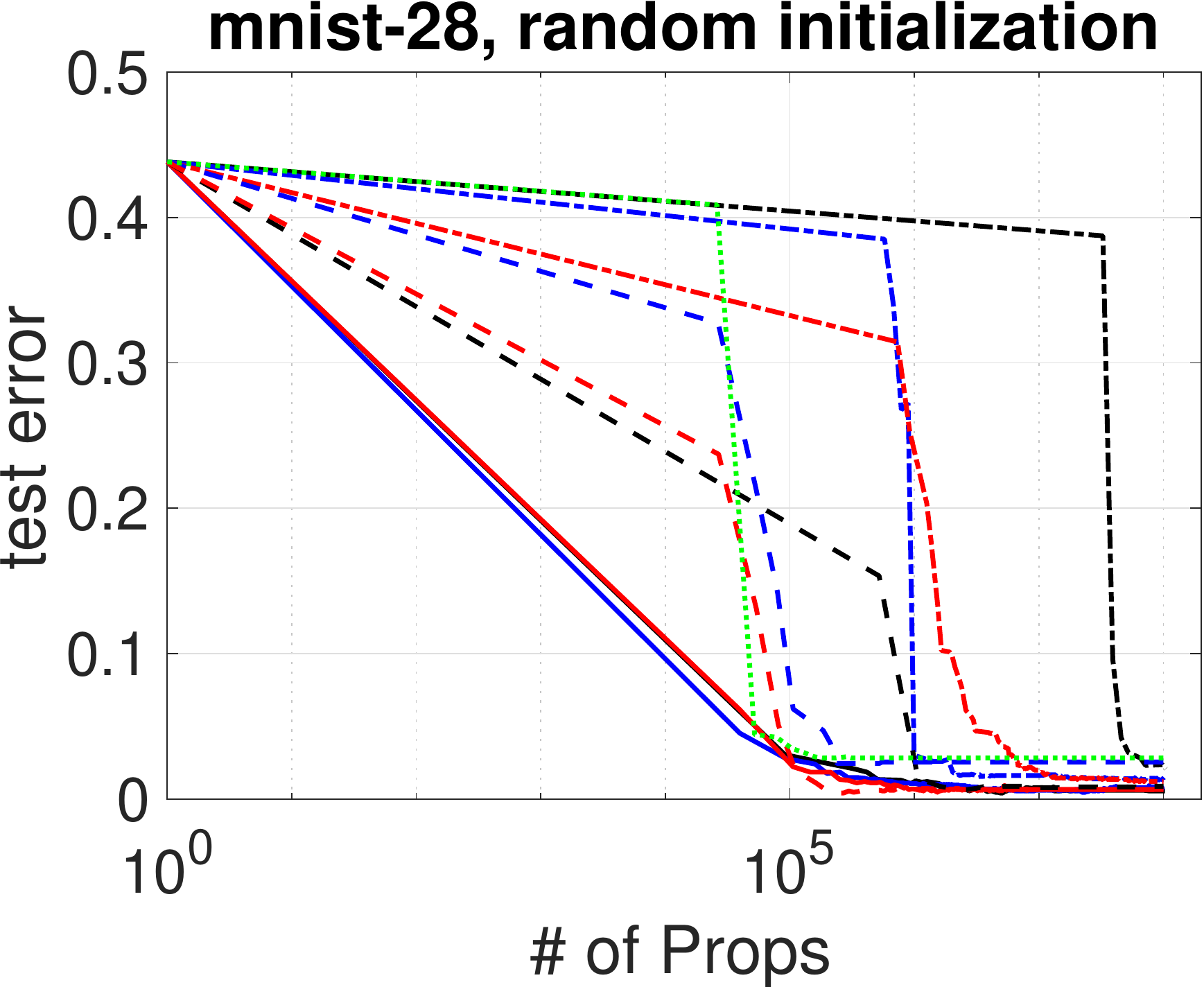}
}
\label{fig:mnist-28-1}
\includegraphics[width=0.6\textwidth]{figs/blc_legend_5}\vspace{-2mm}
\subfigure[use $5\%$ of training data]{
\includegraphics[width=0.24\textwidth]{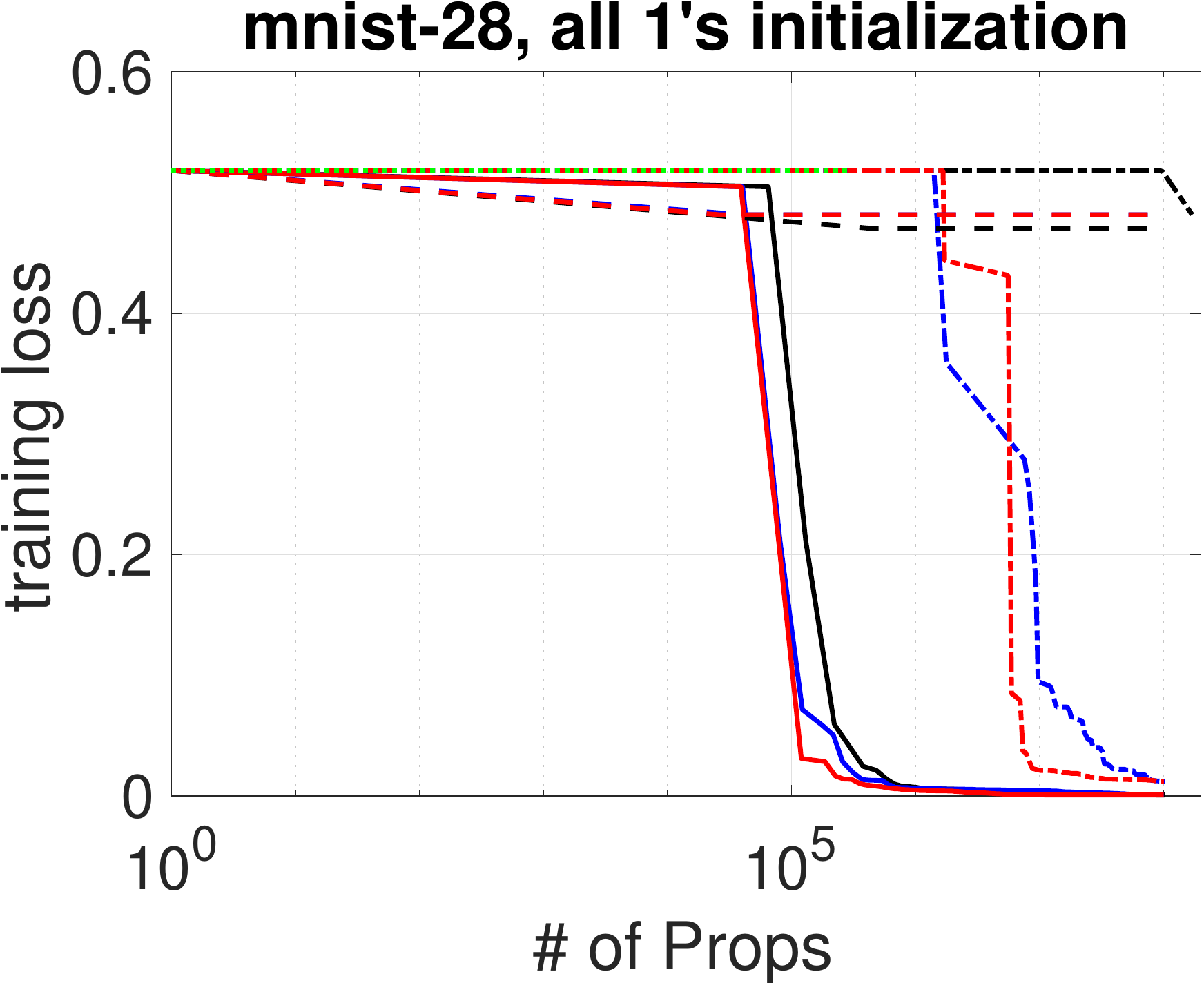}
\includegraphics[width=0.24\textwidth]{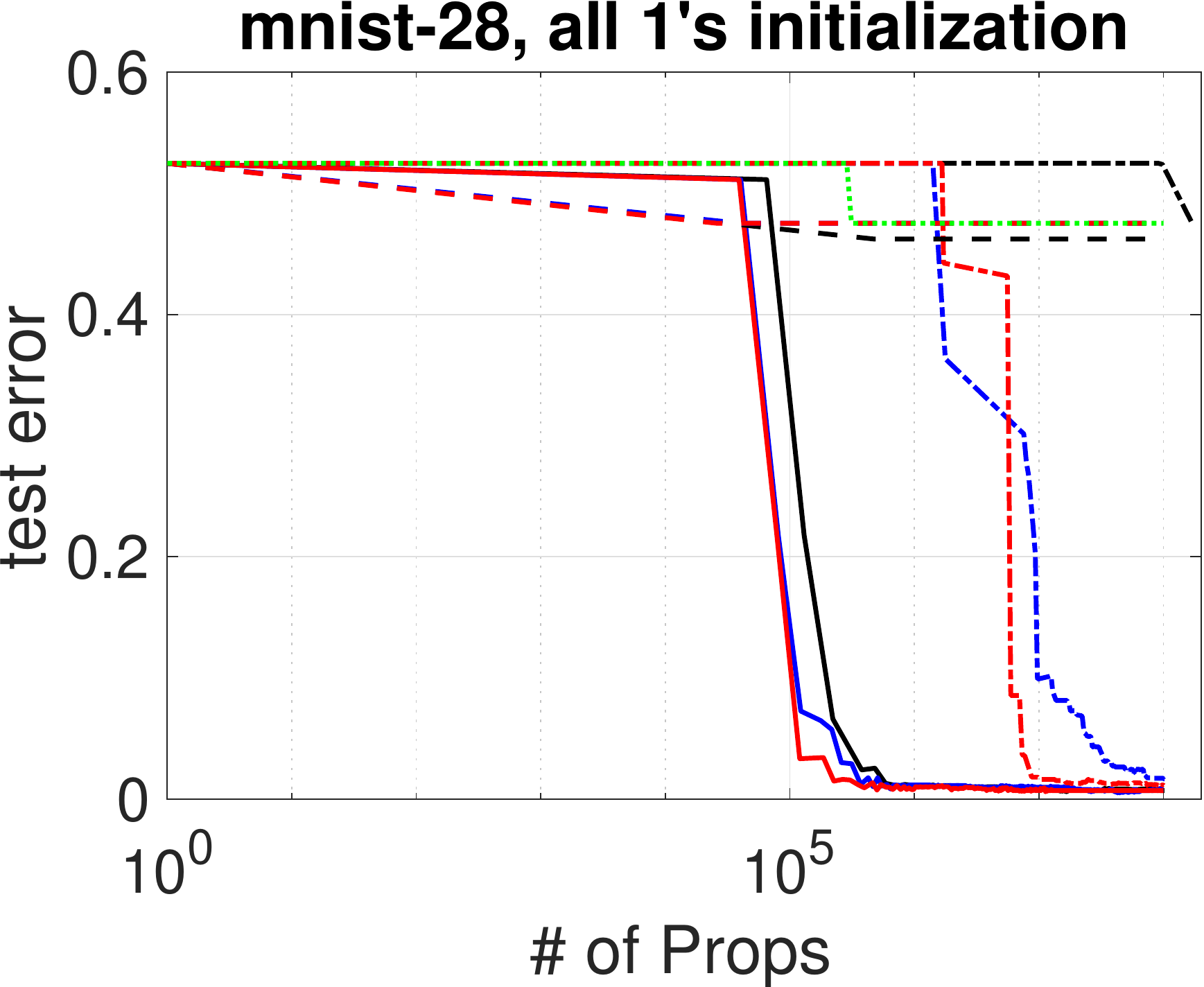}
\includegraphics[width=0.24\textwidth]{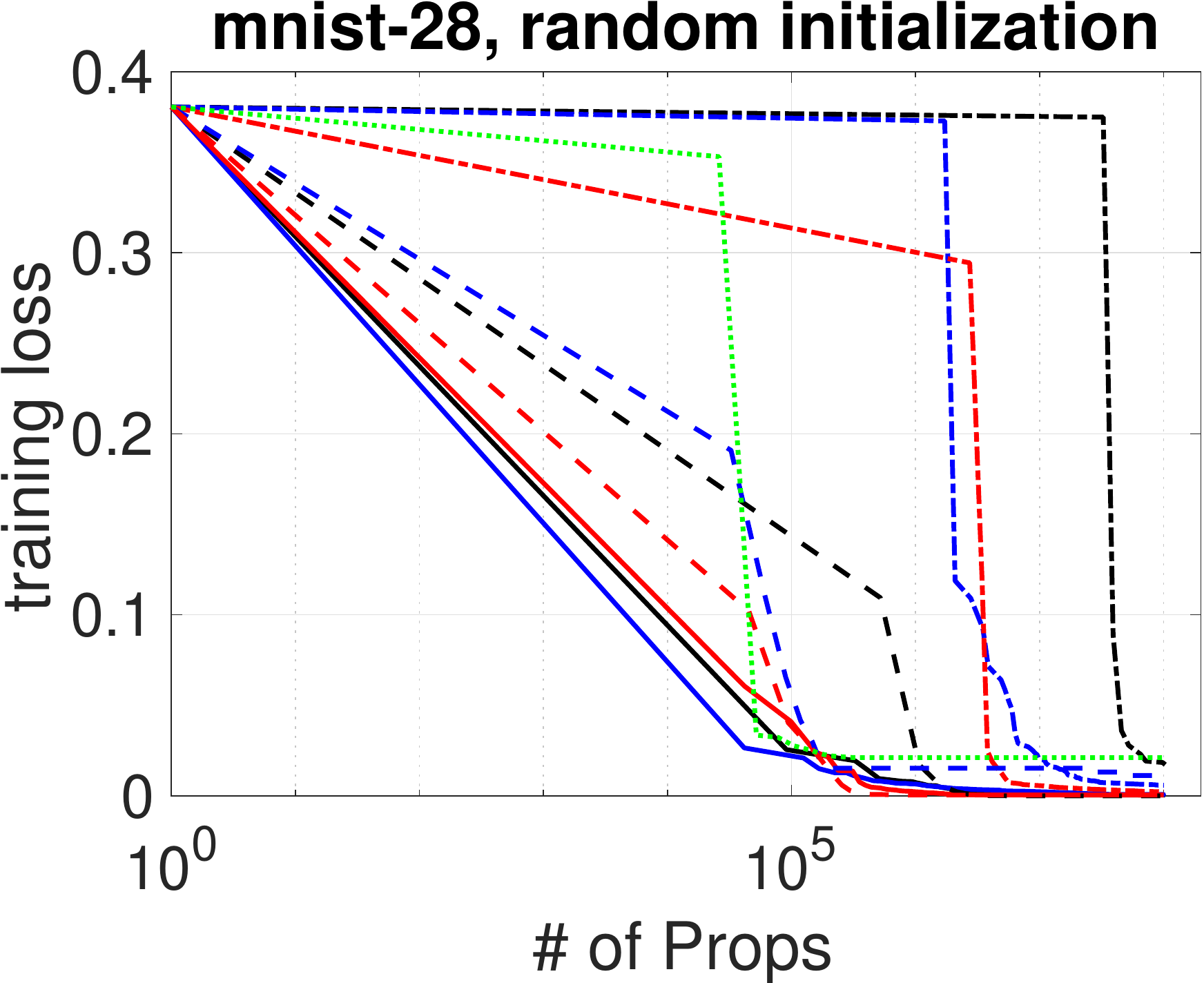}
\includegraphics[width=0.24\textwidth]{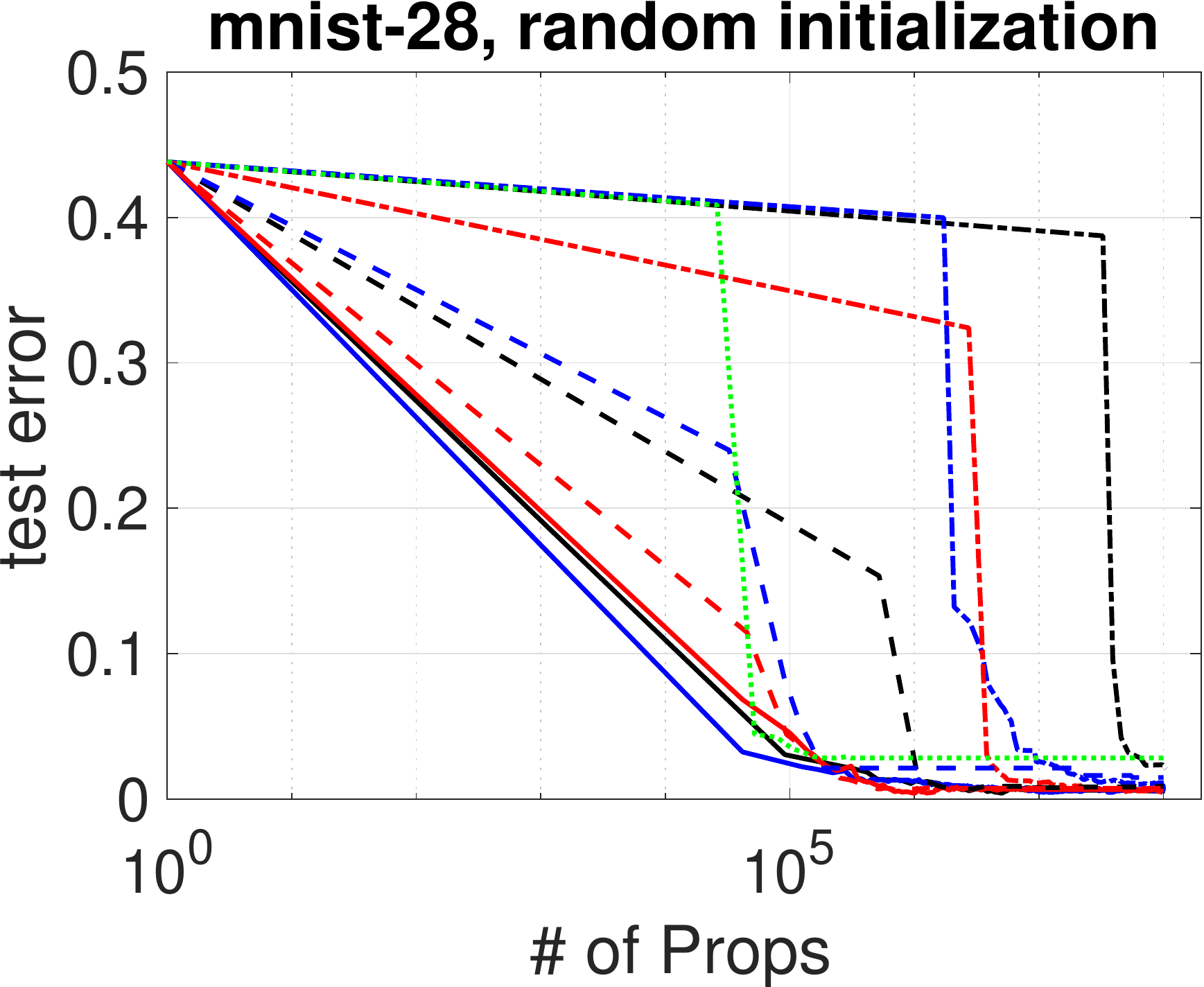}
}
\caption{Binary Classification on \texttt{mnist-28}(digit 2 and 8) with different sub-sampling sizes. We clearly see the advantages of non-uniform sampling over uniform sampling. In all cases here, sub-sampled methods with non-uniform sampling are at least as fast, if not faster,  than the corresponding ones with uniform sampling. Note that increasing the sample size from $1\%$ to $5\%$ does greatly affect the performance of algorithms with uniform sampling (in particular for ARC variants), whereas those employing non-uniform sampling appear to be unaffected. This suggest that, although sufficient for non-uniform sampling,  $1\%$ sample size is indeed not large enough when sampling is done uniformly.}
\label{fig:mnist-28}
\end{figure}

\begin{figure}[htbp]
\centering
\includegraphics[width=0.6\textwidth]{figs/blc_legend_1}\vspace{-2mm}
\subfigure[use $1\%$ of training data]{
\includegraphics[width=0.24\textwidth]{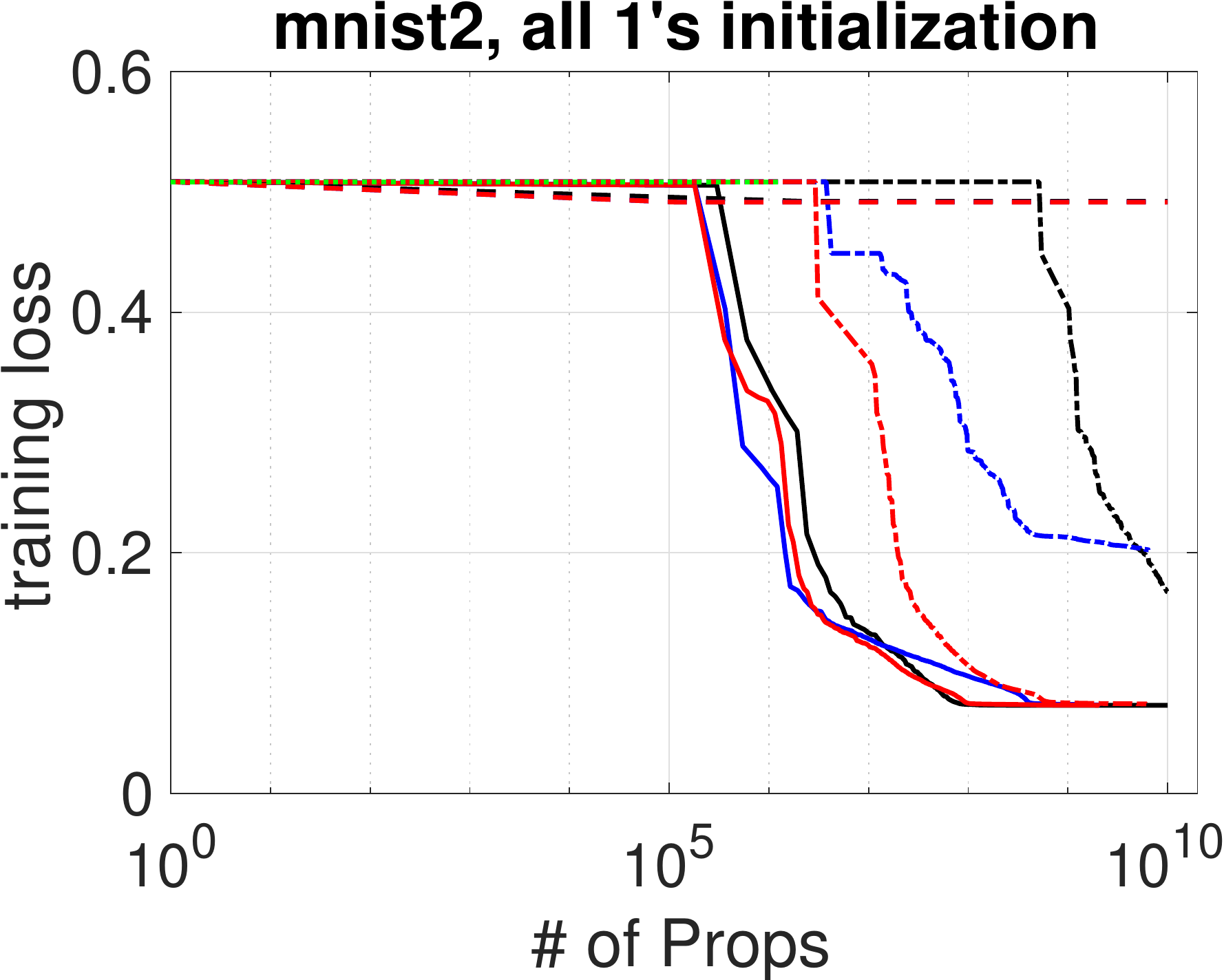}
\includegraphics[width=0.24\textwidth]{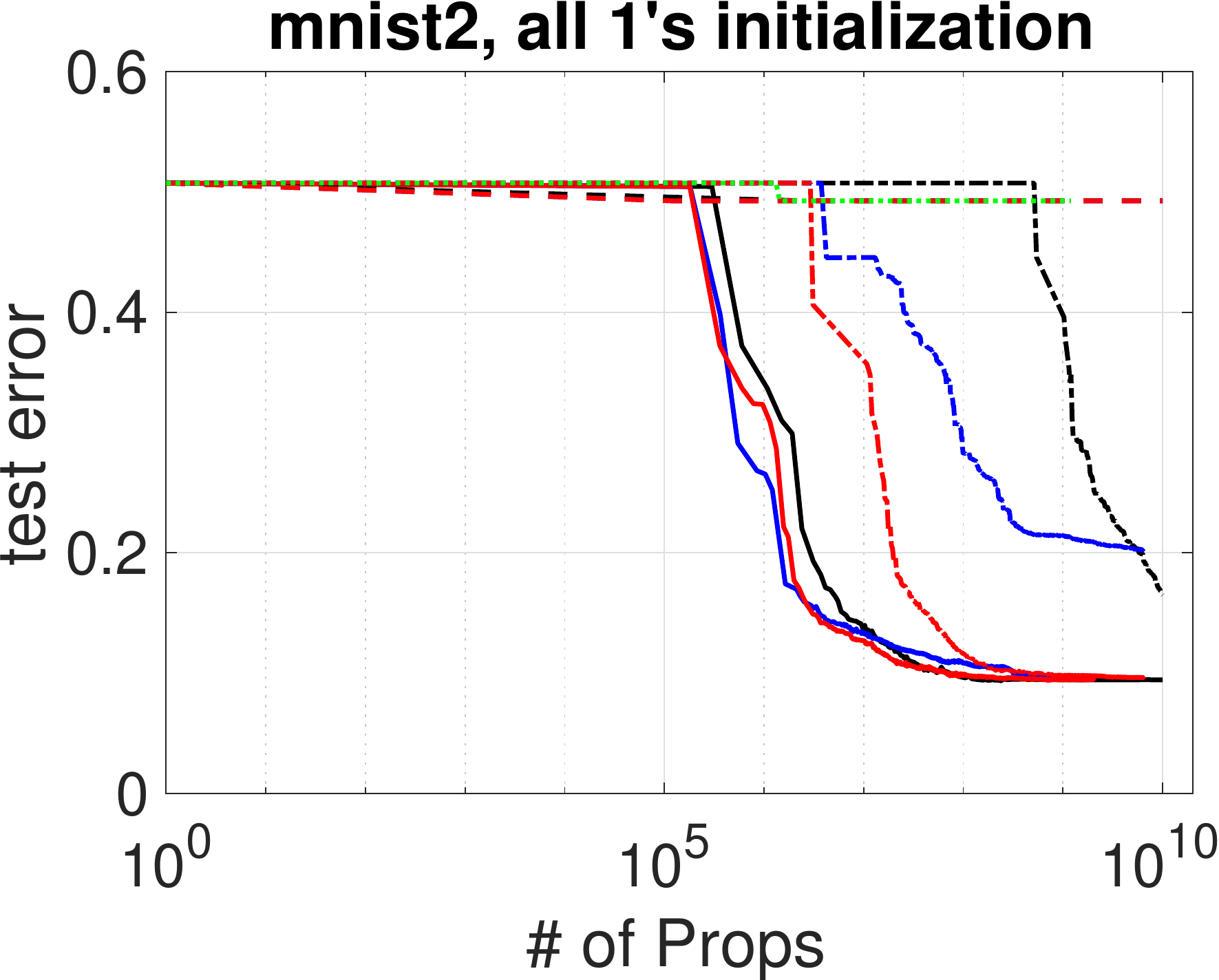}
\includegraphics[width=0.24\textwidth]{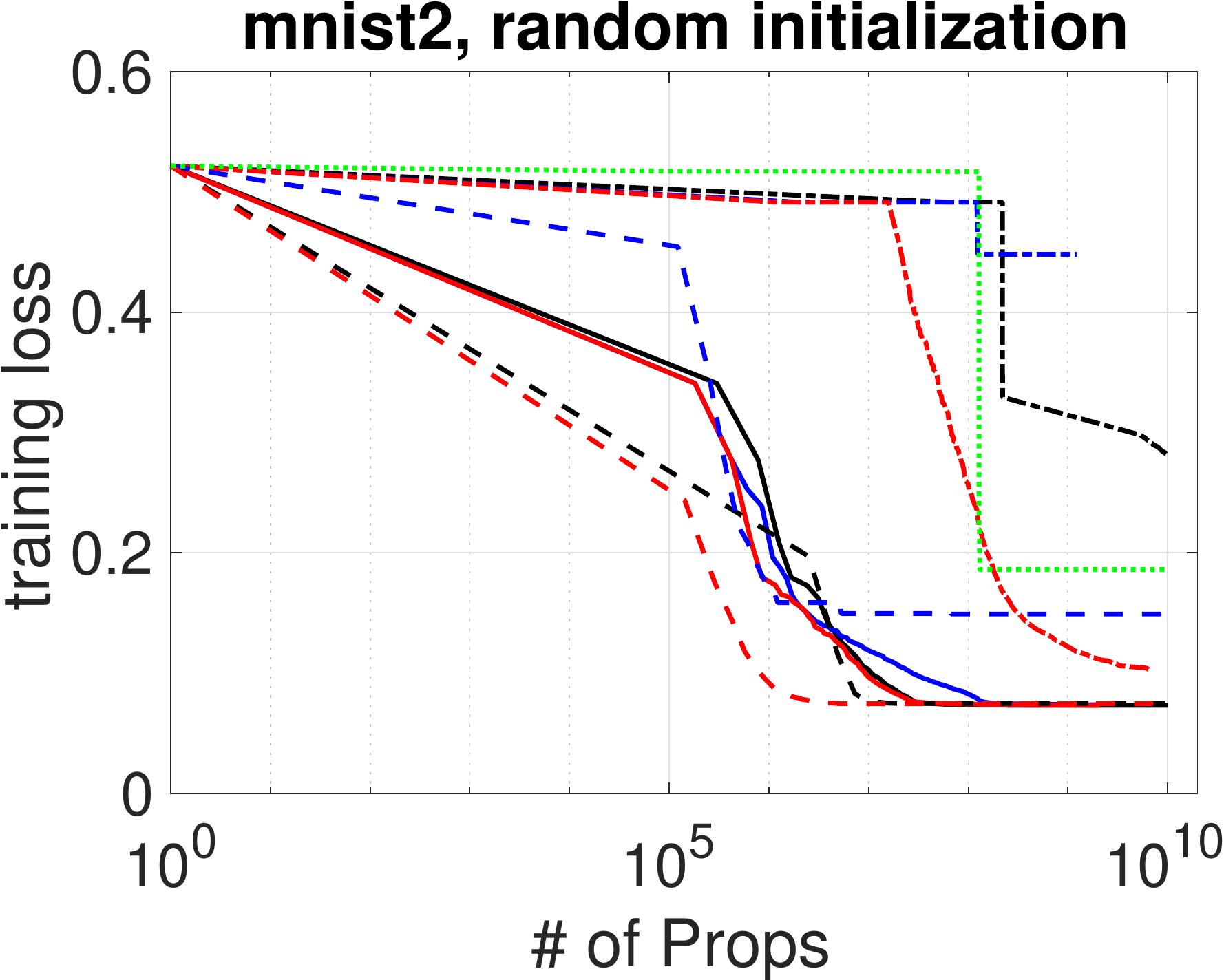}
\includegraphics[width=0.24\textwidth]{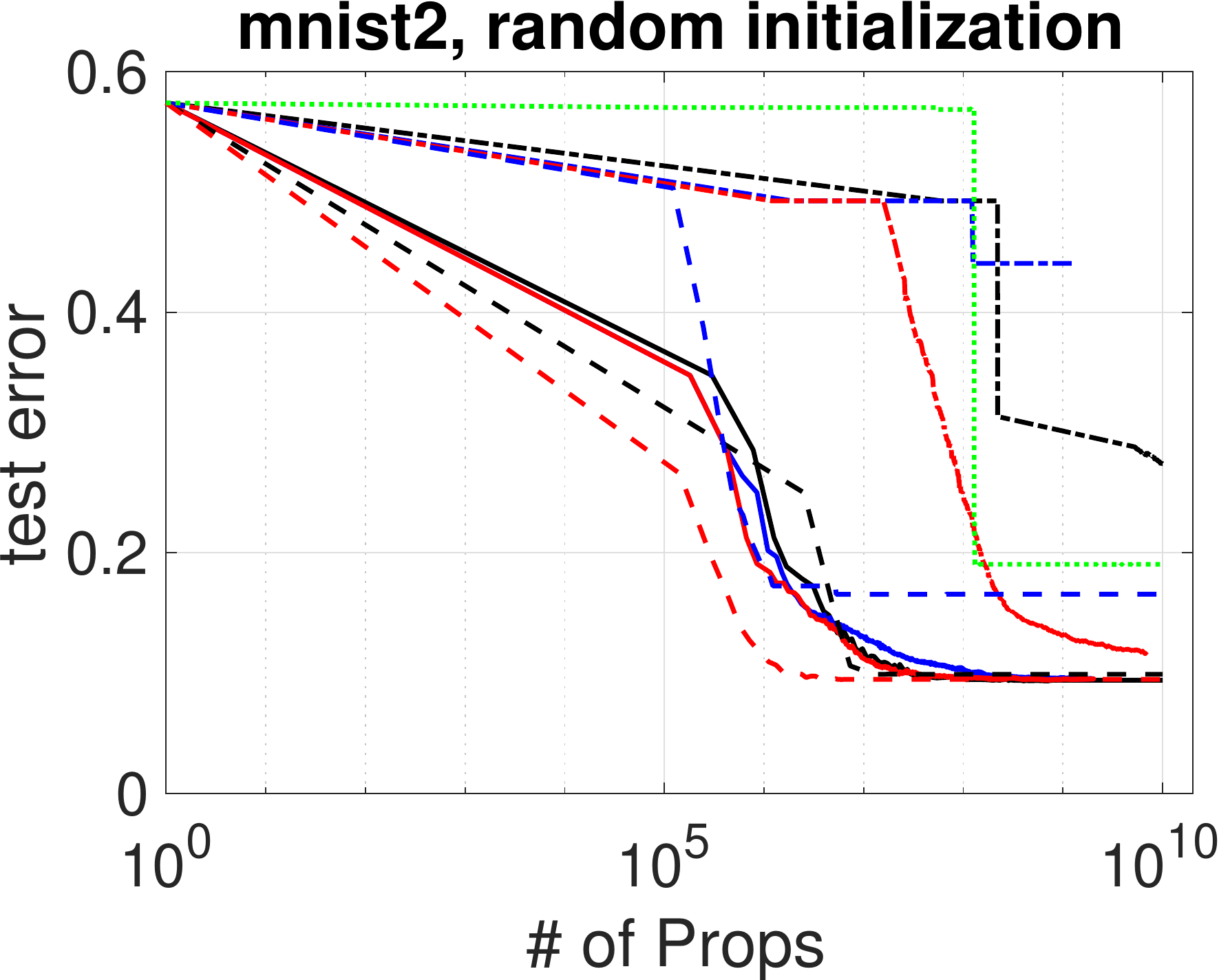}
}
\label{fig:mnist2-1}
\includegraphics[width=0.6\textwidth]{figs/blc_legend_5}\vspace{-2mm}
\subfigure[use $5\%$ of training data]{
\includegraphics[width=0.24\textwidth]{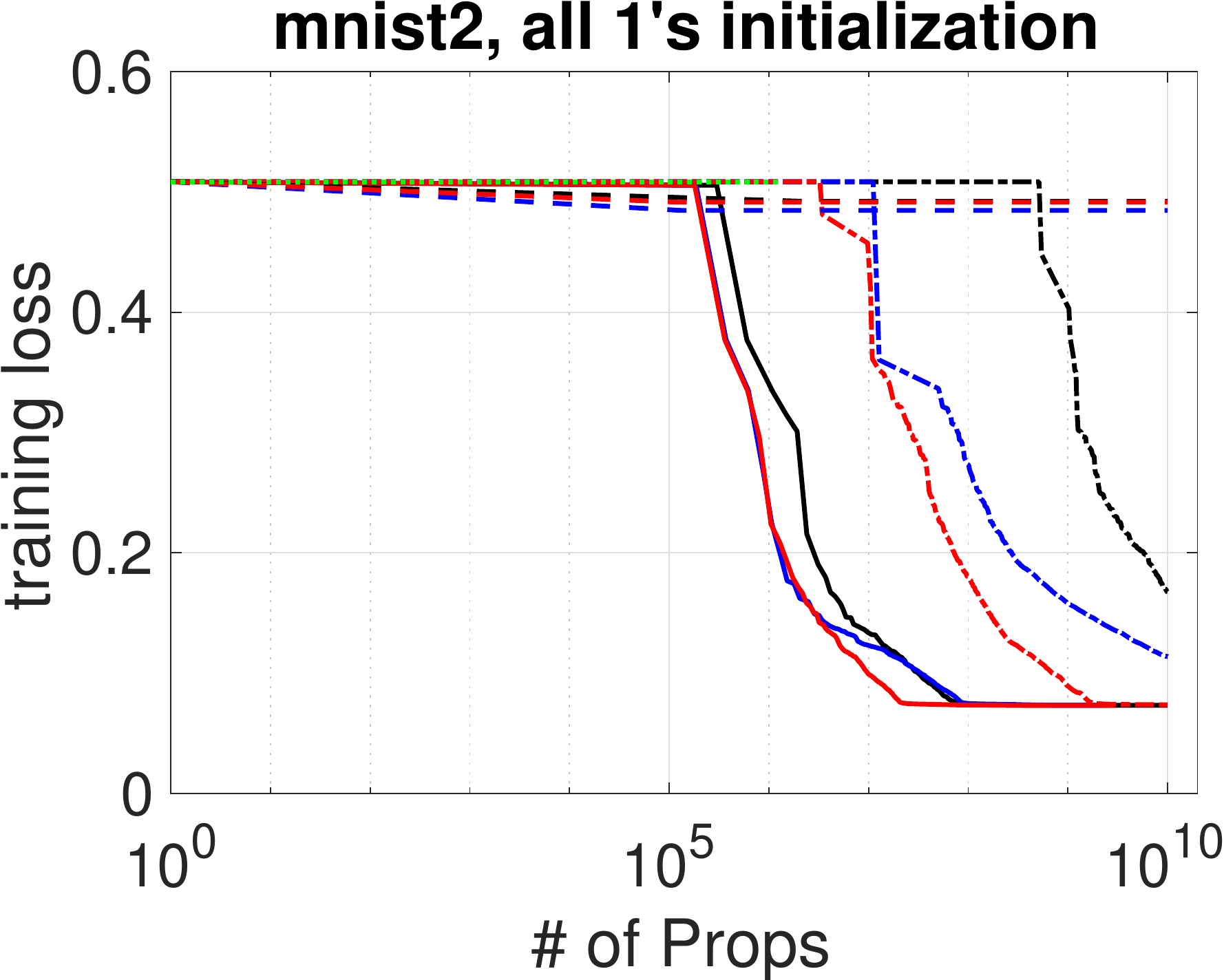}
\includegraphics[width=0.24\textwidth]{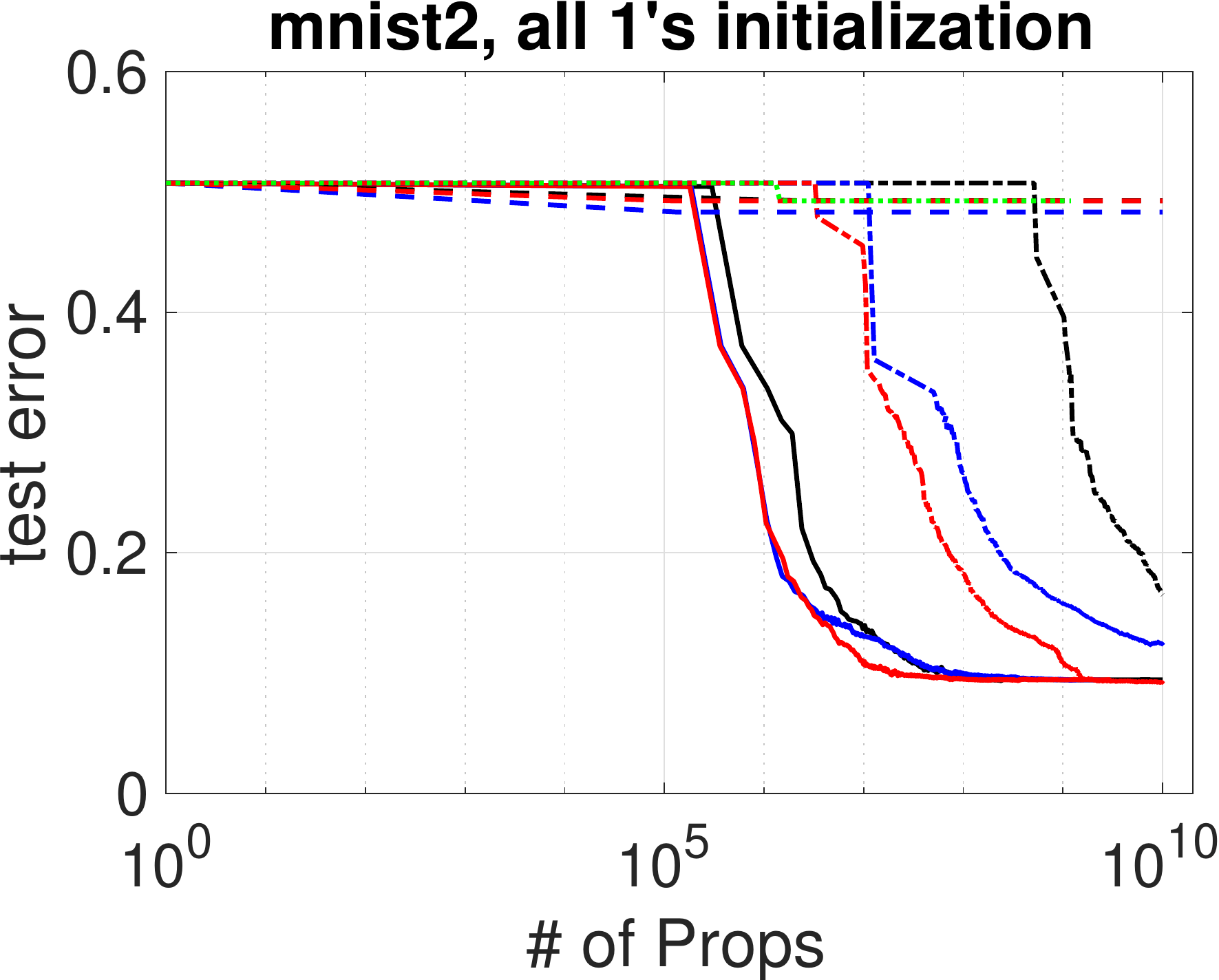}
\includegraphics[width=0.24\textwidth]{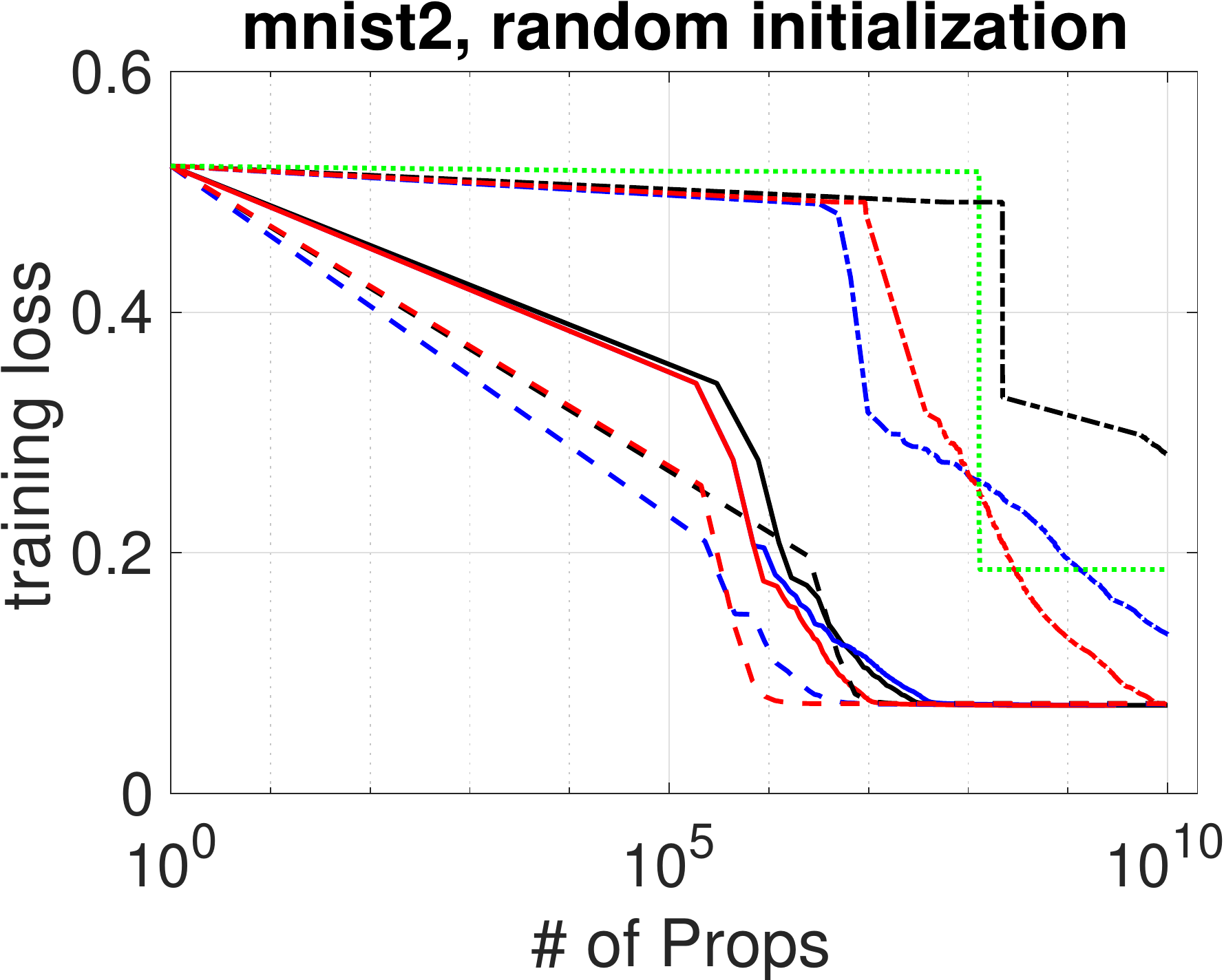}
\includegraphics[width=0.24\textwidth]{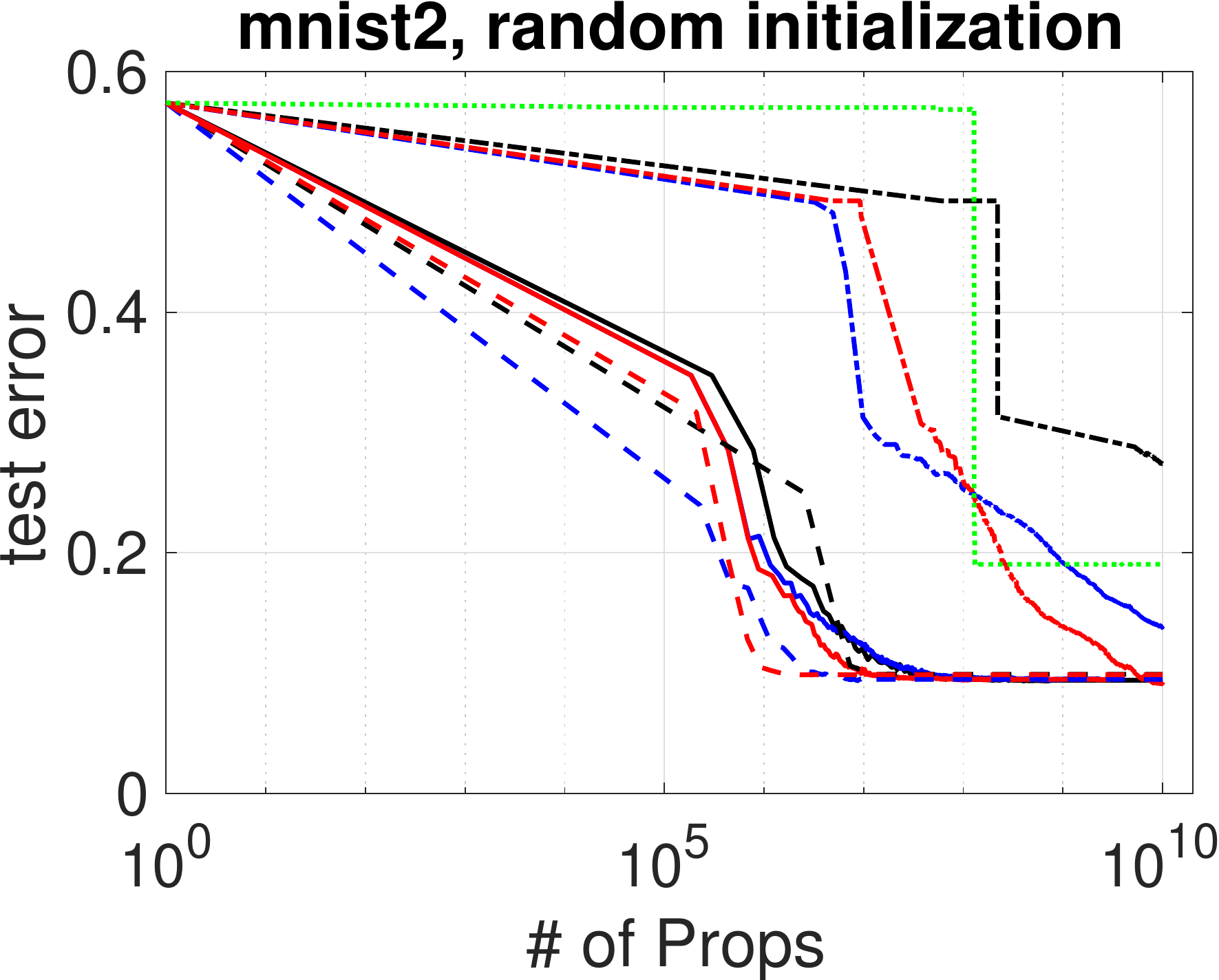}
}
\caption{Binary Classification on \texttt{mnist2}(odd and even digits) with different sub-sampling sizes. ARC Full outperforms ARC Uniform with $ 1\% $ sample size, suggesting that $ 1\% $ might be too small a sample size for uniform sampling. With $ 5\% $ sample size, however, ARC Non-Uniform converges faster than ARC-Full. Increasing the sample size from $1\%$ to $5\%$ does greatly affect the performance of algorithms with uniform sampling, whereas those employing non-uniform sampling appear to be unaffected.}
\label{fig:mnist2-5}
\end{figure}

\begin{figure*}[htbp]
	\centering
	\includegraphics[scale=0.4]{figs/blc_legends.pdf}
	
	\includegraphics[width=0.24\textwidth]{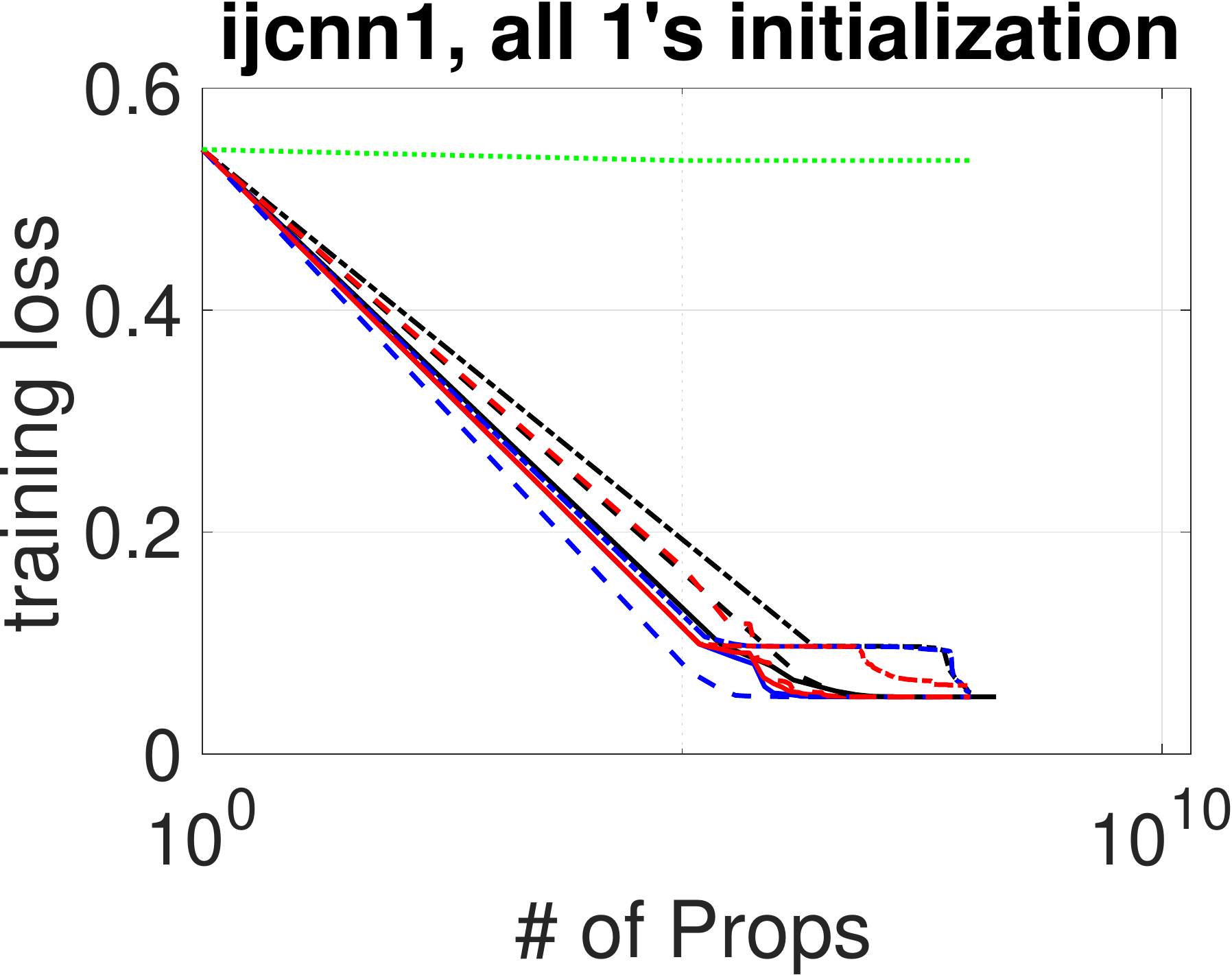}
	\includegraphics[width=0.24\textwidth]{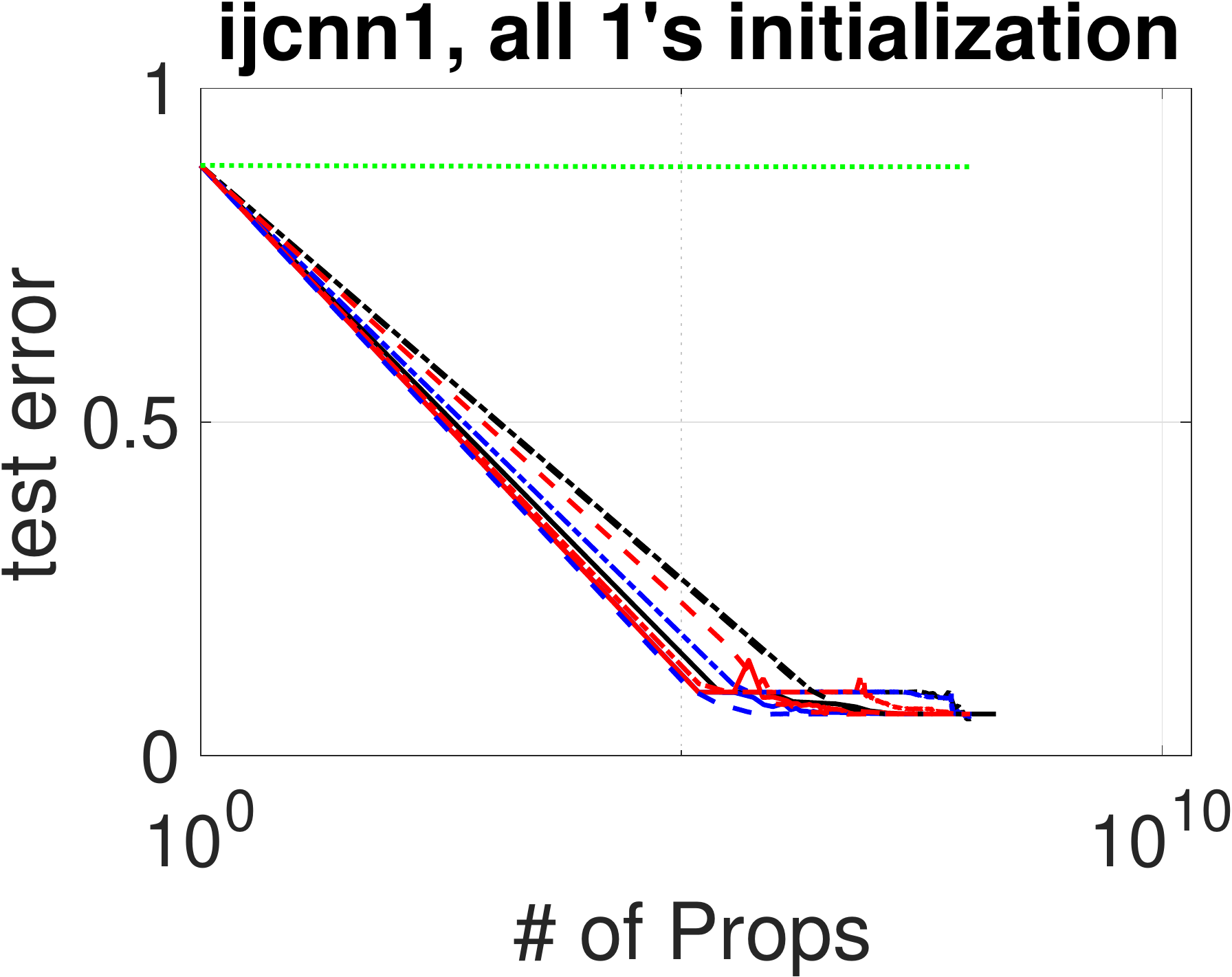}
	\includegraphics[width=0.24\textwidth]{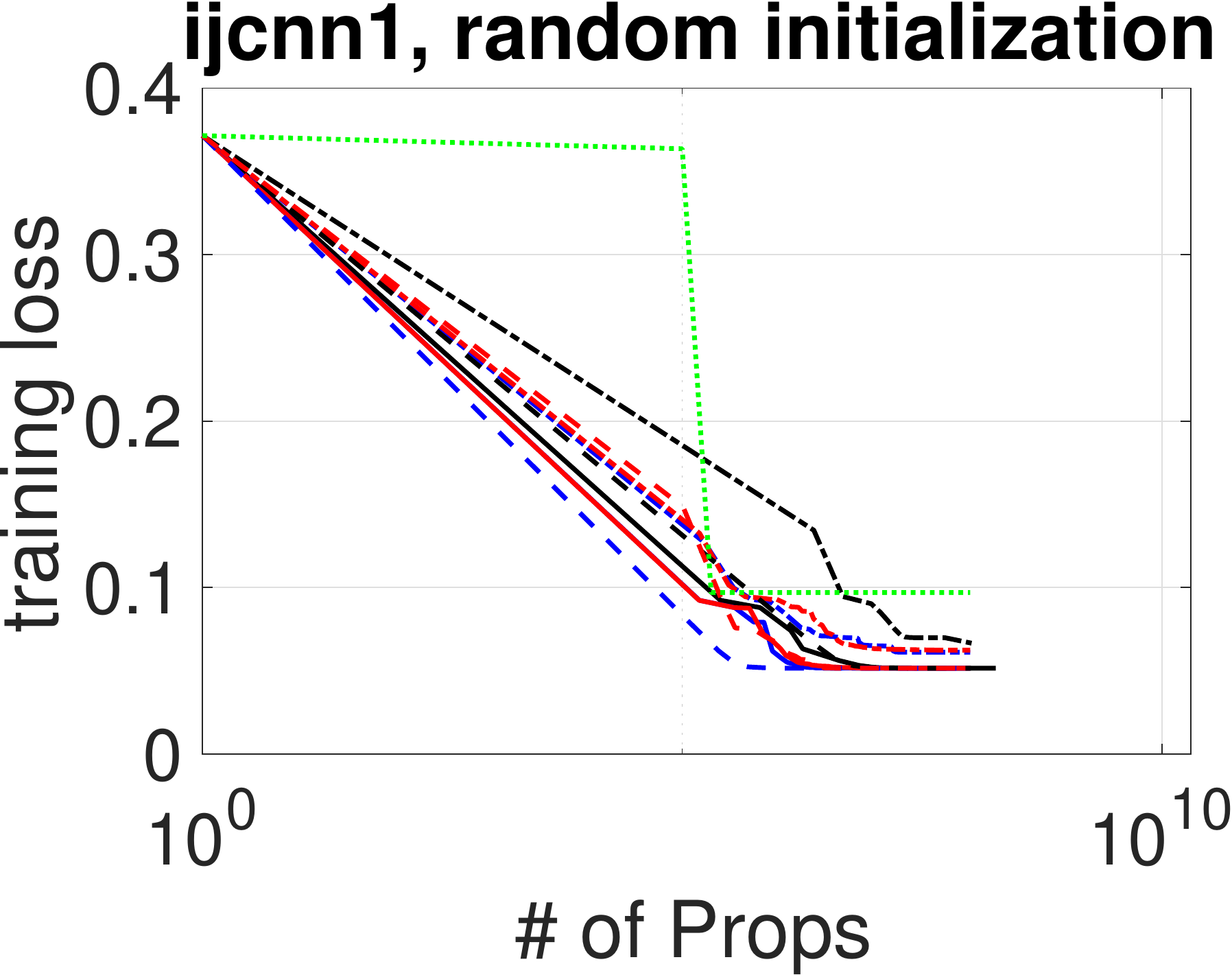}
	\includegraphics[width=0.24\textwidth]{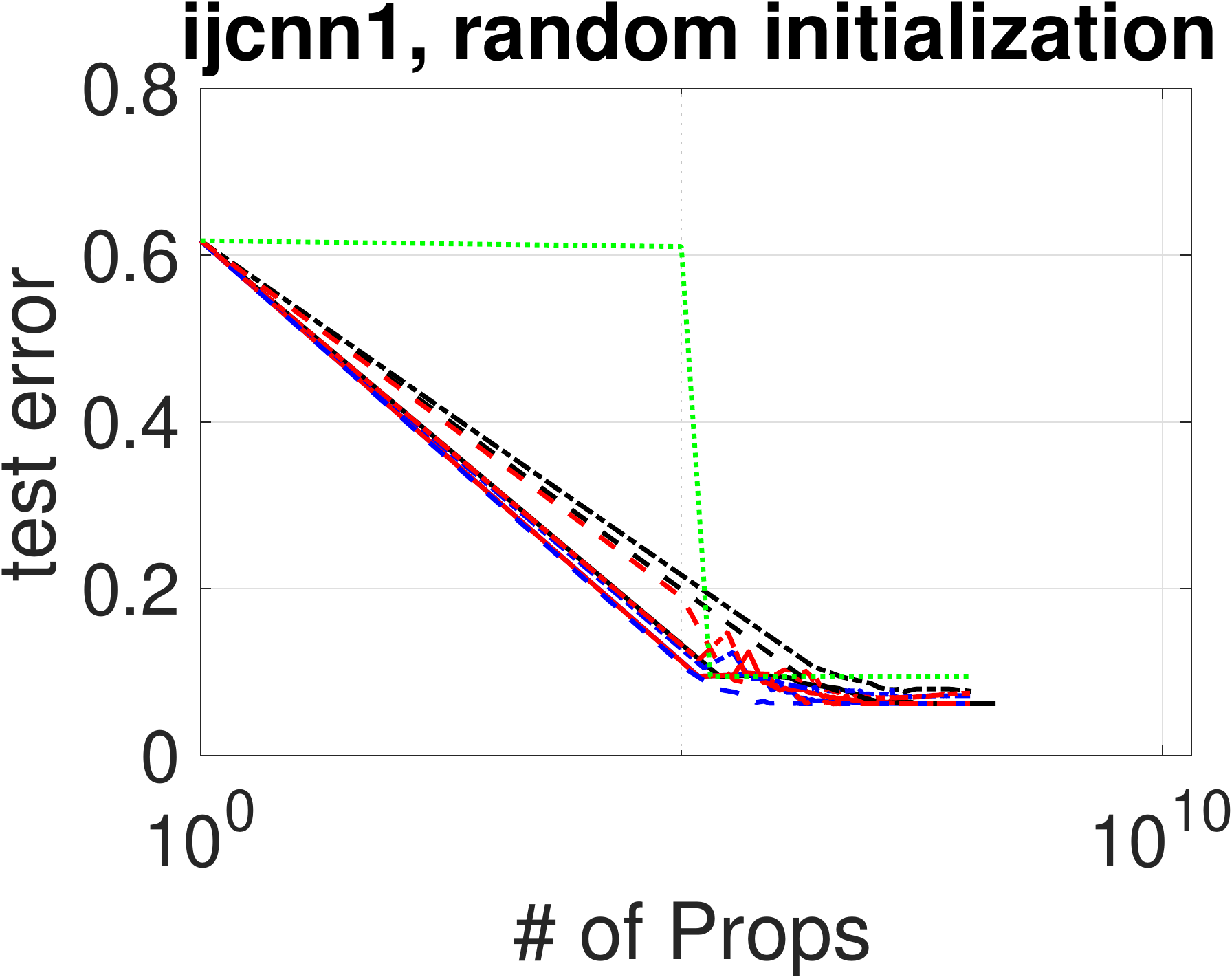}
	
	\caption{Binary Classification on Different \texttt{ijcnn1} with all $1$'s initialization (left two columns) and randon initialization (right two columns). The x-axis is the number of propagations in log-scale.}
	\label{fig:blc_ijcnn1}
\end{figure*}


\end{document}